\documentclass[a4paper,10pt]{article}
\usepackage{amsthm,amsmath,amssymb}
\usepackage{epsfig}
\usepackage{geometry}
\usepackage{listings}

\newtheorem{assumption}{Assumption}
\newtheorem{remark}{Remark}
\newtheorem{theorem}{Theorem}
\newtheorem{lemma}{Lemma}
\newtheorem{prop}{Proposition}

\title{A Milstein scheme for
SPDEs}

\author{Arnulf Jentzen$^1$ 
and 
Michael R\"{o}ckner$^{2}$
\bigskip
\\
\small{$^1$Program in Applied 
and Computational Mathematics, 
Princeton University,}
\\
\small{Princeton, 
NJ 08544-1000, 
USA, 
e-mail:
ajentzen@math.princeton.edu}
\smallskip
\\
\small{$^2$Faculty 
of Mathematics, 
Bielefeld University,
33501 Bielefeld,} 
\\
\small{Germany,
e-mail: roeckner@math.uni-bielefeld.de}
}

\begin{document}

\maketitle

\begin{abstract}
This article studies an infinite 
dimensional analog of Milstein's 
scheme for finite dimensional 
stochastic ordinary differential 
equations (SODEs). 
The Milstein scheme is known to be impressively efficient for SODEs 
which fulfill a certain commutativity
type condition.
This article introduces the infinite
dimensional analog of this
commutativity type condition
and observes that a certain class
of semilinear stochastic partial differential 
equation (SPDEs) with 
multiplicative trace class noise
naturally fulfills the 
resulting infinite
dimensional commutativity condition.
In particular, a suitable 
infinite dimensional
analog of Milstein's algorithm can be
simulated efficiently for such SPDEs
and requires less 
computational operations
and random variables than 
previously considered algorithms 
for simulating such SPDEs.
The analysis is supported by numerical 
results for a stochastic heat equation
and
stochastic reaction diffusion equations
showing significant computational
savings.
\end{abstract}
\section{Introduction}\label{sec:intro}

In this article an infinite 
dimensional analog of 
Milstein's scheme for 
finite dimensional 
stochastic ordinary differential 
equations (SODEs) is studied. 
In order to get a better
understanding of this Milstein type
scheme in infinite dimensions,
we first briefly review Milstein's method
for finite dimensional SODEs
and then concentrate on the 
case of infinite
dimensional stochastic
partial differential equations
(SPDEs) in the
rest of this introductory section.

Let 
$ T \in (0,\infty) $ 
be a real number, 
let
$ 
  d, m \in \mathbb{N} := \{1,2,\ldots\} 
$
be natural numbers,
let 
$ ( \Omega, \mathcal{F}, \mathbb{P} ) 
$ 
be a probability space 
with a normal filtration 
$ 
  ( \mathcal{F}_t )_{ t \in [0,T] } 
$
and let 
$ 
  w = ( w^1, \ldots, w^m ) 
  \colon [0,T] \times \Omega 
  \rightarrow \mathbb{R}^m 
$ 
be 
an $ m $-dimensional standard 
$ 
  ( \mathcal{F}_t )_{ t \in [0,T] } 
$-Brownian motion.
Moreover, let 
$ x_0 \in \mathbb{R}^d $ 
and 
let 
$ 
  \mu = ( \mu_1, \ldots, \mu_d ) 
  \colon \mathbb{R}^d 
  \rightarrow \mathbb{R}^d 
$ 
and 
$ 
  \sigma = ( \sigma_{ i, j } 
  )_{ 
    i \in \{ 1, \ldots, d \}, 
    j \in \{ 1, \ldots, m \} 
  } 
  \colon 
  \mathbb{R}^d \rightarrow 
  \mathbb{R}^{ d \times m } 
$ 
be two smooth functions 
satisfying 
suitable
Lipschitz assumptions
(see
condition~(3.21) 
in Theorem~10.3.5 
in P.\ E.\ Kloeden and 
E.\ Platen~\cite{kp92} 
for details). 
The SODE
\begin{equation}\label{eq:SODE_A}
  d X_t
  =
  \mu( X_t ) \, dt
  +
  \sigma( X_t ) \, dw_t,
  \qquad 
  X_0 = x_0
\end{equation}
for all $ t \in [0,T] $ then admits a unique solution. 
More precisely, there exists an up to 
indistinguishability unique adapted stochastic 
process 
$ 
  X \colon [0,T] \times \Omega 
  \rightarrow \mathbb{R}^d 
$ 
with continuous 
sample paths which satisfies
\begin{align}\label{eq:Cprocess_X}
  X_t
 &=
  x_0
  +
  \int_0^t
  \mu( X_s ) \, ds
  +
  \int_0^t
  \sigma( X_s ) \, dw_s
\\&=
\nonumber
  x_0
  +
  \int_0^t
  \mu( X_s ) \, ds
  +
  \sum_{ i = 1 }^{ m }
  \int_0^t
  \sigma_i( X_s ) \, dw^i_s
\end{align}
$ \mathbb{P} $-a.s.\ for 
all $ t \in [0,T] $. 
Here 
$ 
  \sigma_i \colon \mathbb{R}^d 
  \rightarrow \mathbb{R}^d 
$ 
is given by 
$ 
  \sigma_i( x ) 
  = ( \sigma_{ 1,i }( x ), \ldots, 
  \sigma_{ d,i }( x ) ) 
$ 
for all 
$ x \in \mathbb{R}^d $ 
and all $ i \in \{ 1, \ldots, m \} $. 
Milstein's method 
(see, e.g., (3.3) in Section~10.3 
in P.\ E.\ Kloeden 
and E.\ Platen~\cite{kp92} and 
also G.\ N.\ Milstein's original 
article~\cite{m74}) applied to the SODE~\eqref{eq:SODE_A} 
is then given by 
$ \mathcal{F} $/$ \mathcal{B}( \mathbb{R}^d ) $-measurable 
mappings 
$ 
  y_n^N \colon
  \Omega \rightarrow \mathbb{R}^d 
$, 
$ 
  n \in \{0,1,\ldots,N\} 
$, 
$ 
  N \in \mathbb{N} 
$, 
with 
$ 
  y_0^N = x_0 
$ 
and
\begin{equation}
\label{eq:Milstein}
\begin{split}
  y_{ n+1 }^N
& =
  y_n^N
  +
  \frac{ T }{ N }
  \cdot 
  \mu( y_n^N )
  +
  \sum_{ i = 1 }^{ m }
  \sigma_i( y_n^N )
  \cdot 
  \left( 
    w^i_{ \frac{ (n+1)T }{ N } }
    -
    w^i_{ \frac{ n T }{ N } }
  \right)
\\ & \quad
  +
  \sum_{ i,j = 1 }^{ m }
  \sum_{ k = 1 }^{ d }
  \left( 
  \frac{ \partial }{ \partial x_k } \sigma_i \right)
  \!( y_n^N ) \cdot
  \sigma_{ k,j }( y_n^N )
  \cdot
  \int_{ \frac{ n T }{ N } 
  }^{ \frac{ (n+1)T }{ N } }
  \int_{ \frac{ n T }{ N } }^s
  dw_u^j \, dw_s^i
\end{split}
\end{equation}
$\mathbb{P}$-a.s.\ for 
all $ n \in \{0,1,\ldots,N-1\} $ 
and all $ N \in \mathbb{N} $. 
Although Milstein's scheme is known to 
converge significantly faster than many 
other methods such as the Euler-Maruyama scheme,
it is only of limited use due 
to difficult simulations 
of the iterated stochastic integrals
$
  \int_{ \frac{ n T }{ N } 
  }^{ \frac{ (n+1)T }{ N } }
  \int_{ \frac{ n T }{ N } }^s
  dw_u^j \, dw_s^i
$
for $ i, j \in \{1,\ldots, m\} $ with $ i \neq j $, $ n \in \{0,1,\ldots N-1 \} $ and $ N \in \mathbb{N} $
in~\eqref{eq:Milstein}. 
In the special situation of 
so called {\it commutative 
noise} (see (3.13) in 
Section~10.3 in \cite{kp92}), 
i.e.,
\begin{equation}
\label{eq:comm_noise}
  \sum_{ k = 1 }^{ d }
  \left( \frac{ \partial }{ \partial x_k } \sigma_i \right)
  \!( x ) \cdot
  \sigma_{ k,j }( x )
  =
  \sum_{ k = 1 }^{ d }
  \left( \frac{ \partial }{ \partial x_k } \sigma_j \right)
  \!( x ) \cdot
  \sigma_{ k,i }( x )
\end{equation}
for all $ x \in \mathbb{R}^d $ and all $ i,j \in \{1,\ldots, m\} $, the Milstein scheme can be simplified 
and complicated iterated stochastic integrals in~\eqref{eq:Milstein} 
can be avoided. More precisely, 
in case~\eqref{eq:comm_noise},
Milstein's 
scheme~\eqref{eq:Milstein} 
reduces to
\begin{equation}
\label{eq:commMilstein}
\begin{split}
  y_{ n + 1 }^N
& =
  y_n^N
  +
  \frac{ T }{ N }
  \cdot 
  \mu( y_n^N )
  +
  \sum_{ i = 1 }^{ m }
  \sigma_i( y_n^N )
  \cdot 
  \left( 
    w^i_{ \frac{ (n+1)T }{ N } }
    -
    w^i_{ \frac{ n T }{ N } }
  \right)
\\ & \quad +
  \frac{1}{2}
  \sum_{ i,j = 1 }^{ m }
  \sum_{ k = 1 }^{ d }
  \left( 
    \frac{ \partial }{ \partial x_k } 
    \sigma_i 
  \right)
  \!( y_n^N ) \cdot
  \sigma_{ k,j }( y_n^N )
  \cdot
  \left( 
    w^i_{ \frac{ (n+1)T }{ N } }
    -
    w^i_{ \frac{ n T }{ N } }
  \right)
  \cdot
  \left( 
    w^j_{ \frac{ (n+1)T }{ N } }
    -
    w^j_{ \frac{ n T }{ N } }
  \right)
\\ & \quad -
  \frac{T}{2 N}
  \sum_{ i = 1 }^{ m }
  \sum_{ k = 1 }^{ d }
  \left( 
    \frac{ \partial }{ \partial x_k } \sigma_i 
  \right)
  \!( y_n^N ) \cdot
  \sigma_{ k,i }( y_n^N )
\end{split}
\end{equation}
$\mathbb{P}$-a.s.\ for 
all $ n \in \{ 0,1, \ldots, N-1 \} $ 
and all $ N \in \mathbb{N} $ 
(see (3.16) in Section~10.3 in \cite{kp92}). 
For instance, in the case $ d = m = 1 $, condition~\eqref{eq:comm_noise} is obviously fulfilled 
and the Milstein scheme~\eqref{eq:commMilstein} 
can then be written as
\begin{equation}
\label{eq:Milstein_rewritten}
\begin{split}
  y_{ n+1 }^N
& =
  y_n^N
  +
  \frac{ T }{ N }
  \cdot 
  \mu( y_n^N )
  +
  \sigma( y_n^N )
  \cdot 
  \left( 
    w_{ \frac{ (n+1)T }{ N } }
    -
    w_{ \frac{ n T }{ N } }
  \right)
\\ & \quad
  +
  \frac{1}{2} \cdot
  \sigma'( y_n^N ) \cdot \sigma( y_n^N )
  \cdot
  \left(
    \left( 
      w_{ \frac{ (n+1)T }{ N } }
      -
      w_{ \frac{ n T }{ N } }
    \right)^2
    -
    \frac{T}{N}
  \right)
\end{split}
\end{equation}
$\mathbb{P}$-a.s.\ for 
all $ n \in \{0,1,\ldots,N-1\} $ 
and all $ N \in \mathbb{N} $ 
(see (3.1) in Section~10.3 in \cite{kp92}). 
Of course, 
\eqref{eq:Milstein_rewritten} can be simulated very efficiently. 
However, \eqref{eq:comm_noise} is in 
the case of a 
multidimensional SODE seldom fulfilled 
and even if it is fulfilled,
Milstein's method~\eqref{eq:commMilstein} 
becomes less useful if $ d, m \in \mathbb{N}$ are large.
For example, if $ d = m = 20 $ holds, 
then the middle term
in~\eqref{eq:commMilstein} contains 
$ 20^3 = 8 000 $ summands.
So, 
more than $ 8 000$ additional 
arithmetic operations
are needed to compute $ y^N_{n+1} $ 
from $ y^N_n $ 
for $ n \in \left\{0,1,\dots,N-1\right\}$, 
$N \in \mathbb{N}$ via~\eqref{eq:commMilstein}
in the case $d=m=20$ in general,
which makes
Milstein's scheme less efficient.
This suggests that there is no
hope to expect
that an infinite dimensional analog
of Milstein's method can be simulated 
efficiently
in the case of infinite dimensional
state spaces 
such as 
$ 
  L^2\!\left( (0,1), \mathbb{R}
  \right) 
$
instead of $ \mathbb{R}^d $ and
$ \mathbb{R}^m $ respectively.
One purpose of this article is to 
demonstrate that this is not true
in the case of a suitable class
of semilinear SPDEs with multiplicative 
trace class noise.
We now illustrate this in more detail.

Let $ H = L^2( (0,1), \mathbb{R} ) $ be the $ \mathbb{R} $-Hilbert space of equivalence classes of 
Lebesgue square integrable 
functions from $ (0,1) $ to $ \mathbb{R} $ and
let 
$ 
  f, b \colon (0,1) \times \mathbb{R}   
  \rightarrow \mathbb{R} 
$ 
be two smooth functions 
satisfying appropriate 
Lipschitz assumptions (see~\eqref{eq:diff_term_assumpt} and~\eqref{eq:partial} for details).
As usual we do not 
distinguish between a
Lebesgue square integrable
function
from $(0,1)$ to $ \mathbb{R}$ 
and its equivalence class in $H$.
Moreover, let $ \kappa \in (0,\infty) $
be a real number,
let 
$ 
  \xi \colon
  [0,1] \rightarrow \mathbb{R} 
$ 
with 
$ 
  \xi(0) = \xi(1) = 0 
$ 
be a smooth function 
and let 
$ 
  W \colon [0,T] \times 
  \Omega \rightarrow H 
$ 
be a standard $ Q $-Wiener process 
with respect to $ \mathcal{F}_t $, $ t \in [0,T] $,
where
$ Q \colon H \rightarrow H $
is a trace 
class operator 
(see, for instance, 
Definition~2.1.12 in \cite{pr07}).
It is a classical result (see, e.g., Proposition~2.1.5 in \cite{pr07}) 
that the covariance
operator 
$ Q \colon H \rightarrow H 
$ 
of the
Wiener process 
$ 
  W \colon [0,T] \times \Omega
  \rightarrow H 
$
has an orthonormal basis 
$ g_j \in H $, 
$ j \in \mathbb{N} $,
of eigenfunctions
with summable 
eigenvalues $ \eta_j \in [0,\infty) $,
$ j \in \mathbb{N} $.
In order to have a more concrete
example, we consider the choice
$ g_j(x) = \sqrt{2} \sin( j \pi x ) $
and $ \eta_j = \frac{1}{j^2} $
for all $ x \in (0,1) $
and all $ j \in \mathbb{N} $
in the following and refer
to Section~\ref{sec:settings}
for our general setting
and to Section~\ref{sec:examples}
for further possible examples.
Then we consider the SPDE
\begin{equation}\label{eq:SPDE2}
  dX_t(x)
  =
  \left[
    \kappa   
    \frac{\partial^2}{
      \partial x^2 
    } X_t(x) 
    + f( x, X_t(x) )
  \right] \! dt
  +
  b( x, X_t(x) ) \, dW_t(x)
\end{equation}
with
$ X_t(0) = X_t(1) = 0 $
and $ X_0(x) = \xi(x) $
for $x \in (0,1) $ and $ t \in [0,T] $ on $ H $. Under the assumptions above the SPDE~\eqref{eq:SPDE2} has a unique mild solution. Specifically, there exists an up to indistinguishability unique adapted
stochastic process 
$ 
  X \colon [0,T] \times \Omega   
  \rightarrow H 
$ 
with continuous 
sample paths which satisfies
\begin{equation}\label{eq:Dprocess_X}
  X_t
  =
  e^{ A t } \, \xi
  +
  \int_0^t
  e^{ A( t - s ) } F( X_s ) \, ds
  +
  \int_0^t
  e^{ A( t - s ) } B( X_s ) \, dW_s
\end{equation}
$ \mathbb{P} $-a.s.\ for 
all $ t \in [0,T] $, 
where 
$ 
  A \colon D( A ) \subset H \rightarrow H 
$ 
is the Laplacian with Dirichlet boundary conditions times
the constant 
$ \kappa \in (0,\infty) $
and 
where 
$ 
  F \colon H \rightarrow H 
$ 
and 
$ 
  B \colon H \rightarrow HS(U_0,H) 
$ 
are given by $ ( F( v ) )( x ) = f( x, v( x ) ) $ and $ ( B( v ) u )( x ) = b( x, v( x ) ) \cdot u( x ) $ for all $ x \in (0,1) $, $ v \in H $ and all $ u \in U_0 $. 
Here 
$ 
  U_0 = Q^{ 1/2 }( H ) 
$
with
$ 
  \left< v, w \right>_{ U_0 } 
  = 
  \left< 
    Q^{ - 1 / 2 } v, 
    Q^{ - 1 / 2 } w 
  \right>_H 
$
for all
$ v, w \in U_0 $
is 
the image 
$ \mathbb{R} $-Hilbert space 
of 
$ 
  Q^{ 1 / 2 } 
$ 
(see Appendix C in \cite{pr07}).
(Note that $ A $ und $ Q $
commutate and even more
satisfy
$ - \kappa \pi^2 \! A^{-1} 
= Q $ in our example 
SPDE~\eqref{eq:Dprocess_X}
in this introductory section
but our general 
setting below does not require that 
these conditions
are fulfilled; 
see Section~\ref{sec:settings}
and Subsection~\ref{sec:reacdiff2}.)

Then our goal is to solve 
the strong approximation 
problem (see 
Section~9.3 in \cite{kp92})
of the SPDE~\eqref{eq:SPDE2}. 
More precisely, we want to 
compute an 
$ \mathcal{F} 
$/$ 
  \mathcal{B}( H ) 
$-measurable numerical 
approximation 
$ 
  Y \colon \Omega \rightarrow H 
$ 
such that
\begin{equation}\label{eq:strongpr}
  \left(
    \mathbb{E}\!\left[
      \int_0^1
      \left| X_T( x ) - Y( x ) \right|^2 dx
    \right]
  \right)^{ \! 1/2 }
  <
  \varepsilon
\end{equation}
for a given precision 
$ \varepsilon > 0 $ with the 
least possible computational effort (number of computational operations 
and independent 
standard normal random variables 
needed to compute 
$ 
  Y \colon \Omega \rightarrow H 
$).
A computational operation is here an arithmetic operation 
(addition, subtraction, multiplication, division),
a trigonometric operation (sine, cosine)
or an evaluation of 
$ 
  f \colon (0,1) \times 
  \mathbb{R} \rightarrow \mathbb{R} 
$,
$ 
  b \colon (0,1) \times 
  \mathbb{R} \rightarrow \mathbb{R} 
$
or the exponential function.

In order to be able to simulate such a numerical approximation on a computer both the time interval $ [0,T] $ and the infinite dimensional space $ H = L^2( (0,1), \mathbb{R} ) $ have to be discretized. While for temporal discretizations the linear implicit Euler scheme (see
\cite{d11,dp09,dz02,g99,h02,h03a,h07a,s99,w05b,y05}) and the linear implicit Crank-Nicolson scheme 
(see \cite{h02,h03a,s99,w05b})
are often used,
spatial discretizations are usually achieved with finite elements (see
\cite{anz98,btz04,dp09,dz02,gkl09,h07a,kkl07,kz09,kll10,kls10,lt10,w05b,y05}), 
finite differences 
(see \cite{g98,gm06,mm05,ps05,r02b,r06b,r07b,s99,w06}) 
and spectral Galerkin methods 
(see 
\cite{gk96,h03a,j08b,jk09b,ks01,lr04,ls06,mr07b,mrw07,mrw08}). 
For instance, 
the linear implicit Euler scheme 
combined with spectral Galerkin 
methods which we denote by 
$ \mathcal{F} $/$ \mathcal{B}( H ) 
$-measurable mappings 
$ 
  Z_n^N \colon \Omega \rightarrow H 
$, 
$ n \in \{ 0,1,\ldots, N^3 \} 
$, 
$ N \in \mathbb{N} $, 
is given by
$ 
  Z^N_0 = P_N( \xi ) 
$ 
and
\begin{equation}\label{eq:implEuler}
  Z_{ n+1 }^N 
  =
  P_N \!
  \left( I - \frac{T}{N^3} A \right)^{ \!\!-1 } \!\!
  \left(
    Z_n^N
    +
    \frac{ T }{ N^3 } \cdot
    f( \cdot, Z_n^N )
    +
    b( \cdot, Z_n^N ) \cdot
    \left( 
      W_{ \frac{ (n+1)T }{ N^3 } }^N
      -
      W_{ \frac{ n T }{ N^3 } }^N
    \right)
  \right)
\end{equation}
$ \mathbb{P} $-a.s.\ for 
all $ n \in \{ 0,1,\ldots, N^3 - 1 \} $ 
and all 
$ N \in \mathbb{N} $. 
Here the bounded linear operators 
$ 
  P_N \colon H \rightarrow H 
$, 
$ N \in \mathbb{N} $, 
and the Wiener processes 
$ 
  W^N \colon 
  [0,T] \times \Omega \rightarrow H 
$, 
$ N \in \mathbb{N} $, 
are given by
\begin{equation}
  ( P_N( v ) )( x )
  =
  \sum_{ n = 1 }^N
  2 \sin( n \pi x )
  \int_0^1
  \sin( n \pi y ) \,
  v( y ) \, dy
\end{equation}
for all $ x \in (0,1) $,
$ v \in H $, $N \in \mathbb{N}$
and by
$
  W_t^N( \omega )
  =
  P_N\!\left(
    W_t( \omega )
  \right)
$
for all 
$ t \in [0,T] $, 
$ \omega \in \Omega $,
$ N \in \mathbb{N} $. 
Moreover, 
the notations 
$ 
  v \cdot w 
  \colon 
  (0,1) \rightarrow \mathbb{R} 
$,
$ 
  v^2 \colon 
  (0,1) \rightarrow \mathbb{R} 
$ 
and
$ 
  \varphi( \cdot, v ) \colon 
  (0,1) \rightarrow \mathbb{R} 
$
given by
\begin{equation}
  \big( v \cdot w \big)(x)
  = v(x) \cdot w(x),
\qquad
  \big( v^2 \big)(x)
  = \big( v(x) \big)^2,
\qquad
  \big( 
    \varphi( \cdot, v ) 
  \big)(x)
  = \varphi( x , v(x) )
\end{equation}
for all $ x \in (0,1)$ and all functions 
$ 
  v, w \colon (0,1) \rightarrow \mathbb{R} 
$,
$ 
  \varphi \colon (0,1) \times \mathbb{R} 
  \rightarrow \mathbb{R} 
$
are used
here and below.
In~\eqref{eq:implEuler}
the infinite 
dimensional $ \mathbb{R} $-Hilbert space $ H $ is projected down to the $ N $-dimensional $ \mathbb{R} $-Hilbert space $ P_N( H ) $ with $ N \in \mathbb{N} $ and the infinite dimensional Wiener process 
$ 
  W \colon [0,T] \times \Omega 
  \rightarrow H 
$ 
is approximated by the finite 
dimensional processes 
$ 
  W^N \colon [0,T] \times 
  \Omega \rightarrow H 
$, 
$ 
  N \in \mathbb{N} 
$, 
for the spatial discretization. 
For the temporal discretization in 
the scheme $ Z_n^N $, $ n \in \{0,1,\ldots,N^3 \} $, 
above
the time interval $ [0,T] $ is divided into $ N^3 $ subintervals, i.e., 
$ 
  N^3 
$ 
time steps are used, 
for $ N \in \mathbb{N} $. 
The exact solution 
$ 
  X \colon [0,T] \times \Omega 
  \rightarrow H 
$ 
of the SPDE~\eqref{eq:SPDE2} 
has values in 
$ 
  D( (-A)^{ \gamma } ) 
$ 
and satisfies 
$ 
  \mathbb{E}\big[
    \| (-A)^{ \gamma } X_T \|^2_H 
  \big]
  < \infty 
$ 
for all 
$ 
  \gamma \in (0,\frac{3}{4} ) 
$ 
(see Section~4.3 in \cite{jr11}). 
This shows
\begin{align*}
  \left(
  \mathbb{E}\!\left[ \left\|
    X_T - P_N( X_T )
  \right\|_H^2
  \right]
  \right)^{ 1 / 2 }
 &\leq 
  \left(
    \mathbb{E}\!\left[
    \left\|
      (-A)^{ \gamma } X_T
    \right\|^2_H
    \right]
  \right)^{ 1 / 2 }
  \left\| 
    (-A)^{ -\gamma }( I - P_N ) 
  \right\|_{ L( H ) }
\\&\leq 
  \left(
    \mathbb{E}\!\left[
    \left\|
      (-A)^{ \gamma } X_T
    \right\|_H^2
    \right]
  \right)^{ 1 / 2 }
  \left( 1 + \kappa^{-1}
  \right)
  N^{ - 2 \gamma }
< \infty
\end{align*}
for all $ N \in \mathbb{N} $ 
and all $ \gamma \in (0,\frac{3}{4} ) $. 
So, $ P_N( X_T ) $ converges in the root
mean square sense to $ X_T $ with 
order $ \frac{3}{2}- $ as $ N $ goes to 
infinity. 
(For a real number 
$ \delta \in (0,\infty) $,
we write $ \delta - $ for the 
convergence order
if the convergence order is higher than
$ \delta - r $ for all arbitrarily small
$ r \in (0,\delta) $.)
Additionally, 
the solution process
of the SPDE~\eqref{eq:SPDE2}
is known to be
$ \frac{ 1 }{ 2 } $-H\"{o}lder
continuous in the root 
mean square sense
(see, e.g.,
Theorem~1 
in \cite{jr11})
and therefore,
the linear implicit Euler scheme
converges temporally 
in the root mean square
to the exact solution
of the SPDE~\eqref{eq:SPDE2}
with order $ \frac{1}{2} $ 
(see, e.g., Theorem~1.1 in \cite{y05}). 
Combining the convergence 
rate $ \frac{3}{2}- $ for the spatial 
discretization and the convergence 
rate $ \frac{1}{2} $ for the temporal 
discretization indicates that it is 
asymptotically optimal to use 
the cubic number $ N^3 $ 
of time steps 
in the linear implicit Euler 
scheme $ Z_n^N $, $ n \in \{ 0,1,\ldots, N^3 \} $, above.

We now review how efficiently 
the numerical method~\eqref{eq:implEuler} solves the strong 
approximation problem~\eqref{eq:strongpr} 
of the SPDE~\eqref{eq:SPDE2}. 
Standard results in the 
literature (see, for instance, 
Theorem~2.1 in \cite{h03a}) yield 
the existence of real 
numbers $ C_r > 0 $,
$ r \in (0,\frac{3}{2}) $, such that
\begin{equation}\label{eq:ArealConstant}
  \left(
    \mathbb{E}\!\left[
      \int_0^1
      \left| 
        X_T( x ) - 
        Z_{ N^3 }^N( x ) 
      \right|^2 dx
    \right]
  \right)^{ \! 1 / 2 }
  \leq 
  C_r \cdot 
  N^{ ( r - \frac{3}{2} ) }
\end{equation}
for all $ N \in \mathbb{N} $
and all arbitrarily small 
$ r \in (0,\frac{3}{2})$. 
The linear implicit
Euler approximation
$ Z^N_{ N^3 } $ thus converges
in the root mean square sense
to $ X_T $ with 
order $ \frac{3}{2}- $
as $N$ goes to infinity.
Moreover, since $ P_N( H ) $ is $ N $-dimensional and since $ N^3 $ 
time steps are used in~\eqref{eq:implEuler}, 
$ O( N^4 \log(N) ) $ computational operations 
and random variables are needed 
to compute $ Z_{ N^3 }^N $. 
The logarithmic 
term in $ O( N^4 \log(N) ) $  
arises due to computing the nonlinearities $ f $ and $ b $
with fast Fourier transform where aliasing errors are neglected here and below. Combining the computational 
effort $ O( N^4 \log(N) ) $ and 
the convergence order $ \frac{3}{2}- $ 
in \eqref{eq:ArealConstant} shows 
that the linear implicit Euler 
scheme needs about 
$ 
  O(\varepsilon^{ -\frac{8}{3} })
$ 
computational operations and 
independent standard normal 
random variables
to achieve the desired 
precision $ \varepsilon > 0 $ 
in \eqref{eq:strongpr}.
In fact, we have demonstrated 
that Euler's 
method~\eqref{eq:implEuler} needs 
$ O( \varepsilon^{ -(\frac{8}{3} + r) }) $
computational operations and 
independent standard normal
random variables
to solve~\eqref{eq:strongpr} for all arbitrarily 
small $ r \in (0, \infty ) $. 
However, we write about 
$ O(\varepsilon^{ - \frac{8}{3} }) $
computational operations 
and independent standard normal
random variables 
for simplicity here and below.

Having reviewed Euler's 
method~\eqref{eq:implEuler}, 
we now derive and study
an infinite 
dimensional analog of Milstein's scheme. 
In the finite dimensional SODE case, 
the Milstein scheme~\eqref{eq:Milstein} is 
derived by applying It{\^o}'s 
formula to the integrand process 
$ \sigma( X_t ) $, $ t \in [0,T] $, 
in~\eqref{eq:Cprocess_X}. 
This approach is based on the fact that the 
diffusion coefficient $ \sigma $ 
is a smooth test function and 
that the solution process of~\eqref{eq:SODE_A} 
is an It{\^o} process. This strategy is 
not directly available in infinite
dimensions since~\eqref{eq:SPDE2} does in 
general not admit a strong solution 
to which 
the standard It{\^o} formula 
could be applied.
Recently, in \cite{jk09c} in the case of
additive noise and in 
\cite{j10} in the 
general case, this problem has been 
overcome by first applying Taylor's 
formula in Banach spaces to the 
diffusion coefficient $ B $ in the 
mild integral equation~\eqref{eq:Dprocess_X} 
and by then inserting a lower order 
approximation recursively 
(see Section~4.3 in \cite{j10}).
More formally, using
$ F( X_s ) \approx F( X_0 ) $ and 
$ B( X_s ) \approx B( X_0 ) + B'( X_0 )( X_s - X_0 ) $ 
for small $ s \in [0,T] $ 
in~\eqref{eq:Dprocess_X}
shows
\begin{align*}
  X_t
 &\approx
  e^{ A t } \xi
  +
  \int_0^t
  e^{ A( t - s ) }
  F( X_0 ) \, ds
  +
  \int_0^t
  e^{ A( t - s ) }
  B( X_0 ) \, dW_s
  +
  \int_0^t
  e^{ A( t - s ) }
  B'( X_0 ) ( X_s - X_0 ) \, dW_s
\\&\approx
  e^{ A t }
  \left(
    X_0 + t \cdot F( X_0 )
    +
    \int_0^t
    B( X_0 ) \, dW_s
    +
    \int_0^t
    B'( X_0 ) ( X_s - X_0 ) \, dW_s
  \right)
\end{align*}
for small $ t \in [0,T] $. 
(We would like to remark that $ B $ 
is, in general, not Fr\'{e}chet 
differentiable on $ H $ but on
a suitable dense subspace of $ H $; 
see Assumption~4 below 
for details.
In this
introductory section
we simply write 
$ B' $ for this Fr\'{e}chet derivative
on a suitable dense subspace
of $ H $ 
and refer to 
Section~\ref{sec:settings} below
for the precise
handling of this issue.)
The estimate 
$ X_s \approx X_0 + 
\int_0^s B( X_0 ) \, dW_u $
for small $ s \in [0,T] $ then gives
\begin{equation}\label{eq:tempApprox}
  X_t
  \approx
  e^{ A t } \!
  \left(
    X_0 + t \cdot F( X_0 )
    +
    \int_0^t B( X_0 ) \, dW_s
    +
    \int_0^t \!
    B'( X_0 ) 
    \left( 
      \int_0^s
      B( X_0 ) \, dW_u \!
    \right) dW_s \!
  \right)
\end{equation}
for small $ t \in [0,T] $.
Using It{\^o}'s formula
this temporal approximation
has already been obtained
in (1.12) in \cite{ms06}
under additional smoothness
assumptions of the driving noise
process of the 
SPDE~\eqref{eq:Dprocess_X}
(see Assumption C in \cite{ms06})
which guarantee the existence of
a strong solution and thus allow
the application of 
the standard It{\^o} formula.
Combining the temporal 
approximation~\eqref{eq:tempApprox} 
and the spatial 
discretization in~\eqref{eq:implEuler} 
suggests the numerical scheme with
$ \mathcal{F} $/$ \mathcal{B}( H ) $-measurable
mappings 
$ 
  Y_n^N \colon 
  \Omega \rightarrow H 
$,
$ n \in \{0,1,\ldots,N^2\} $,
$ N \in \mathbb{N} $, given by
$ Y_0^N = P_N( \xi ) $ and
\begin{equation}
\label{eq:numScheme}
\begin{split}
  Y_{ n + 1 }^N
& =
  P_N \, e^{ A \frac{ T }{ N^2 } }
  \Bigg(
    Y_n^N
    +
    \frac{ T }{ N^2 } \cdot F( Y_n^N )
    +
    B( Y_n^N )
    \left( 
      W_{ \frac{ (n+1)T }{ N^2 } }^N
      -
      W_{ \frac{ n T }{ N^2 } }^N
    \right)
\\ & \quad
    +
    \int_{ \frac{ n T }{ N^2 } }^{ \frac{ (n+1)T }{ N^2 } }
    B'( Y_n^N )
    \left(
      \int_{ \frac{ n T }{ N^2 } }^s
      B( Y_n^N ) \, dW_u^N
    \right) dW_s^N
  \Bigg)
\end{split}
\end{equation}
$ \mathbb{P} $-a.s.\ for 
all $ n \in \{0,1,\ldots,N^2-1\} $ 
and all $ N \in \mathbb{N} $. 
Now, we are at a stage similar to 
the finite dimensional 
case~\eqref{eq:Milstein}: a higher
order Milstein
type method seems to be derived which nevertheless 
seems to be of 
limited use due to the iterated
high dimensional stochastic
integral in~\eqref{eq:numScheme}.
However, a key observation here
is the formula
\begin{multline}
\label{eq:keyObservation}
  \int_{ \frac{ n T }{ N^2 } 
  }^{ \frac{ (n+1)T }{ N^2 } }
  B'( Y_n^N )
  \left(
    \int_{ \frac{ n T }{ N^2 } }^s
    B( Y_n^N ) \, dW_u^N
  \right) dW_s^N
\\
  =
  \frac{1}{2} 
  \left( \frac{ \partial }{ \partial y }b 
  \right)\!( \cdot, Y_n^N ) 
  \cdot b( \cdot, Y_n^N )
  \cdot
  \left(
    \left( 
      W_{ \frac{ (n+1)T }{ N^2 } }^N
      -
      W_{ \frac{ n T }{ N^2 } }^N
    \right)^{  \! 2 }
    -
    \frac{T}{N^2}
    \sum_{ i = 1 }^{ N }
    \eta_i ( g_i )^2
  \right)
\end{multline}
$ \mathbb{P} $-a.s.\ for 
all $ n \in \{0,1,\ldots,N^2-1\} $ 
and all $ N \in \mathbb{N} $ (see 
Subsection~\ref{sec:centralob} 
for the proof of the iterated
integral identity~\eqref{eq:keyObservation} and
see below for a heuristic explanation of this fact).
So, the iterated high dimensional stochastic
integral in~\eqref{eq:numScheme} reduces to a
simple product of functions. 
The function
$ 
  \frac{\partial}{\partial y } b 
  \colon 
  (0,1)
  \times
  \mathbb{R} \rightarrow \mathbb{R} 
$ 
is here the partial
derivative $ ( \frac{\partial}{\partial y} b )( x, y ) $
for $ x \in (0,1) $ and $ y \in \mathbb{R} $.
Using~\eqref{eq:keyObservation} the numerical scheme~\eqref{eq:numScheme} thus reduces to
\begin{equation}
\label{eq:reducedScheme}
\begin{split}
  Y_{ n+1 }^{ N }
& =
  P_N \, e^{ A \frac{ T }{ N^2 } }
  \Bigg( 
    Y_n^N 
    + 
    \frac{ T }{ N^2 } \cdot
    f( \cdot, Y_n^N )
    +
    b( \cdot, Y_n^N ) \cdot
    \left( 
      W_{ \frac{ (n+1)T }{ N^2 } }^N
      -
      W_{ \frac{ n T }{ N^2 } }^N
    \right)
\\ & \quad
    +
    \frac{1}{2} 
    \left( \frac{ \partial }{ \partial y }b 
    \right)\!( \cdot, Y_n^N ) 
    \cdot b( \cdot, Y_n^N ) \cdot
    \left(
      \left( 
        W_{ \frac{ (n+1)T }{ N^2 } }^N
        -
        W_{ \frac{ n T }{ N^2 } }^N
      \right)^{ \! 2 }
      - 
      \frac{ T }{ N^2 }
      \sum_{ i = 1 }^{ N }
      \eta_i ( g_i )^2
    \right)
  \Bigg)
\end{split}
\end{equation}
$ \mathbb{P} $-a.s.\ for 
all $ n \in \{ 0,1, \ldots, N^2 - 1 \} $ and all $ N \in \mathbb{N} $. 
Note 
that only increments of the finite
dimensional Wiener processes 
$ 
  W^N \colon [0,T] \times \Omega 
  \rightarrow H 
$, 
$ 
  N \in \mathbb{N} 
$, 
are used in \eqref{eq:reducedScheme}.
Moreover, observe that, as in the case of~\eqref{eq:implEuler}, 
the infinite dimensional 
$ \mathbb{R} $-Hilbert space 
$ H $ is projected down to the 
$ N $-dimensional 
$ \mathbb{R} 
$-Hilbert space $ P_N( H ) $ with 
$ N \in \mathbb{N} $ and the 
infinite dimensional Wiener process 
$ 
  W \colon [0,T] \times \Omega 
  \rightarrow H 
$ 
is approximated by the 
finite dimensional Wiener processes 
$ 
  W^N \colon [0,T] \times 
  \Omega \rightarrow H 
$, 
$ 
  N \in \mathbb{N} 
$, 
for the 
spatial discretization
in \eqref{eq:reducedScheme}. For the temporal discretization in the scheme $ Y_n^N $, $ n \in \{ 0,1, \ldots, N^2 \} $,
above the time 
interval $ [0,T] $ is divided into $ N^2 $ subintervals, i.e., 
$ N^2 $ 
instead of $ N^3 $ time steps 
are used in~\eqref{eq:reducedScheme}, 
for $ N \in \mathbb{N} $.
In the following we explain why it is 
crucial to use $ N^2 $ time steps 
in \eqref{eq:reducedScheme} instead of
$ N^3 $ time steps in the case of the
linear implicit Euler scheme~\eqref{eq:implEuler}.

More formally, we now illustrate 
how efficiently the Milstein
type algorithm~\eqref{eq:reducedScheme}
solves the strong approximation problem~\eqref{eq:strongpr} of the SPDE~\eqref{eq:SPDE2}.
Theorem~\ref{thm:mainresult} (see Section~\ref{secscheme} below) 
gives the existence of real 
numbers
$ C_r > 0 $, 
$ r \in (0, \frac{3}{2} ) $,
such that
\begin{equation}
\label{eq:realNumberCepsilon}
  \left(
    \mathbb{E}\!\left[
      \int_0^1 
      \left|
        X_T( x ) - Y_{ N^2 }^N( x )
      \right|^2 dx
    \right]
  \right)^{ \! 1 / 2 }
  \leq 
  C_r 
  \cdot
  N^{ ( r - \frac{3}{2} ) }
\end{equation}
for all 
$ N \in \mathbb{N} $ 
and all arbitrarily small 
$ r \in (0,\frac{3}{2} ) $. 
The approximation
$ Y^N_{ N^2 } $ 
thus converges
in the root mean square sense
to $ X_T $ with 
order $ \frac{3}{2}- $
as $N$ goes to infinity.
The expression
\begin{equation}\label{eq:additionalInfo}
  \frac{1}{2} 
  \left( \frac{ \partial }{ \partial y }b 
  \right)\!( \cdot, Y_n^N ) 
  \cdot b( \cdot, Y_n^N )
  \cdot
  \left(
    \left(
      W_{ \frac{ (n+1)T }{ N^2 } }^N
      -
      W_{ \frac{ n T }{ N^2 } }^N
    \right)^{ \! 2 }
    -
    \frac{ T }{ N^2 }
    \sum_{ i = 1 }^{ N }
    \eta_i ( g_i )^2
  \right)
\end{equation}
for $ n \in \{ 0,1,\ldots, N^2 - 1\} $ and $ N \in \mathbb{N} $ 
in the Milstein type approximation
\eqref{eq:reducedScheme} 
contains additional information of
the solution process 
of \eqref{eq:SPDE2} and this
allows us to use less time steps,
$ N^2 $ in \eqref{eq:reducedScheme} 
instead of $ N^3 $
in \eqref{eq:implEuler}, to achieve the same
convergence rate as the linear implicit Euler
scheme~\eqref{eq:implEuler} 
(compare
\eqref{eq:ArealConstant} and 
\eqref{eq:realNumberCepsilon}).
Nonetheless,
\eqref{eq:additionalInfo} and 
hence the numerical method~\eqref{eq:reducedScheme}
can be simulated easily.
The function 
$ \frac{ T }{ N^2 } 
\sum^N_{
  i = 1
} \eta_i \left( g_i \right)^2
$ in \eqref{eq:additionalInfo}
can be computed once in advance
for which $ O( N^2 ) $
computational operations are needed.
Having computed 
$ \frac{ T }{ N^2 } 
\sum^N_{ i = 1 }
\eta_i \left( g_i \right)^2
$, $ O( N \log(N) ) $
further computational operations
and random variables are needed
to 
compute 
\eqref{eq:additionalInfo}
from $ Y^N_{n} $
for one fixed 
$ n \in \left\{ 0,1, \dots, N^2 -1 
\right\} $ by using fast Fourier
transform.
Since 
$ O( N \log(N) ) $
computational operations and 
independent standard normal
random variables are needed
for one time step and 
since $ N^2 $ time steps are used
in \eqref{eq:reducedScheme},
$ O( N^3 \log(N) ) $
computational operations and
random variables are needed
to compute $ Y^N_{ N^2 } $.
Combining the computational effort
$ O( N^3 \log(N) ) $
and the convergence order $ \frac{3}{2}- $
in \eqref{eq:realNumberCepsilon}
shows that the Milstein type 
method~\eqref{eq:reducedScheme} 
needs about
$ O(\varepsilon^{-2}) $
computational operations and 
independent standard
normal random variables 
to achieve the desired 
precision $ \varepsilon > 0 $ 
in \eqref{eq:strongpr}.
To sum up, 
the Milstein type
algorithm~\eqref{eq:reducedScheme}
requires 
about $ O( \varepsilon^{-2} ) $
and the linear-impicit Euler
scheme~\eqref{eq:implEuler}
requires about
$ O( \varepsilon^{-\frac{8}{3}} ) $
computational operations and
independent standard normal
random variables for solving
the strong approximation
problem~\eqref{eq:strongpr} for 
the SPDE~\eqref{eq:SPDE2}.

\begin{figure}[ht]
\begin{center}
\includegraphics[width=10cm]{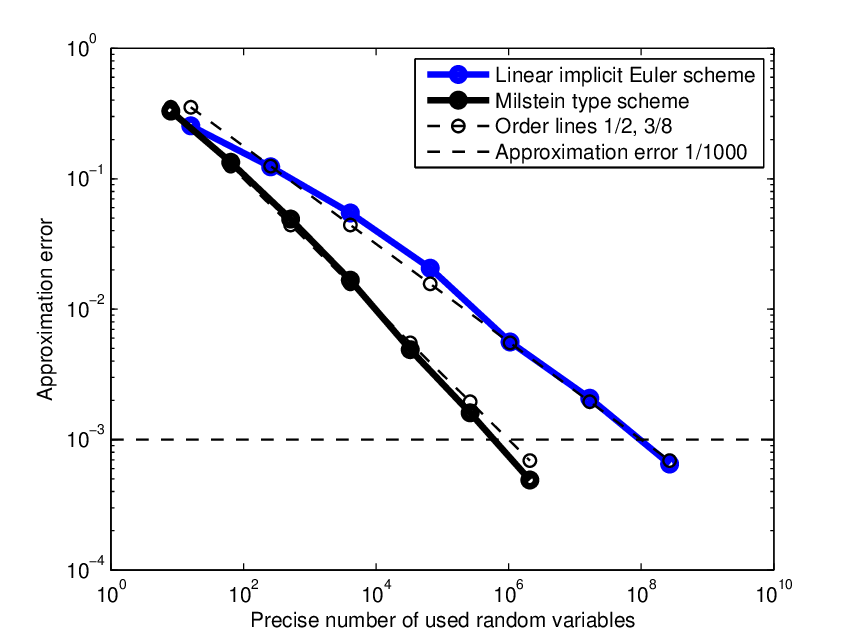}
\caption{\label{fig1}
SPDE~\eqref{eq:SPDE_reduced}:
Approximation error
in the sense of \eqref{eq:strongpr}
of the linear implicit Euler
approximation $ Z^N_{N^3} $
(see \eqref{eq:implEuler})
and of the Milstein type
approximation $ Y^N_{N^2}$
(see \eqref{eq:reducedScheme})
against the 
precise number of independent 
standard normal random variables
needed to compute the corresponding
approximation for 
$ N \in \left\{ 
2, 4, 8, 16,
32, 64, 128 \right\} $.
}
\end{center}
\end{figure}

The convergence rates
$ O( \varepsilon^{-\frac{8}{3}} ) $ and
$ O( \varepsilon^{-2} ) $ are both asymptotic
results as $ \varepsilon > 0 $ tends to zero.
Therefore, from a more 
practical point of view,
one may ask whether the 
Milstein type
algorithm~\eqref{eq:reducedScheme} 
solves the strong
approximation 
problem~\eqref{eq:strongpr}
more efficiently than the linear implicit 
Euler scheme~\eqref{eq:implEuler} for a given
concrete $ \varepsilon > 0 $ and a given
example of the form~\eqref{eq:SPDE2}.
In order to study this question we compare
both methods in the case of a simple
stochastic reaction diffusion equation.
More formally, let
$ \kappa = \frac{1}{100}$,
let 
$ 
  \xi \colon [0,1] \rightarrow \mathbb{R} 
$
be given by 
$ \xi( x ) = 0 $ for all
$ x \in [0,1] $ and suppose that
$ 
  f, b 
  \colon (0,1) \times \mathbb{R} 
  \rightarrow \mathbb{R} 
$
are given by $ f( x, y ) = 1 - y $ and
$ b( x, y ) = \frac{ 1 - y }{ 1 + y^2 } $
for all $ x \in (0,1) $, $ y \in \mathbb{R} $.
The SPDE~\eqref{eq:SPDE2} thus reduces to
\begin{equation}\label{eq:SPDE_reduced}
  dX_t( x )
  =
  \left[
    \frac{1}{100}
    \frac{ \partial^2 }{ 
      \partial x^2 } X_t( x )
    +
    1 - X_t( x )
  \right] \! dt
  +
  \frac{ 1 - X_t( x ) }{ 1 + X_t( x )^2 
  } \, dW_t( x )
\end{equation}
with $ X_t( 0 ) = X_t( 1 ) = 0 $ and $ X_0 = 0 $ for $ x \in (0,1) $ and $ t \in [0,1] $ (see also Section~\ref{sec:reacdiff} for more
details concerning
this example). 
Additionally, assume that
\eqref{eq:strongpr} 
for the 
SPDE~\eqref{eq:SPDE_reduced}
should be solved 
with the
precision of say 
three decimals, i.e.,
with the precision 
$ \varepsilon = 
\frac{ 1 }{ 1000 } $
in \eqref{eq:strongpr}.
In Figure~\ref{fig1}
the approximation error
in the sense of \eqref{eq:strongpr}
of the linear implicit Euler
approximation $ Z^N_{N^3} $
(see \eqref{eq:implEuler})
and of the  
approximation $ Y^N_{N^2}$
(see \eqref{eq:reducedScheme})
is plotted against
the precise number of independent 
standard normal random variables
needed to compute the corresponding
approximation for 
$ N \in \left\{ 
2, 4, 8, 16,
32, 64, 128 \right\} $:
It turns
out that $ Z_{ 128^3 }^{ 128 } $ 
in the case of the linear implicit Euler scheme~\eqref{eq:implEuler}
and 
that $ Y_{ 128^2 }^{ 128 } $ in the
case of the 
Milstein type
algorithm~\eqref{eq:reducedScheme} 
achieve the desired
precision
$ \varepsilon = \frac{ 1 }{ 1000 } $ in \eqref{eq:strongpr}
for the 
SPDE~\eqref{eq:SPDE_reduced}.
The {\sc Matlab} codes for simulating 
$ Z_{ 128^3 }^{ 128 } $ 
via \eqref{eq:implEuler}
and $ Y_{ 128^2 }^{ 128 } $ 
via \eqref{eq:reducedScheme}
for the SPDE~\eqref{eq:SPDE_reduced} 
are presented below
in Figure~\ref{code1} and Figure~\ref{code2} respectively. 
The
differences of the codes and the additional code
needed for the 
Milstein type
algorithm~\eqref{eq:reducedScheme} are
printed bold in Figure~\ref{code2}.
The {\sc Matlab} code 
in Figure~\ref{code1} requires
on our {\sc Intel Pentium D}
running at 3.0 GHz
a CPU time of about
$ 15 $~minutes and $25.03$ seconds
($925.03$ seconds) while the 
code in Figure~\ref{code2} 
requires a CPU time of 
about 
$8.93$ seconds
to be evaluated on the same 
computer.
So, on the above computer
the Milstein type
algorithm~\eqref{eq:reducedScheme} 
is for the 
SPDE~\eqref{eq:SPDE_reduced} 
about hundred times faster
than the linear implicit Euler 
scheme~\eqref{eq:implEuler}
in order to achieve a precision 
of three decimals 
in \eqref{eq:strongpr}.
Further numerical examples for the
Milstein type 
algorithm~\eqref{eq:reducedScheme}
can be found in 
Section~\ref{sec:examples}
below.
\begin{figure}[ht]
\begin{lstlisting}[firstnumber=1]
N = 128; M = N^3; A = -pi^2*(1:N).^2/100; Y = zeros(1,N);
mu = (1:N).^-2; f = @(x) 1-x; b = @(x) (1-x)./(1+x.^2);
for m=1:M
  y = dst(Y) * sqrt(2);
  dW = dst( randn(1,N) .* sqrt(mu*2/M) );
  y = y + f(y)/M + b(y).*dW;
  Y = idst( y ) / sqrt(2) ./ ( 1 - A/M );
end
plot( (0:N+1)/(N+1), [0,dst(Y)*sqrt(2),0] );
\end{lstlisting}
\caption{{\sc Matlab} code for 
simulating
the linear implicit Euler
approximation 
$ Z^{N}_{N^3} $
with $N = 128$
(see \eqref{eq:implEuler})
for the SPDE~\eqref{eq:SPDE_reduced}.}
\label{code1}
\end{figure}
\begin{figure}[ht]
\begin{lstlisting}[firstnumber=1]
N = 128; M = |N^2|; A = -pi^2*(1:N).^2/100; Y = zeros(1,N);
mu = (1:N).^-2; f = @(x) 1-x; b = @(x) (1-x)./(1+x.^2);
|bb = @(x) (1-x).*(x.^2-2*x-1)/2./(1+x.^2).^3; g = zeros(1,N);|
|for n=1:N|
  |g = g+2*sin(n*(1:N)/(N+1)*pi).^2*mu(n)/M;|
|end|
for m=1:M
  y = dst(Y) * sqrt(2);
  dW = dst( randn(1,N) .* sqrt(mu*2/M) );
  y = y + f(y)/M + b(y).*dW + |bb(y).*(dW.^2 - g)|;
  Y = |exp( A/M )| .* idst( y ) / sqrt(2);
end
plot( (0:N+1)/(N+1), [0,dst(Y)*sqrt(2),0] );
\end{lstlisting}
\caption{{\sc Matlab} code for 
simulating
the Milstein type
approximation 
$ Y^{N}_{N^2} $
with $N = 128$
(see \eqref{eq:reducedScheme})
for the SPDE~\eqref{eq:SPDE_reduced}.}
\label{code2}
\end{figure}

Having illustrated the efficiency of
the method~\eqref{eq:reducedScheme}, we
now take a short look at 
the literature of numerical 
analysis for SPDEs.
First, it should be mentioned that any
combination of finite elements, finite
differences or spectral Galerkin methods
for the spatial discretization and the 
linear implicit Euler scheme or also the
linear implicit Crank-Nicolson scheme for
the temporal discretization do
not reduce the computational
effort 
$ O( \varepsilon^{ -\frac{8}{3} } ) $
for the problem~\eqref{eq:strongpr} 
in case of the SPDE~\eqref{eq:SPDE2}. 
However,
the splitting-up method 
(see 
\cite{bgr90,bgr92,fl91,gk03b,gk03a,gk05,ir00}
and the references therein)
converges with a higher temporal
order than the linear 
implicit Euler scheme.
The key idea of the splitting-up 
method is to split the considered SPDE 
into appropriate subequations 
that are easier to solve
than the original SPDE, e.g.,
that can be solved
explicitly.
The applicability of the
splitting-up method thus 
essentially depends on the 
simplicity of the resulting 
subequations 
and can therefore 
in general not be 
used efficiently for 
nonlinear SPDEs.
Nonetheless, in the case of an 
appropriate class of linear SPDEs,
the splitting-up method 
and the Milstein type
scheme here
require nearly the same computational
effort for solving the strong
approximation problem
(see Section~\ref{sec:twoheat} 
for a more detailed comparison
of the splitting-up method 
and the Milstein type
algorithm here).
Additionally, in the case $ f = 0 $ 
in \eqref{eq:SPDE2}, 
T. M{\"u}ller-Gronbach and 
K. Ritter proposed a new 
scheme which
reduces the number of 
independent standard normal random
variables
needed
for solving a similar
problem as \eqref{eq:strongpr} 
from about
$ O( \varepsilon^{ - \frac{8}{3} } ) $ to
about $ O( \varepsilon^{ - 2 } ) $ (see
\cite{mr07b} and also
\cite{mr07a}). 
Nonetheless, the number of
computational operations needed 
and thus the overall 
computational effort 
could not be reduced by 
the algorithm in \cite{mr07b}. 
Moreover, Milstein type schemes for
SPDEs have been considered in \cite{gk96,ks01,lcp10,ms06}. 
In \cite{gk96}, W. Grecksch 
and P. E. Kloeden proposed
a Milstein like scheme 
for an SPDE driven by a scalar
one-dimensional Brownian 
motion (see also \cite{ks01}). 
In view 
of \eqref{eq:Milstein_rewritten}, 
the Milstein type scheme in \cite{gk96,ks01}
can be simulated efficiently 
since the driving
noise process is one-dimensional. 
Furthermore, 
in the case of a suitable
linear SPDE, 
A.~Lang, P.-L.~Chow
and J.~Potthoff constructed 
in the
interesting article \cite{lcp10}
a scheme 
similar to \eqref{eq:numScheme}
but with an additional term. 
(The additional
term may be useful for 
decreasing the error constant but 
turns out not to be needed in order to
achieve the higher approximation order due to Theorem~\ref{thm:mainresult} here.) 
In order to simulate
the iterated stochastic integral in their scheme, 
they then suggest to 
omit the summands 
in the double 
sum in \eqref{eq:commMilstein} for
which $ i \neq j $ holds 
(see (10) in \cite{lcp10}
and also \cite{lcp11}). Their
idea thus yields a scheme that can be simulated very
efficiently, but does in general not converge with a higher
order anymore, except for 
a linear SPDE driven by a scalar
one-dimensional Brownian motion.
Finally, based on It{\^o}'s formula,
Y. S. Mishura and 
G. M. Shevchenko 
proposed in \cite{ms06} the temporal
approximation \eqref{eq:tempApprox}
under additonal smoothness
assumptions of the driving noise
process of the 
SPDE~\eqref{eq:Dprocess_X}
which garantuee the existence of
a strong solution and thus allow
the application of It{\^o}'s formula
(see Assumption C in \cite{ms06}).
The simulation of the iterated
stochastic integrals in 
the Milstein type approximation
in \cite{ms06}
remained an open question
(see Remark 1.1 in \cite{ms06}).
To sum up, 
in the general setting of the 
possibly nonlinear SPDE~\eqref{eq:SPDE2},
the Milstein type
algorithm~\eqref{eq:reducedScheme}
is -- to the best of our
knowledge -- the first numerical 
approximation 
method which has
been shown to require
asymptotically
less computational operations 
and independent standard
normal random variables
than the 
required 
$ 
  O( \varepsilon^{ - \frac{8}{3} } ) 
$
of the
linear-implicit Euler scheme
in order to solve the
strong approximation
problem~\eqref{eq:strongpr} 
of the SPDE~\eqref{eq:SPDE2}.

The rest of this article is 
organized as follows. In
Section~\ref{sec:settings} the setting and the 
assumptions used are formulated. The numerical method and its 
convergence result (Theorem~\ref{thm:mainresult}) are
presented in Section~\ref{secscheme}. In Section~\ref{sec:examples} several examples of
Theorem~\ref{thm:mainresult} including a stochastic heat equation and stochastic reaction diffusion equations in one
and two dimensions
are considered. 
The
proof of Theorem~\ref{thm:mainresult} is postponed
to Section~\ref{secproofs}.

Next let us add some concluding
remarks.
There are a
number of directions for 
further research arising
from this work. One is to analyze whether the
exponential term 
in \eqref{eq:reducedScheme} can be 
replaced by a simpler mollifier such as 
$ ( I - \frac{T}{N^2} A )^{ -1 } $ for $ N \in \mathbb{N} $.
This would make the scheme 
even simpler to
simulate.
A second direction is to combine the temporal approximation 
in \eqref{eq:reducedScheme} with
other spatial discretizations such as finite elements.
This makes it possible to handle more complicated
multidimensional domains on which the eigenfunctions
of the Laplacian are not known explicitly. A third
direction is to reduce the 
Lipschitz assumptions in
Section~\ref{sec:settings} in order to 
handle further classes of SPDEs 
with non-globally Lipschitz 
nonlinearities such as 
stochastic Burgers equations,
stochastic porous medium 
equations and hyperbolic SPDEs.
Next a combination of the 
multilevel Monte Carlo approach 
(see \cite{h98,g08b}) 
with the Milstein type algorithm here
should result in a faster approximation
of statistical quantities 
of the solution
process of the SPDE~\eqref{eq:SPDE2}.
After a first preprint of this work
has appeared, 
a few research articles related 
to this work
have appeared; see 
\cite{bl11b,bl11c,w11,bls11,djr11}.
Some of the above outlined 
future research directions
and other issues such as
Runge-Kutta type schemes
for SPDEs 
based on the Milstein type scheme here
and
Milstein type schemes for SPDEs
driven by non-Gaussian noise
have been investigated in these articles.

Finally, let us point out limitations
of the Milstein type algorithm 
presented here.
There are essentially two assumptions
which need to be fulfilled so that
the above outlined Milstein type
algorithm can be applied.
First, the noise must be of 
{\it trace class} type.
Indeed, the covariance
operator 
$ 
  Q \colon H
  \rightarrow H 
$ 
of the Wiener
process is assumed to
be a trace class operator here
and we do not know how to
treat the case where 
$ 
  Q \colon H \rightarrow H
$
is bounded but without a finite
trace with
the above outlined Milstein type
algorithm.
Second, the diffusion coefficient
$ B $ must satisfy a certain
{\it commutativity type condition};
see Remark~\ref{rem:commutativity} below for details.
In the setting of the 
SPDE~\eqref{eq:SPDE2}
this condition is always fulfilled
and no restriction is required
even when the domain of the SPDE
becomes multi-dimensional dimensional
(see, e.g., 
Section~\ref{sec:twoheat})
but in the more general setting
of systems of SPDEs this results
in a serious restriction on the 
diffusion coefficient of the
considered system of SPDEs.
In addition to these two essential
restrictions,
several further restrictive
assumptions such as
the semilinear structure of the
SPDE,
the explicit knowledge of the
eigenfunctions of $ A $ and $ Q $
as well as
Lipschitz assumptions
on the nonlinearities $ F $
and $ B $
are used in this
article.
These further assumptions
are, however, no serious restrictions
since they can be relaxed or have
already been relaxed 
(see \cite{w11,bl11b,bl11c,djr11}).

\section{Setting and assumptions}\label{sec:settings}
Throughout this article suppose that the
following setting and the following assumptions
are fulfilled.
Let 
$ T \in (0, \infty) $, 
let 
$ 
  \left( \Omega, \mathcal{F}, \mathbb{P}   
  \right) 
$ 
be a probability space with a normal filtration
$( \mathcal{F}_t )_{t \in [ 0, T] }$ and let
$ 
  \left( 
    H, 
    \left< \cdot , \cdot \right>_H, 
    \left\| \cdot \right\|_H
  \right) $ 
and
$ 
  \left( 
    U, 
    \left< \cdot , \cdot \right>_U, 
    \left\| \cdot \right\|_U 
  \right) 
$ 
be two separable
$\mathbb{R}$-Hilbert spaces.
Moreover, let 
$ 
  Q \colon U \rightarrow U 
$ 
be a trace class operator and let 
$ 
  W \colon [0,T] \times \Omega 
  \rightarrow U 
$ 
be a standard 
$ Q $-Wiener process with respect to
$ ( \mathcal{F}_t )_{ t \in [0,T] } $.
\begin{assumption}[Linear operator A]\label{semigroup}
Let $ \mathcal{I} $ be a 
finite or countable set
and let $ \left( \lambda_i 
\right)_{i \in \mathcal{I} } $ $ \subset $ $(0, \infty) $
be a family of real numbers with 
$ \inf_{ i \in \mathcal{I} } \lambda_i \in (0,\infty) $. 
Moreover, let
$ (e_i)_{ i \in \mathcal{I} } $
be an orthonormal basis of $H$ and 
let 
$ 
  A \colon D(A) \subset H \rightarrow H 
$
be a linear operator with
\begin{equation}
  Av 
  =
  \sum_{ i \in \mathcal{I} }
    - \lambda_i 
    \left<
      e_i, v
    \right>_H
    e_i
\end{equation}
for all $v \in D(A)$ and with 
$
  D(A) 
  = 
  \big\{
    w \in H
    \colon
    \sum_{ i \in \mathcal{I} }
      \left|
        \lambda_i
      \right|^2
      \left|
        \left<
          e_i, w
        \right>_H
      \right|^2
    <
    \infty
  \big\}
$.
\end{assumption}
By 
$ 
  H_r := 
  D\!\left( \left( - A \right)^r \right) 
$
equipped with the norm
$ 
  \left\| v \right\|_{ H_r }
  := 
  \left\| ( - A )^r v \right\|_H 
$
for all $ v \in H_r $ 
and all $ r \in [0,\infty) $ we denote
the $ \mathbb{R} $-Hilbert spaces
of domains of fractional powers 
of the linear operator $ -A : D(A) \subset H
\rightarrow H $.
\begin{assumption}[Drift 
coefficient $F$]\label{drift}
Let $\beta \in [0,1) $ be a real number and 
let 
$ 
  F \colon H_{\beta} \rightarrow H 
$ 
be a twice
continuously Fr\'{e}chet differentiable
mapping with
$ 
  \sup_{ v \in H_{\beta} } 
  \left\| F'(v) \right\|_{ L(H) }
  < \infty
$
and 
$ 
  \sup_{ v \in H_{\beta} } 
  \left\| F''(v) \right\|_{ L^{(2)}(H_{\beta},H) }
  < \infty
$.
\end{assumption}
For formulating the 
assumption on the diffusion coefficient of our SPDE, denote by	
$ 
  \left( 
    U_0, 
    \left< \cdot , \cdot \right>_{ U_0 }, 
    \left\| \cdot \right\|_{ U_0 }
  \right) 
$ 
the separable
$\mathbb{R}$-Hilbert space 
$ U_0 := Q^{ 1/2 }( U ) $
with 
$ 
  \left< v, w \right>_{ U_0 } =  
  \left< Q^{ - 1 / 2 } v, 
  Q^{ - 1 / 2 } w \right>_U 
$ 
for all 
$ v, w \in U_0 $
(see, for example, Section~2.3.2 in \cite{pr07}).
For an arbitrary bounded linear operator $ S \in L( U ) $, we denote by 
$ 
  S^{-1} 
  \colon \text{im}( S ) 
  \subset U \rightarrow U 
$ 
the pseudo inverse of $ S $ (see, for instance, Appendix C in \cite{pr07}).
\begin{assumption}[Diffusion 
coefficient $B$]\label{diffusion}
Let 
$ 
  B \colon 
  H_{ \beta } \rightarrow HS(U_0,H) 
$ 
be a twice continuously Fr\'{e}chet differentiable
mapping with $ \sup_{ v \in H_{\beta} } \! 
\| B'( v ) 
\|_{ L( H, HS( U_0, H ) ) } \!
< \infty $ and 
$ \sup_{ v \in H_{\beta} } \!
\| B''( v ) 
\|_{ L^{(2)}( H_{\beta}, HS( U_0, H ) ) }
$ 
$<$ 
$\infty $. 
Moreover, 
let 
$ \alpha, c \in (0,\infty) $,
$ \delta, \vartheta 
\in ( 0, \frac{1}{2} )$,
$ \gamma \in [ \max(\delta,\beta), 
\delta + \frac{1}{2} ) $ 
be real numbers,
let $ B( H_{ \delta } ) $
$\subset$ $ HS(U_0, H_{ \delta } ) $
and suppose that
\begin{equation}\label{condition_AA}
  \left\| B(u) \right\|_{ HS( U_0, H_{ \delta } ) } 
  \leq
  c 
  \left( 
    1 + 
    \left\| u \right\|_{ H_{ \delta } } 
  \right) ,
\end{equation}
\begin{equation}\label{condition_A}
    \left\|
      B'\!\left( v \right) B\!\left( v \right)
      -
      B'\!\left( w \right) B\!\left( w \right)
    \right\|_{ HS^{ (2) }( U_0, H ) }
  \leq c
    \left\| v - w \right\|_H ,
\end{equation}
\begin{equation}\label{condition_AAA}
  \left\|
    ( - A )^{ - \vartheta } 
    B(v) Q^{ - \alpha }
  \right\|_{ HS( U_0, H ) }
  \leq
  c 
  \left( 
    1 + 
    \| v \|_{ H_{ \gamma } } 
  \right)
\end{equation}
for all $ u \in H_{ \delta } $ 
and all $ v, w \in H_{ \gamma } $.
Finally, let the bilinear Hilbert-Schmidt
operator $ B'(v) B(v) \in HS^{(2)}(U_0, H) $
be symmetric for all $ v \in H_{ \beta } $.
\end{assumption}
We now add some comments
concerning Assumption~\ref{diffusion}.
First, we note that 
Assumption~\ref{diffusion}
implies 
$ \beta \leq \delta + \frac{1}{2} $.
Indeed, $ \beta > \delta + \frac{1}{2} $
implies 
$ [ \max( \delta, \beta ), 
\delta + \frac{1}{2} ) = \emptyset $,
which contradicts to
$ \gamma \in [ \max( \delta, \beta ), 
\delta + \frac{1}{2} ) $
in Assumption~\ref{diffusion}.
Furthemore, we observe that
the above assumption
$ \sup_{ v \in H_{\beta} } 
\| B'( v ) 
\|_{ L( H, HS( U_0, H ) ) }
$ 
$ < $ $ \infty $
and the fact 
that $ H_{ \beta } $
is dense in $ H $ imply
that 
$ 
  B \colon H_{ \beta }
  \rightarrow HS(U_0, H) 
$
can be continuously extended 
to a globally Lipschitz continuous
mapping 
$ 
  \tilde{B} \colon H \rightarrow HS(U_0,H) 
$
from $H$ to $HS(U_0,H) $.
Here and below we do not distinguish
between 
$ 
  B \colon H_{\beta }
  \rightarrow HS(U_0, H) 
$ 
and its extension
$ 
  \tilde{B} \colon H \rightarrow
  HS( U_0, H ) 
$ 
for simplicity
of presentation.
Additionally, we note that
the operator 
$ 
  B'(v) B(v) 
  \colon
  U_0 \times U_0 \rightarrow H 
$ given by
\begin{equation}
\label{def:bilinear}
  \Big( B'(v) B(v) \Big) 
  ( u, \tilde{u} )
  =
  \Big(
    B'(v) \big( B(v) u \big)
  \Big) 
  ( \tilde{u} )
\end{equation}
for all $ u, \tilde{u} \in U_0 $ 
is a bilinear Hilbert-Schmidt operator in 
\begin{equation}
  HS^{(2)}( U_0, H ) \cong 
  HS( \overline{ U_0 \otimes U_0 }, H ) 
\end{equation}
for all $ v \in H_{\beta} $. 
Next we add a remark
on the symmetry assumption
on these bilinear Hilbert-Schmidt
operators
in Assumption~\ref{diffusion}.
\begin{remark}[Commutative noise
in infinite dimensions]
\label{rem:commutativity}
The assumed symmetry of 
the bilinear Hilbert-Schmidt
operator
$ 
  B'(v)B(v) \in HS^{ (2) }( U_0, H ) 
$ 
in Assumption~\ref{diffusion} reads as
\begin{equation}
\label{eq:commutativity}
  \Big(
    B'(v) \big( B(v) u \big)
  \Big) 
  ( \tilde{u} )
=
  \Big(
    B'(v) \big( B(v) \tilde{u} \big)
  \Big) 
  ( u )
\end{equation}
for all $ u, \tilde{u} \in U_0 $
and all $ v \in H_{ \beta } $.
Note that 
\eqref{eq:commutativity} is the 
abstract possibly infinite dimensional
coordinate free analog 
of \eqref{eq:comm_noise}. 
More formally,
if $ H = \mathbb{R}^d $, $ U = \mathbb{R}^m $ and
$ Q = I $ with $ d, m \in \mathbb{N} $ holds,
then \eqref{eq:commutativity} 
reduces to \eqref{eq:comm_noise}
(with $ \sigma $ replaced by $ B $).
\end{remark}

In Section~\ref{sec:examples} 
below
we describe a natural class
of examples of SPDEs which
satisfy the commutativity
condition~\eqref{eq:commutativity}
(see \eqref{eq:linopB} 
and \eqref{eq:checkcomm}
for details).
Finally, we emphasize that we do 
not assume that the linear operators 
$ 
  A \colon D( A ) \subset H \rightarrow H 
$ 
and 
$ 
  Q \colon U \rightarrow U 
$ 
are simultaneously diagonalizable. 
\begin{assumption}[Initial value $\xi$]\label{initial}
Let 
$ 
  \xi \colon \Omega 
  \rightarrow H_{\gamma} 
$ 
be 
an $ \mathcal{F}_0 
$/$ \mathcal{B}\left( 
H_{\gamma} \right) 
$-measurable 
mapping
with 
$ 
  \mathbb{E}\big[
    \| \xi \|^4_{ H_{\gamma} } 
  \big]
  < \infty 
$.
\end{assumption}
These assumptions suffice to ensure the existence of an up to modifications
unique solution of the SPDE~\eqref{eq:SPDE}.
\begin{prop}[Existence, uniqueness 
and regularity of 
solutions]\label{lem:existence}
Let Assumptions~\ref{semigroup}-\ref{initial} in Section~\ref{sec:settings} be fulfilled.
Then there exists an up to modifications 
unique predictable stochastic process
$
  X \colon [0,T] \times \Omega 
  \rightarrow H_{\gamma} 
$
which fulfills
$
  \sup_{ t \in [0,T] } 
  \mathbb{E}\big[
    \|
      X_t
    \|_{ H_{ \gamma } }^4 
  \big]
  < \infty
$,
$
  \sup_{ t \in [0,T] } 
  \mathbb{E}\big[
    \|
      B(X_t)
    \|_{ HS(U_0,H_{ \delta }) }^4
  \big]
  < \infty
$
and
\begin{equation}
\label{eq:SPDE}
    X_t
    = 
    e^{At} \xi 
    +
    \int_0^t e^{A(t-s)} F(X_s) \, ds 
    + 
    \int_0^t e^{A(t-s)} B(X_s) \, dW_s
\end{equation}
$ \mathbb{P} $-a.s.\ for 
all $ t \in [0,T] $.
Moreover, we have
\begin{equation}
\label{eq:holder}
  \sup_{ 
    \substack{
      t_1, t_2 \in [0,T] 
    \\
      t_1 \neq t_2
    }
  }
  \frac{
    \left(
      \mathbb{E}\big[
        \| X_{ t_2 } - X_{ t_1 }
        \|_{ H_{ r } }^4
      \big]
    \right)^{ \frac{1}{4} }
  }{
    \left| t_2 - t_1 
    \right|^{ 
      \min( \gamma - r, \frac{1}{2} ) 
    }
  }
< \infty
\end{equation}
for all $ r \in [0,\gamma] $.
Additionally, the solution process
$ X_t $, $ t \in [0,T] $,
is continuous with respect to
$ 
  \big( 
    \mathbb{E}\big[
      \| 
        \cdot
      \|^p_{ H_{ \gamma } }
    \big]
  \big)^{ 1/p }
$.
\end{prop}
\noindent
Proposition~\ref{lem:existence} immediately 
follows from Theorem~1 in \cite{jr11}. 
\section{Numerical scheme and main result}\label{secscheme}
In this section the numerical method 
is introduced and its convergence 
result is stated.
To this end let 
$ 
  \mathcal{J} 
$ 
be a set, 
let 
$ \left( g_j \right)_{ j \in \mathcal{J} } 
\subset U $ be an orthonormal basis of 
eigenfunctions of 
$ 
  Q \colon U \rightarrow U 
$ 
and 
let 
$ 
  \left( \eta_j \right)_{ j \in \mathcal{J} }   
  \subset 
  [0, \infty ) 
$ 
be the corresponding 
family of eigenvalues 
(such an orthonormal basis of eigenfunctions exists since 
$ 
  Q \colon U \rightarrow U 
$ 
is a trace class operator, 
see, e.g., Proposition~2.1.5 
in \cite{pr07}). 
In particular, 
we have 
\begin{equation}\label{eq:Qu} 
  Q u = 
  \sum_{ j \in \mathcal{J} } 
  \eta_j \left< g_j, u \right>_U g_j 
\end{equation}
for all $ u \in U $. 
Additionally, 
let $ \left( \mathcal{I}_N \right)_{ N \in \mathbb{N} } $ 
and $ \left( \mathcal{J}_K \right)_{ K \in \mathbb{N} } $ be 
sequences of finite subsets of $ \mathcal{I} $ and $ \mathcal{J} $ 
respectively.
Then we define the linear projection 
operators
$   
  P_N \colon H \rightarrow H 
$, 
$N \in \mathbb{N} $, 
by
$ 
  P_N( v ) := 
  \sum_{ i \in \mathcal{I}_N } 
  \left< e_i, v \right>_H e_i 
$
for all 
$ 
  v \in H 
$ 
and all 
$ 
  N \in \mathbb{N} 
$.
Furthermore, we define Wiener 
processes 
$ 
  W^K \colon [0,T] \times \Omega 
  \rightarrow U_0 
$, 
$ 
  K \in \mathbb{N} 
$, 
by 
\begin{equation}
W_t^K( \omega ) := 
\sum_{ 
  \substack{
    j \in \mathcal{J}_K 
  \\
    \eta_j \neq 0
  }
} 
\left< g_j, W_t( \omega ) \right>_U g_j 
\end{equation}
for all $ t \in [0,T] $, $ \omega \in \Omega $ 
and all $ K \in \mathbb{N} $.
We also use 
the 
$ 
  \mathcal{F} 
$/$ 
  \mathcal{B}(U_0) 
$-measurable
mappings 
$ 
  \Delta W^{M,K}_m 
  \colon \Omega \rightarrow U_0 
$,
$ 
  m \in 
  \left\{ 0, 1, \dots, M - 1 
  \right\} 
$,
$ 
  M, K \in \mathbb{N} 
$, 
given by
$ 
  \Delta W^{M,K}_m(\omega) := 
  W^K_{ \frac{ (m + 1) T }{ M } }(\omega) 
  - 
  W^K_{ \frac{ m T }{ M } }(\omega) 
$
for all 
$ \omega \in \Omega $,
$ 
  m \in 
  \left\{ 0, 1, \dots, M - 1 
  \right\} 
$ 
and all 
$ 
  M, K \in \mathbb{N} 
$.  
The numerical scheme 
which we denote by 
$ 
  \mathcal{F} 
$/$ 
  \mathcal{B}(H) 
$-measurable
mappings 
$ 
  Y_m^{N,M,K} \colon \Omega 
  \rightarrow H_N 
$, 
$ 
  m \in \{ 0, 1, \ldots, M \} 
$, 
$ 
  N,M,K \in \mathbb{N} 
$, 
is then given by 
$ 
  Y_0^{N,M,K} := P_N( \xi ) 
$ 
and
\begin{equation}
\label{eq:scheme}
\begin{split}
  Y_{m+1}^{N,M,K} 
& := 
  P_N \, e^{ A \frac{ T }{ M } }
  \Bigg(
    Y_m^{N,M,K} 
    +
    \frac{T}{M} \cdot 
    F\!\left( Y_m^{N,M,K} \right)
    +
    B\!\left( Y_m^{N,M,K} 
    \right) \!
    \Delta W_m^{M,K}
\\& \quad
    +
    \frac{1}{2}
    B'\!\left( Y_m^{N,M,K} \right) \!
    \Big( \!
      B\!\left( Y_m^{N,M,K} \right) \!
      \Delta W_m^{M,K}
    \Big) 
    \Delta W_m^{M,K}
\\ & \quad
    - \frac{T}{2 M} \!\!
      \sum_{ 
        \substack{
          j \in \mathcal{J}_K 
        \\
          \eta_j \neq 0
        }
      } 
        \!\!
        \eta_j
        B'\!\left( Y_m^{N,M,K} \right) \!
        \Big( \!
          B\!\left( Y_m^{N,M,K} \right) 
          g_j
        \Big) 
        g_j
  \Bigg)
\end{split}
\end{equation}
for all $ m \in \{0,1,\ldots,M-1\} $ and 
all $ N,M,K \in \mathbb{N} $.
Note that only increments of the 
Wiener processes 
$ 
  W^K \colon [0,T] \times \Omega 
  \rightarrow U_0 
$,
$ 
  K \in \mathbb{N} 
$, 
are used 
in the scheme above and 
we emphasize that for many 
SPDEs the method~\eqref{eq:scheme} 
is easy to simulate 
and 
to implement (see 
Sections~\ref{sec:intro} and
\ref{sec:examples}
for a few examples). 
Moreover, observe that
the scheme~\eqref{eq:scheme}
can also be written as
\begin{equation}
\label{eq:scheme}
\begin{split}
  Y_{m+1}^{N,M,K} 
& = 
  P_N \, e^{ A \frac{ T }{ M } }
  \Bigg(
    Y_m^{N,M,K} 
    +
    \frac{T}{M} \cdot 
    F\!\left( Y_m^{N,M,K} \right)
    +
    B\!\left( Y_m^{N,M,K} 
    \right) \!
    \Delta W_m^{M,K}
\\& \quad
    +
    \frac{1}{2}
    B'\!\left( Y_m^{N,M,K} \right) \!
    \Big( \!
      B\!\left( Y_m^{N,M,K} \right) \!
      \Delta W_m^{M,K}
    \Big) 
    \Delta W_m^{M,K}
\\ & \quad
    -
    \frac{1}{2} 
    \mathbb{E}\!\left[
    B'\!\left( Y_m^{N,M,K} \right) \!
    \Big( \!
      B\!\left( Y_m^{N,M,K} \right) \!
      \Delta W_m^{M,K}
    \Big) 
    \Delta W_m^{M,K}
    \, \big| \,
    \mathcal{F}_{ \frac{ m T }{ M } }
    \right]
  \Bigg)
\end{split}
\end{equation}
for all $ m \in \{0,1,\ldots,M-1\} $ and 
all $ N,M,K \in \mathbb{N} $.
In the next step the convergence result 
for the scheme~\eqref{eq:scheme}
is presented.
\begin{theorem}[Main
result]\label{thm:mainresult}
Let Assumptions~\ref{semigroup}-\ref{initial} in Section~\ref{sec:settings} be 
fulfilled. Then there is 
a real number $ C \in (0,\infty) $ such that
\begin{equation}
\label{eq:mainresult}
    \left(
      \mathbb{E}\!
       \left[
      \big\| 
        X_{ \frac{m T}{M} }
        - Y^{N,M,K}_m
      \big\|_H^2
       \right]
    \right)^{ \frac{1}{2} }
    \leq 
    C 
  \left(
      \left(
        \inf_{ 
          i \in \mathcal{I} \backslash
          \mathcal{I}_N
        }
        \lambda_i
      \right)^{ \! - \gamma }
    +
     \bigg(
      \sup_{ 
         j \in \mathcal{J} 
        \backslash \mathcal{J}_K 
      } 
      \eta_j 
      \bigg)^{ \! \alpha }
    +
      M^{ 
        - \min\left( 
          2 \left( \gamma - \beta \right) , 
          \gamma 
        \right) 
      }
  \right)
\end{equation}
for 
all $m \in \left\{ 0, 1, \dots, M \right\}$ 
and all $ N,M,K \in \mathbb{N} $.
\end{theorem}
We now explain the result of Theorem~\ref{thm:mainresult} 
in more detail. 
The strong root mean square difference 
$ 
  \big( 
    \mathbb{E}\big[ 
      \| 
        X_{ \frac{ m T }{ M } } -  
        Y_m^{ N,M,K } 
      \|_H^2 
    \big] 
  \big)^{ 1 / 2 } 
$ for $ m \in \{0,1,\ldots,M \}$ and for $ N,M,K \in \mathbb{N} $ of the exact solution of the SPDE~\eqref{eq:SPDE} and of the numerical solution~\eqref{eq:scheme} is estimated in \eqref{eq:mainresult} by a constant times the sum of three terms. The first term, i.e.,
$
  (
     \inf_{ 
        i \in \mathcal{I} \backslash
        \mathcal{I}_N
      }
      \lambda_i
  )^{ - \gamma }
$ 
for $ N \in \mathbb{N} $,
arises due to discretizing 
the exact solution spatially, 
i.e., 
due to 
$ 
  \mathbb{E}\big[ 
    \| X_t - P_N( X_t ) \|_H^2 
  \big]
$ 
for 
$ N \in \mathbb{N} $ 
and 
$ t \in [0,T] $. 
The second expression, i.e., 
$
  (
    \sup_{ 
      j \in \mathcal{J} 
      \backslash \mathcal{J}_K 
    } 
    \eta_j 
  )^{ \alpha }
$ 
for $ K \in \mathbb{N} $,
occurs due to discretizing the noise spatially, i.e., 
due to 
$ 
  \mathbb{E}\big[
    \| W_t - W_t^K \|_H^2 
  \big]
$ 
for $ K \in \mathbb{N} $ 
and $ t \in [0,T] $. 
If $ U_0 \subset U $ 
is finite dimensional we choose 
$ \mathcal{J}_K := \left\{ j \in \mathcal{J} \big| 
\eta_j \neq 0 \right\} $ 
for all $ K \in \mathbb{N} $ 
and obtain
$
  (
    \sup_{ 
      j \in \mathcal{J} 
      \backslash \mathcal{J}_K 
    } 
    \eta_j 
  )^{ \alpha }
  =
  0
$ 
for all $ K \in \mathbb{N} $ in 
that case. The third 
term, i.e., 
$ 
  M^{ - \min\left( 
    2 \left( \gamma - \beta \right) , 
    \gamma \right) 
  } 
$ 
for 
$ 
  M \in \mathbb{N} 
$, 
corresponds to the temporal 
discretization error and 
converges to zero as the number 
of time steps $ M \in \mathbb{N} $ 
goes to infinity.

\section{Examples}\label{sec:examples}
In this section Theorem~\ref{thm:mainresult} 
is illustrated with various examples. To this end let
$ d \in \{1,2,3\} $ and let $ H = U = L^2( (0,1)^d, \mathbb{R} ) $
be the $ \mathbb{R} $-Hilbert space of equivalence
classes of
Lebesgue square integrable 
functions from 
$ (0,1)^d $ to $ \mathbb{R} $. 
As usual we do not distinguish between 
a 
Lebesgue square integrable 
function from $ (0,1)^d $ to $ \mathbb{R} $ and its equivalence 
class in $ H $. 
The scalar product and the
norm in $ H $ and $ U $ are given by
\begin{equation}
  \left< v,w \right>_H
  =
  \left< v,w \right>_U
  =
  \int_{ (0,1)^d }
  v( x ) \cdot w( x ) \, dx
\end{equation}
and
\begin{equation}
  \left\| v \right\|_H
  =
  \left\| v \right\|_U
  =
  \left(
    \int_{ (0,1)^d }
    \left| v( x ) \right|^2 dx
  \right)^{ \frac{1}{2} }
\end{equation}
for all $ v,w \in H = U $.
Additionally, the notations
\begin{equation}
  \left\| v \right\|_{ C\left( (0,1)^d, \mathbb{R} \right) }
  :=
  \sup_{ x \in (0,1)^d }
  \left| v( x ) \right|
  \in [0,\infty]
\end{equation}
and
\begin{equation}
  \left\| v \right\|_{ C^r\left( (0,1)^d, \mathbb{R} \right) }
  :=
  \sup_{ x \in (0,1)^d }
  \left| v( x ) \right|
  +
  \sup_{ \substack{ x, y \in (0,1)^d \\ x \neq y } }
  \frac{ \left| v( x ) - v( y ) \right| }{
  \left\| x - y \right\|^r_{ \mathbb{ R }^d } }
  \in [0,\infty]
\end{equation}
are used troughout this section for all functions
$ 
  v \colon (0,1)^d 
  \rightarrow \mathbb{ R } 
$ 
and 
all $ r \in (0,1] $.
Here and below we use 
the Euclidean 
norms 
$ 
  \left\| x \right\|_{ \mathbb{R}^n }
  := 
  ( | x_1 |^2 + \ldots + | x_n |^2 )^{1/2} 
$
for all 
$ 
  x = \left( x_1, \dots, x_n \right)
  \in \mathbb{R}^n 
$ 
and all 
$ 
  n \in \mathbb{N} 
$.
Concerning the Wiener process 
$ 
  W \colon [ 0,T ] \times \Omega 
  \rightarrow U 
$ 
we 
assume that the 
eigenfunctions $ g_j \in U $, $ j \in \mathcal{J} $, 
of the covariance operator 
$ 
  Q \colon U \rightarrow U 
$ 
are continuous and satisfy
\begin{equation}\label{fcond}
  \sup_{ j \in \mathcal{J} }
  \left\| g_j \right\|_{ C\left( (0,1)^d, \mathbb{R} \right) }
  < \infty 
\qquad
  \text{and}
\qquad
  \sum_{ j \in \mathcal{J} }
  \left(
  \eta_j
  \left\| g_j \right\|^2_{ C^{ \rho }\left( (0,1)^d, 
  \mathbb{R} \right) } 
  \right)
  < \infty 
\end{equation}
for some 
$ \rho \in \left( 0, 1 \right) $
in this section.
We will give some concrete examples for the $ \left( g_j \right)_{ j \in \mathcal{J} } $ that 
fulfill \eqref{fcond} later.
Additionally, we denote by $ s_n( D ) \in [0,\infty) $, $ n \in \mathbb{N} $, the sequence of characteristic numbers of a compact operator 
$ 
  D \colon H \rightarrow H 
$ 
(see, e.g., 
Section~9 in Chapter XI in \cite{ds88}). 
Finally, we define the Schatten norms by
\begin{equation}
  \left\| D \right\|_{ S_p( H ) }
  :=
  \left(
    \sum_{ n=1 }^{ \infty }
    \left| s_n( D ) \right|^p
  \right)^{ \! \frac{1}{p} }
  \in [0,\infty]
\end{equation}
for all compact operators 
$ 
  D \colon H \rightarrow H 
$ 
and all 
$ p \in [1,\infty) $ 
(see also the above named reference).

We now present a prominent example
of {\bf the linear operator $A$ in
Assumption 1}.
Let $ \mathcal{I} = \mathbb{N}^d $ and 
let $ e_i \in H $ for $ i \in \mathcal{I} $ 
be given by
\begin{equation}
  e_i( x )
  =
  2^{ \frac{d}{2} } 
  \sin( i_1 \pi x_1 ) 
  \cdot \ldots \cdot 
  \sin( i_d \pi x_d )
\end{equation}
for all $ x = ( x_1, \ldots, x_d ) \in (0,1)^d $ and all $ i= ( i_1, \ldots, i_d ) \in \mathbb{N}^d $.
Additionally, 
let $ \kappa \in (0,\infty) $ be a 
fixed real number and let 
$ 
  \left( \lambda_i 
  \right)_{ i \in \mathcal{I} } 
$ 
be given by
\begin{equation}
  \lambda_i
  = 
  \kappa \, 
  \pi^2 \left( (i_1)^2 + \ldots + (i_d)^2 \right) 
\end{equation}
for all $ i=(i_1, \ldots, i_d) \in \mathbb{N}^d $. Hence, the linear operator 
$ 
  A \colon 
  D( A ) \subset H \rightarrow H 
$ 
in Assumption~\ref{semigroup} 
reduces to the Laplacian with Dirichlet boundary conditions times the constant 
$ \kappa \in (0,\infty) 
$, 
i.e., 
\begin{equation}
  A v = \kappa \cdot \Delta v 
  =
  \kappa 
  \left\{
    \left( 
      \frac{ \partial^2 }{ \partial x_1^2 }
    \right) \! v
    + \ldots +
    \left( 
      \frac{ \partial^2 }{ \partial x_d^2 }
    \right) \! v
  \right\}
\end{equation}
for all $ v \in D( A ) $ 
in this section. 
Furthermore, let 
$ \left( \mathcal{I}_N \right)_{ N \in \mathbb{N} } $ 
be given by $ \mathcal{I}_N = \{1,\ldots,N\}^d $ 
for all $ N \in \mathbb{N} $.

In order to describe a natural
candidate for 
{\bf the drift coefficient in Assumption~\ref{drift}}, 
let $ \beta = \frac{ d }{ 5 } $ and 
let 
$ 
  f \colon (0,1)^d \times 
  \mathbb{R} \rightarrow \mathbb{R} 
$ 
be a twice continuously differentiable function 
with $ \int_{ (0,1)^d } \left| f( x,0 ) \right|^2 dx < \infty $
and 
$
  \sup_{ x \in (0,1)^d }
  \sup_{ y \in \mathbb{R} }
  \big| 
    \big( 
      \frac{ \partial^n 
    }{
      \partial y^n
    } 
    f 
    \big)( x , y ) 
  \big|
  < \infty
$
for all $ n \in \left\{ 1,2 \right\} $. 
Then the (in general nonlinear)
operator 
$ 
  F \colon H_{ \beta } \rightarrow H 
$ 
given by
\begin{equation}
  \big( F( v ) \big)( x )
  =
  f( x, v( x ) )
\end{equation}
for all $ x \in (0,1)^d $
and all $ v \in H_{ \beta } $
satisfies Assumption~\ref{drift} since 
$ 
  H_{ \beta } 
  = 
  H_{ \frac{d}{5} } 
  \subset L^5( (0,1)^d, \mathbb{R} ) 
$ 
continuously 
(see~Remark~6.94 
in \cite{rr93}).

In the next step a natural example
for {\bf the diffusion 
coefficient in Assumption~\ref{diffusion}} 
is given. 
Let 
$ 
  b \colon 
  (0,1)^d \times \mathbb{R} 
  \rightarrow \mathbb{R} 
$ 
be a twice 
continuously differentiable function 
with 
\begin{equation}
\label{eq:diff_term_assumpt}
  \left| b( x,0 ) \right| \leq q,
  \qquad
  \left| \left( \frac{ \partial^n }{
    \partial y^n
  } b \right)\!( x , y ) 
  \right|
  \leq q ,
  \qquad
  \left\| \left( \frac{ \partial }{
    \partial x
  } b \right)\!( x , y ) 
  \right\|_{ L( \mathbb{R}^d, \mathbb{R} ) }
  \leq q
\end{equation}
and
\begin{equation}
\label{eq:partial}
  \left|
    \left(
      \frac{ \partial }{ \partial y } b 
    \right)\!(x,y) \cdot b(x,y)
    -
    \left(
      \frac{ \partial }{ \partial y } b 
    \right)\!(x,z) \cdot b(x,z)
  \right|
  \leq q \left| y - z \right| 
\end{equation}
for all $ x \in (0,1)^d $, $ y, z \in \mathbb{R} $, $ n \in \{ 1,2 \} $ and 
some given $ q \in (0,\infty) $.
We refer to 
Subsections~\ref{sec:reacdiff}-\ref{sec:twoheat}
below for concrete functions
$ 
  b \colon (0,1)^d \times \mathbb{R}
  \rightarrow \mathbb{R} 
$
satisfying 
\eqref{eq:diff_term_assumpt}
and \eqref{eq:partial}.
Now let 
$ 
  B \colon H_{ \beta } \rightarrow 
  HS( U_0, H ) 
$ 
be
the (in general nonlinear) operator
\begin{equation}\label{eq:linopB}
  \Big( B(v) u \Big)( x )
  =
  b( x, v( x ) )
  \cdot u( x )
\end{equation}
for all $ x \in (0,1)^d $,
$ v \in H_{ \beta } $ 
and all $ u \in U_0 \subset U = H $. 
We now check step by step 
that 
$ 
  B \colon H_{ \beta } 
  \rightarrow HS( U_0, H ) 
$ 
given by \eqref{eq:linopB} satisfies 
Assumption~\ref{diffusion}. First of all, 
$ B $ is well defined. 
More precisely, we have
\begin{align*}
&
  \left\| B( v ) \right\|_{ HS( U_0, H ) }^2
=
  \sum_{ j \in \mathcal{J} }
  \left\| B(v)
    \sqrt{ \eta_j } g_j \right\|_{ H }^2
=
  \sum_{ j \in \mathcal{J} }
  \eta_j 
  \left\| B(v) g_j 
  \right\|_{ H }^2
\\&=
  \sum_{ j \in \mathcal{J} }
  \eta_j \left(
    \int_{ (0,1)^d }
    \left| b( x , v( x ) ) \cdot g_j( x ) \right|^2 dx
  \right)
\leq 
  \sum_{ j \in \mathcal{J} }
  \eta_j
    \left(
      \int_{ (0,1)^d }
      \left| b( x , v( x ) ) \right|^2 dx
    \right)
    \left(
      \sup_{ x \in (0,1)^d }
      \left| g_j( x ) \right|^2
    \right)
\end{align*}
and hence
\begin{align*}
  \left\| B( v ) \right\|_{ HS( U_0, H ) }^2
&\leq 
  \sum_{ j \in \mathcal{J} }
  \eta_j
  \left(
    \int_{ (0,1)^d }
      \left( 
        \left| b( x , v( x ) )  - b( x, 0 ) 
        \right| 
      + \left| b( x, 0 ) \right|
      \right)^2 dx
  \right)
  \left\| g_j \right\|_{ C\left( (0,1)^d, 
      \mathbb{R} \right) }^2
\\&\leq 
  \sum_{ j \in \mathcal{J} }
  q^2 \eta_j
  \left(
    \int_{ (0,1)^d }
    \left(
    \left| v( x ) \right|
    + 1 \right)^2 dx
  \right)
  \left\| g_j 
  \right\|_{ C\left( (0,1)^d, \mathbb{R} \right) }^2
\end{align*}
and finally
\begin{align}\label{eq:bounded_B}
\nonumber
 \left\| B( v ) \right\|_{ HS( U_0, H ) } 
  &\leq 
  q 
    \left(
      \left\| v \right\|_H + 1
    \right)
  \left(
    \sum_{ j \in \mathcal{J} }
    \eta_j
    \left\| g_j 
    \right\|_{ C\left( (0,1)^d, \mathbb{R} \right) }^2
  \right)^{ \frac{1}{2} }
\\&\leq 
  q
  \left(
    \left\| v \right\|_H + 1
  \right)
  \left(
    \sum_{ j \in \mathcal{J} }
    \eta_j
  \right)^{ \frac{ 1 }{ 2 } }
  \left(
    \sup_{ j \in \mathcal{J} }
    \left\| g_j 
    \right\|_{ C\left( (0,1)^d, \mathbb{R} \right) }
  \right)
\\&=
\nonumber
  q
  \sqrt{ \text{Tr}( Q ) }
  \left(
    \sup_{ j \in \mathcal{J} }
    \left\| g_j 
    \right\|_{ 
      C\left( (0,1)^d, \mathbb{R} \right) 
    }
  \right)
  \left(
    \left\| v \right\|_H + 1 
  \right)
  < \infty
\end{align}
for all 
$ v \in H_{\beta} $ 
which indeed
shows that $ B $ is well defined.
Moreover, $ B $ is twice continuously 
Fr{\'e}chet differentiable and we have
\begin{align*}
&
  \left\| B'( v ) u \right\|_{ HS( U_0, H ) }^{2}
=
  \sum_{ j \in \mathcal{J} }
  \left\| 
    \left( B'( v ) u \right) 
    \sqrt{ \eta_j } g_j 
  \right\|_{ H }^{2}
=
  \sum_{ j \in \mathcal{J} }
  \eta_j
  \left\| 
    \left( B'( v ) u \right) 
    g_j 
  \right\|_{ H }^{2}
\\&=
  \sum_{ j \in \mathcal{J} }
  \eta_j
  \left(
    \int_{ (0,1)^d }
    \left| 
      \left( \frac{ \partial }{ \partial y } b \right)
      \!\left( x, v( x ) \right)
      \cdot
      u( x )
      \cdot
      g_j( x )
    \right|^2 dx
  \right)
\leq 
  \sum_{ j \in \mathcal{J} }
  q^2 \eta_j
  \left(
    \int_{ (0,1)^d }
    \left| 
      u( x ) \cdot g_j( x )
    \right|^2 dx
  \right)
\end{align*}
and hence
\begin{align*}
&
  \left\| B'( v ) u \right\|_{ HS( U_0, H ) }
\leq
  \left(
  \sum_{ j \in \mathcal{J} }
  q^2 \eta_j
  \left\| u \right\|_H^2
  \left\| g_j \right\|_{ C\left( (0,1)^d,
    \mathbb{R} \right) }^2
  \right)^{ \! 1 / 2 }
\\&\leq
  q
  \left\| u \right\|_H
  \left(
    \sum_{ j \in \mathcal{J} }
    \eta_j
  \right)^{ \! 1 / 2 }
  \left( 
    \sup_{ j \in \mathcal{J} }
    \left\| g_j 
    \right\|_{ C\left( (0,1)^d, \mathbb{R} \right) }
  \right)
=
  q
  \sqrt{ \text{Tr}( Q ) }
  \left( 
    \sup_{ j \in \mathcal{J} }
    \left\| g_j 
    \right\|_{ C\left( (0,1)^d, \mathbb{R} \right) }
  \right)
  \left\| u \right\|_H
\end{align*}
for all 
$ 
  u, v \in H_{\beta} 
$ 
which shows
\begin{equation}
  \sup_{ v \in H_{ \beta } }
  \left\| B'( v ) \right\|_{ 
    L( H, HS( U_0, H ) ) 
  }
  \leq
  q
  \sqrt{ \text{Tr}( Q ) }
  \left( 
    \sup_{ j \in \mathcal{J} }
    \left\| g_j 
    \right\|_{ 
      C\left( (0,1)^d, \mathbb{R} \right) 
    }
  \right)
  < \infty .
\end{equation}
Additionally, we have
\begin{align*}
  \left\| B''( v ) (u,w) 
  \right\|_{ HS( U_0, H ) }^2
&=
  \sum_{ j \in \mathcal{ J } }
  \left\| 
    B''( v ) (u,w) 
    \sqrt{ \eta_j } g_j
  \right\|_{ H }^2
\\&=
  \sum_{ j \in \mathcal{ J } }
  \eta_j
  \left(
    \int_{ (0,1)^d }
    \left|
      \left( 
        \frac{ \partial^2 }{ \partial y^2 } 
        b 
      \right) \!
      \left( x, v(x) \right)
      \cdot u(x) \cdot w(x) \cdot g_j(x) 
    \right|^2 dx
  \right)
\\&\leq
  \sum_{ j \in \mathcal{ J } }
  q^2 \eta_j
  \left(
    \int_{ (0,1)^d }
    \left|
      u(x) \cdot w(x)
    \right|^2 dx
  \right)
  \left\| g_j 
  \right\|^2_{ C\left( (0,1)^d, \mathbb{R} \right) }
\end{align*}
and using 
$ 
  L^5\!\left( (0,1)^d, \mathbb{R} \right) 
  \subset
  L^4\!\left( (0,1)^d, \mathbb{R} \right) 
$
continuously shows
\begin{align*}
 &\left\| B''( v ) (u,w) \right\|_{ HS( U_0, H ) }
\\
  &\leq
  q \sqrt{ \mathrm{Tr}(Q) }
  \left(
    \int_{ (0,1)^d }
    \left|
      u(x)
    \right|^4 dx
  \right)^{\!\!\frac{1}{4}}
  \left(
    \int_{ (0,1)^d }
    \left|
      w(x)
    \right|^4 dx
  \right)^{\!\!\frac{1}{4}}
  \left(
    \sup_{ j \in \mathcal{J} }
    \left\| g_j 
    \right\|_{ C\left( (0,1)^d, \mathbb{R} \right) }
  \right)
\\&
  \leq
  q
  \sqrt{ \text{Tr}( Q ) }
  \left(
    \sup_{ j \in \mathcal{J} }
    \left\| g_j 
    \right\|_{ C\left( (0,1)^d, \mathbb{R} \right) }
  \right)
  \left( 
    \int_{ (0,1)^d }
    \left| u( x ) \right|^5 dx
  \right)^{\!\! \frac{1}{5} }
  \left( 
    \int_{ (0,1)^d }
    \left| w( x ) \right|^5 dx
  \right)^{\!\! \frac{1}{5} }
\end{align*}
for all $ u, v, w \in H_{\beta} $. 
Therefore, 
$ 
  H_{ \beta } \subset 
  L^5\left( (0,1)^d, \mathbb{R} \right) 
$ 
continuously shows
\begin{equation}
  \sup_{ v \in H_{\beta} }
  \left\| 
    B''( v ) 
  \right\|_{ 
    L^{(2)}( H_{\beta}, HS( U_0, H ) ) 
  }
  < \infty 
\end{equation}
due to \eqref{fcond} and hence, it remains to establish \eqref{condition_AA}-\eqref{condition_AAA} and the symmetry of $ B'( v ) B( v ) \in HS^{(2)}( U_0. H ) $ for all 
$ 
  v \in H_{\beta} 
$.
For the latter one, 
note that
\begin{equation}
\label{eq:checkcomm}
  \bigg(
    \Big(
    B'( v ) B( v )
    \Big)
    ( u, \tilde{u} )
  \bigg)( x )
  =
  \bigg(
    B'( v )
    \Big(
      B( v ) u
    \Big)
    \tilde{u}
  \bigg)( x )
=
  \left(
    \frac{ \partial }{ \partial y } b
  \right)\!
  \left( x, v( x ) \right)
  \cdot 
  b\!\left( x, v( x ) \right)
  \cdot
  u( x )
  \cdot 
  \tilde{u}( x )
\end{equation}
for all 
$ x \in (0,1)^d $, 
$ u, \tilde{u} \in U_0 $ 
and all $ v \in H_{\beta} $ 
which immediately shows that 
$ B'( v ) B( v ) \in 
HS^{(2)}( U_0, H ) $ 
is symmetric 
for all $ v \in H_{\beta} $
(see \eqref{eq:commutativity}). 
Moreover, we have
\begin{equation*}
\begin{split}
&
  \left\|
    B'( v ) B( v )
    -
    B'( w ) B( w )
  \right\|^2_{ HS^{(2)}( U_0, H ) }
=
  \sum_{ j,k \in \mathcal{J} }
  \eta_j \eta_k
  \left\|
    B'( v )
    \left( B( v ) g_j \right) g_k
    -
    B'( w )
    \left( B( w ) g_j \right) g_k
  \right\|_H^2
\\ & 
  \leq 
  \sum_{ j, k \in \mathcal{J} }
  \eta_j \eta_k
  \Bigg(
    \int_{ (0,1)^d }
    \Bigg|
      \left( 
        \frac{ \partial }{ \partial y } b 
      \right)\!
      \left( x, v( x ) \right)
      \cdot 
      b\!\left( x, v(x) \right)
      -
      \left( 
        \frac{ \partial }{ \partial y } b 
      \right)\!
      \left( x, w( x ) \right)
      \cdot 
      b\!\left( x, w(x) \right)
    \Bigg|^2 dx
  \Bigg)
  \left[
    \sup_{ l \in \mathcal{J} }
    \left\| g_l \right\|^4_{ 
      C\left( (0,1)^d, \mathbb{R} \right) 
    }
  \right]
\end{split}
\end{equation*}
and using \eqref{eq:partial} yields
\begin{align*}
  \left\|
    B'( v ) B( v )
    -
    B'( w ) B( w )
  \right\|_{ HS^{(2)}( U_0, H ) }
&\leq 
  q
  \left\| v - w \right\|_H
  \left(
    \sum_{ j,k \in \mathcal{J} }
    \eta_j \eta_k
  \right)^{ \frac{1}{2} }
  \left(
    \sup_{ i \in \mathcal{J} }
    \left\| g_j \right\|^2_{ C\left( (0,1)^d, \mathbb{R} \right) }
  \right)
\\&=
  q
  \text{Tr}( Q )
  \left(
    \sup_{ i \in \mathcal{J} }
    \left\| g_j \right\|^2_{ C\left( (0,1)^d, \mathbb{R} \right) }
  \right)
  \left\| v - w \right\|_H
\end{align*}
for all $ v, w \in H_{ \beta } $ which shows that \eqref{condition_A} indeed holds.
Estimates~\eqref{condition_AA}
and \eqref{condition_AAA} will
be verified in the more concrete
examples in 
Subsections~\ref{sec:reacdiff}-\ref{sec:twoheat} below.

Concerning {\bf the initial value 
in Assumption \ref{initial}},
let 
$ 
  x_0 \colon
  [0,1]^d \rightarrow \mathbb{R} 
$ 
be a twice
continuously differentiable function 
with 
$ 
  x_0|_{ \partial (0,1)^d } \equiv 0 
$.
Then the 
$ 
  \mathcal{F}_0 
$/$ 
  \mathcal{B}( H_{\gamma} ) 
$-measurable 
mapping 
$ 
  \xi \colon \Omega \rightarrow 
  H_{ \gamma } 
$
given by 
$ 
  \xi( \omega ) = x_0 
$ 
for all 
$ 
  \omega \in \Omega 
$
fulfills 
Assumption~\ref{initial} 
for 
all $ \gamma \in (0, 1) $.

Having constructed examples for 
Assumptions~\ref{semigroup}-\ref{initial}, 
we now formulate 
the 
SPDE~\eqref{eq:SPDE} in 
the setting of this section. 
More precisely,
in the setting above the 
SPDE~\eqref{eq:SPDE} reduces 
to
\begin{equation}\label{eq:SPDE_reduced2}
  d X_t( x )
  =
  \Big[
    \kappa
    \Delta
    X_t( x )
    +
    f( x, X_t( x ) )
  \Big] dt
  +
  b( x, X_t( x ) ) \, dW_t( x )
\end{equation}
with $ X_{ t \, | \, \partial( 0,1 )^d } \equiv 0 $
and $ X_0( x ) = x_0( x ) $ for $ t \in [0,T] $ and $ x \in (0,1)^d $.
Moreover, we define a family 
$ 
  \beta^j \colon [0,T] \times 
  \Omega 
  \rightarrow \mathbb{R} 
$, 
$ 
  j \in 
  \{ 
    k \in \mathcal{J} 
    \colon
    \eta_k \neq 0 
  \} 
$, 
of independent standard Brownian motions by
\[
  \beta_t^j( \omega )
  :=
  \frac{ 1 }{ \sqrt{ \eta_j } }
  \left< g_j, W_t( \omega ) \right>_U
\]
for all $ \omega \in \Omega $, $ t \in [0,T] $ and all
$ j \in \mathcal{J} $ with $ \eta_j \neq 0 $. Using
this notation, the 
SPDE~\eqref{eq:SPDE_reduced2} 
can be written as
\begin{equation}\label{eq:SPDE19}
  d X_t( x )
  =
  \Big[
    \kappa
    \Delta
    X_t( x )
    +
    f( x, X_t( x ) )
  \Big] dt
  +
  \sum_{ 
    \substack{ j \in \mathcal{J} \\ 
    \eta_j \neq 0 } 
  }
  \Big[
    b( x, X_t( x ) ) 
    \sqrt{ \eta_j } \,
    g_j( x )
  \Big] d\beta_t^j
\end{equation}
with $ X_{ t \, | \, \partial( 0,1 )^d } \equiv 0 $ 
and $ X_0( x ) = x_0( x ) $ for $ t \in [0,T] $ and $ x \in (0,1)^d $.
The Milstein type
algorithm~\eqref{eq:scheme} 
applied to the SPDE~\eqref{eq:SPDE_reduced2} 
then
reduces to 
$ Y_0^{ N, M, K } = P_N( x_0 ) $ 
and
\begin{equation}
\label{eq:red_alg}
\begin{split}
  Y_{ m + 1 }^{ N, M, K }
& =
  P_N \,
  e^{ A \frac{T}{M} }
  \Bigg(
    Y_{ m }^{ N, M, K }
    +
    \frac{T}{M}
    \cdot
    f( \cdot, Y_{ m }^{ N, M, K } )
    +
    b( \cdot, Y_{ m }^{ N, M, K } )
    \cdot
    \Delta W_m^{ M, K }
\\ & \quad
    +
    \frac{1}{2}
    \left(
      \frac{\partial}{\partial y} b
    \right) \!
    ( \cdot, Y_{ m }^{ N, M, K } )
    \cdot
    b( \cdot, Y_{ m }^{ N, M, K } )
    \cdot
    \bigg(
      \left( \Delta W_m^{ M, K } \right)^2
      -
      \frac{T}{M}
      \sum_{ j \in \mathcal{J}_K }
      \eta_j( g_j )^2
    \bigg)
  \Bigg)
\end{split}
\end{equation}
for all $ m \in \{0, 1, \ldots, M-1 \} $ and all $ N, M, K \in \mathbb{N} $. Finally, Theorem~\ref{thm:mainresult} shows the existence of a real number $ C \in (0, \infty) $ such that
\begin{equation}
\label{eq:constant_C}
  \left(
    \mathbb{E}\!\left[
    \int_{ (0,1)^d }
    \big|
      X_{ T }(x) - Y_M^{ N, M, K }(x)
    \big|^2 \, dx
    \right]
  \right)^{ \! 1 / 2 }
\leq
  C\left(
    N^{ -2\gamma }
    +
    \bigg(
      \sup_{ 
        j \in 
        \mathcal{J}\setminus\mathcal{J}_K 
      }
      \eta_j
    \bigg)^{ \alpha }
    +
    M^{ 
      - \! \min( 2( \gamma - \beta ), \gamma ) 
    }
  \right)
\end{equation}
for all 
$ 
  N, M, K \in \mathbb{N} 
$.
We now illustrate 
estimate~\eqref{eq:constant_C} 
in the following 
three more concrete examples. 
We begin with the introductory 
example from Section~\ref{sec:intro} 
(see \eqref{eq:SPDE2} and \eqref{eq:SPDE_reduced}).

\subsection{A one-dimensional stochastic 
reaction diffusion equation}\label{sec:reacdiff}
In this subsection 
let $ d = 1 $, $ T = 1 $, 
$ \kappa = \frac{1}{100} $, 
let 
$ 
  x_0 \colon [0,1] \rightarrow \mathbb{R} 
$
be given by $ x_0( x ) = 0 $ 
for all $ x \in [0,1] $, 
let 
$  
  f, b \colon (0,1) \times 
  \mathbb{R} \rightarrow \mathbb{R} 
$ 
be given by 
$ f( x, y ) = 1 - y $ 
and $ b( x, y ) = \frac{1 - y}{1 + y^2} $
for all $ x \in (0,1) $, 
$ y \in \mathbb{R} $, 
let $ \mathcal{J} = \mathbb{N} $, 
let 
$ \mathcal{J}_K = \left\{ 1, 2, \dots, K
\right\} $ for all $ K \in \mathbb{N} $,
let $ \eta_j = \frac{1}{j^2} $ and 
let $ g_j = e_j $ for 
all $ j \in \mathbb{N} $. 
The SPDE~\eqref{eq:SPDE19} thus 
reduces to
\begin{equation}
\label{eq:SPDE_41}
  d X_t( x )
  =
  \left[
    \frac{1}{100}
    \frac{ \partial^2 
    }{ \partial x^2 } X_t( x )
    +
    1 - X_t( x )
  \right] dt
  +
  \sum_{ j = 1 }^{ \infty }
  \frac{ 1 - X_t( x )  
  }{ 
    1 + X_t( x )^2 
  }
  \frac{ \sqrt{ 2 } }{j}
  \sin( j \pi x ) \, d\beta_t^j
\end{equation}
with $ X_t( 0 ) = X_t( 1 ) = 0 $ and $ X_0( x ) = 0 $ for $ x \in (0,1) $ and $ t \in [0,1] $. The SPDE~\eqref{eq:SPDE_41} is nothing else than equation~\eqref{eq:SPDE_reduced} in the introduction. In order to apply Theorem~\ref{thm:mainresult} it remains to verify \eqref{condition_AA} and \eqref{condition_AAA}. Estimate~\eqref{condition_AA} is 
fulfilled for all 
$ \delta \in ( 0, \frac{1}{4} ) $ 
here due to Subsection~4.3 
in \cite{jr11}.
In order to establish \eqref{condition_AAA} several preparations are needed. More formally, let $ ( \tilde{\Omega}, \tilde{\mathcal{F}}, \tilde{\mathbb{P}} ) $ be a further probability space on which
a sequence 
$ 
  \chi_i 
  \colon
  \tilde{\Omega} \rightarrow   
  \mathbb{R} $, $ i \in \mathbb{N} 
$, 
of $ \tilde{\mathcal{F}} $/$ \mathcal{B}( \mathbb{R} ) $-measurable independent standard normal random variables is defined. 
Then we define 
the $ \tilde{\mathcal{F}} $/$ \mathcal{B}( U_0 ) $-measurable mappings 
$ 
  \chi^{ K, \vartheta } 
  \colon
  \tilde{\Omega} \rightarrow U_0 
$ 
by
\begin{equation}
 \chi^{ K, \vartheta }( \omega, x ) 
:= 
  \sum_{ i = 1 }^{ K } 
  \chi_i( \omega )
  \left( \lambda_i 
  \right)^{ -\vartheta } e_i( x ) 
\end{equation}
for all $ \omega \in \tilde{\Omega} $, 
$ x \in (0,1) $, 
$ K \in \mathbb{N} $ 
and all 
$ \vartheta \in (0,\frac{1}{2}) $. 
It will be essential to estimate 
$ 
  \tilde{\mathbb{E}}\!\left[
  | 
  \chi^{ K, \vartheta }( x ) 
  |^2 
  \right]
$
and 
$ 
  \tilde{\mathbb{E}}\!\left[
  | 
    \chi^{ K, \vartheta }( x ) - 
    \chi^{ K, \vartheta }( y ) 
  |^2 
  \right]
$ 
for $ x, y \in (0,1) $, 
$ K \in \mathbb{N} $ 
and $ \vartheta \in (0,\frac{1}{2}) $ 
in order to 
check \eqref{condition_AAA}. 
(Here 
$ \tilde{ \mathbb{E} }[ Z ]
:= \int_{ \tilde{\Omega} } 
Z( \tilde{\omega} ) \,
\tilde{\mathbb{P}}( \tilde{\omega} ) 
\in [0,\infty] $
for every
$ \tilde{\mathcal{F}} $/$
\mathcal{B}([0,\infty))$-measurable
mapping
$ 
  Z \colon 
  \tilde{\Omega} 
  \rightarrow
  [0,\infty) 
$.)
To this end note that
\begin{equation}
\label{eq:chi_A}
  \tilde{\mathbb{E}}\!\left[
  \left|
    \chi^{ K, \vartheta }( x )
  \right|^2
  \right]
=
  \sum_{ i = 1 }^{ K }
  \left( \lambda_i \right)^{ -2\vartheta }
  \left| e_i( x ) \right|^2
  \leq
  2
  \sum_{ i = 1 }^{ K }
  \left(
    \kappa \pi^2 i^2
  \right)^{ -2\vartheta }
\leq
  2\left( 1 + \kappa^{ -1 } \right)
  \left(
    \sum_{ i = 1 }^{ \infty }
    i^{ -4\vartheta }
  \right)
\end{equation}
for all 
$ x \in (0,1) $, 
$ K \in \mathbb{N} $ 
and all 
$ \vartheta \in ( 0, \frac{1}{2} ) 
$. 
Moreover, we have
\begin{align}
\label{eq:chi_B}
&
  \tilde{\mathbb{E}}\!\left[
  \left|
    \chi^{ K, \vartheta }(x)
    -
    \chi^{ K, \vartheta }(y)
  \right|^2
  \right]
=
  \sum_{ i = 1 }^{ K }
  \left( \lambda_i \right)^{ -2\vartheta }
  \left| e_i(x) - e_i( y ) \right|^2
\nonumber
\\&\leq
\nonumber
  \sum_{ i = 1 }^{ K }
  \left( \lambda_i \right)^{ -2\vartheta }
  \left| e_i(x) - e_i( y ) \right|^{ 2s }
  \left(
    \left| e_i(x) \right| + 
    \left| e_i( y ) \right|
  \right)^{ 2(1-s) }
\leq
  \sum_{ i = 1 }^{ K }
  \left( \lambda_i \right)^{ -2\vartheta }
  \left( 2 \pi^2 i^2 \right)^s
  8^{ (1-s) }
  \left| x - y \right|^{ 2s }
\\&\leq
  8\left(
    \sum_{ i = 1 }^{ K }
    \left( \kappa \pi^2 i^2 \right)^{ -2\vartheta }
    \left( \pi^2 i^2 \right)^s
  \right)
  \left| x - y \right|^{ 2s }
\leq
  \frac{ 8 
  }{ 
    \kappa^{ 2\vartheta } 
  }
  \left(
    \sum_{ i = 1 }^{ K }
    \left( 
      \pi i 
    \right)^{ 
      ( 2 s - 4 \vartheta ) 
    }
  \right)
  \left| x - y \right|^{ 2s }
\\&\leq
\nonumber
  3 \left( 1 + \kappa^{ -1 } \right)
  \left(
    \sum_{ i = 1 }^{ \infty }
    i^{ (2s - 4\vartheta) }
  \right)
  \left| x - y \right|^{ 2s }
\end{align}
for all $ x, y \in (0,1) $, $ K \in \mathbb{N} $, $ \vartheta \in ( \frac{s}{2} + \frac{1}{4}, \frac{1}{2} ) $ and all $ s \in (0, \frac{1}{2} ) $. 
We also use the notation
$$
  \left\|
    v
  \right\|_{ W^{ r, 2 } }
  :=
  \left(
    \int_0^1
    \left| v( x ) \right|^2 dx
    +
    \int_0^1
    \int_0^1
    \frac{ \left| v( x ) - v( y ) \right|^2 }{
           \left| x - y \right|^{ ( 1 + 2r ) } }
    \, dx \, dy
  \right)^{ \frac{1}{2} }
  \in
  [ 0, \infty ]
$$
for all $ \mathcal{B}( (0,1) ) $/$ \mathcal{B}( \mathbb{R} ) $-measurable mapping $ v : (0,1) \rightarrow \mathbb{R} $ and all
$ r \in ( 0, \infty ) $. 
Then we obtain
\begin{align*}
  \mathbb{E}\!\left[
  \left\|
    B( v ) \chi^{ K, \vartheta }
  \right\|^2_{ 
    W^{ r,2 }( (0,1), \mathbb{R} ) 
  }
  \right]
&=
  \int_0^1
  \mathbb{E}\!\left[
  \left|
    b( x, v(x) ) \cdot 
    \chi^{ K, \vartheta }( x )
  \right|^2
  \right] dx
\\&\quad+
  \int_0^1
  \int_0^1
  \frac{ 
    \mathbb{E}\!\left[
    \left|
      b( x, v(x) ) \cdot 
      \chi^{ K, \vartheta }( x )
      -
      b( y, v(y) ) \cdot 
      \chi^{ K, \vartheta }( y )
    \right|^2
    \right]
  }{
    \left| x - y \right|^{ (1+2r) }
  } \, dx \, dy
\\&\leq
  2 \int_0^1
  \left| b( x, v(x) ) \right|^2
  \mathbb{E}\!\left[
  \left| 
    \chi^{ K, \vartheta }( x )
  \right|^2 
  \right] dx
\\&\quad+
  2
  \int_0^1
  \int_0^1
  \frac{ 
    \left| b( x, v(x) ) \right|^2
    \mathbb{E}\!\left[
    \left| 
      \chi^{ K, \vartheta }( x ) - 
      \chi^{ K, \vartheta }( y ) 
    \right|^2
    \right]
  }{
    \left| x - y \right|^{ ( 1 + 2r ) }
  } \, dx \, dy
\\&\quad+
  2
  \int_0^1
  \int_0^1
  \frac{ 
    \left| b( x, v(x) ) - b( y, v(y) ) 
    \right|^2
    \mathbb{E}\!\left[
      \left| 
        \chi^{ K, \vartheta }( y )
      \right|^2
    \right]
  }{
    \left| x - y \right|^{ (1+2r) }
  } \, dx \, dy
\end{align*}
and using \eqref{eq:chi_A} 
and \eqref{eq:chi_B} shows
\begin{align*}
  \mathbb{E}\!\left[
  \left\|
    B( v ) \chi^{ K, \vartheta }
  \right\|^2_{ 
    W^{ r,2 }( (0,1), \mathbb{R} ) 
  }
  \right]
&\leq
  4 \left( 1 + \kappa^{ -1 } \right)
  \left(
    \sum_{ i=1 }^{ \infty }
    i^{ -4\vartheta }
  \right)
  \left\|
    b( \cdot, v )
  \right\|^2_{ 
    W^{ r,2 }( (0,1), \mathbb{R} ) 
  }  
\\&\quad+
  6 \left( 1 + \kappa^{ -1 } \right)
  \left(
    \sum_{ i=1 }^{ \infty }
    i^{ ( 2 s - 4 \vartheta ) }
  \right)
  \int_0^1
  \int_0^1
  \frac{ 
    \left| b( x, v(x) ) \right|^2
  }{
    \left| x - y \right|^{ (1 + 2r - 2s) }
  } \, dx \, dy
\\&\leq
  4 \left( 1 + \kappa^{ -1 } \right)
  \left(
    \sum_{ i=1 }^{ \infty }
    i^{ (2s-4\vartheta) }
  \right)
  \left\|
    b( \cdot, v )
  \right\|^2_{ 
    W^{ r,2 }( (0,1), \mathbb{R} ) 
  }  
\\&\quad+
  12 \left( 1 + \kappa^{ -1 } \right)
  \left(
    \sum_{ i=1 }^{ \infty }
    i^{ ( 2s -4\vartheta ) }
  \right)
  \left\|
    b( \cdot, v )
  \right\|^2_{ H }
  \int_0^1
  y^{ (2s -2r - 1) } \, dy
\\&\leq
  \frac{ 
    10 
    \left( 1 + \kappa^{ -1 } \right) 
  }{ (s-r) }
  \left(
    \sum_{ i=1 }^{ \infty }
    i^{ ( 2 s - 4 \vartheta ) 
    }
  \right)
  \left\|
    b( \cdot, v )
  \right\|^2_{ 
    W^{ r,2 }( (0,1), \mathbb{R} ) 
  }
\end{align*}
for all $ v \in H $, 
$ \vartheta \in ( \frac{s}{2} + \frac{1}{4}, \frac{1}{2} ) $, 
$ s \in (r, \frac{1}{2} ) $, 
$ K \in \mathbb{N} $ and 
all $ r \in (0,\frac{1}{2}) $. 
Therefore, inequality~(25) 
in Section~4 in \cite{jr11} gives
\begin{align}\label{eq:helpsec41}
\nonumber
  \left(
    \sup_{ K \in \mathbb{N} }
    \mathbb{E}\!\left[
    \left\|
      B( v ) \chi^{ K, \vartheta }
    \right\|^2_{ 
      W^{ r,2 }( (0,1), \mathbb{R} ) 
    }
    \right]
  \right)^{ \! \frac{1}{2} }
&\leq
  \frac{ 
    4 \left( 1 + \kappa^{ -1 } \right) 
  }{ 
    (s-r) 
  }
  \sqrt{
    \sum_{ i=1 }^{ \infty }
    i^{ ( 2s -4\vartheta ) }
  }
  \left\|
    b( \cdot, v )
  \right\|_{ 
    W^{ r,2 }( (0,1), \mathbb{R} ) 
  }
\\&\leq
  \frac{ 
    4 \left( 1 + \kappa^{ -1 } \right) 
  }{ 
    (s-r) 
  }
  \left(
    \sum_{ i=1 }^{ \infty }
    i^{ ( 2s -4\vartheta ) }
  \right)
  \frac{ 3 q C_{ \frac{r}{2} } 
  }{ ( 1- r ) }
  \left( 1 + 
    \left\| v 
    \right\|_{ H_{ \frac{r}{2} } 
    } 
  \right)
\\&\leq
\nonumber
  \frac{ 
    12 C_{ \frac{r}{2} } q 
    \left( 1 + \kappa^{ -1 } \right) 
  }{ 
    (s-r)^2 
  }
  \left(
    \sum_{ i=1 }^{ \infty }
    i^{ ( 2s -4\vartheta ) }
  \right)
  \left( 
    1 + 
    \left\| v 
    \right\|_{ H_{ \frac{r}{2} } } 
  \right)
  < \infty
\end{align}
for all 
$ \vartheta \in 
( \frac{s}{2} + \frac{1}{4}, 
\frac{1}{2} ) $, 
$ s \in (r, \frac{1}{2} ) $, 
$ 
  v \in H_{ \frac{r}{2} } 
$ 
and all 
$ 
  r \in (0, \frac{1}{2}) 
$. 
Moreover, we have
\begin{align*}
&
  \left\|
    \left( -A \right)^{ -\vartheta }
    B( v ) Q^{ -\alpha }
  \right\|_{ HS( U_0, H ) }
=
  \left\|
    \left( -A \right)^{ -\vartheta }
    B( v ) Q^{ ( \frac{1}{2} - \alpha ) }
  \right\|_{ HS( H ) }
=
  \left\|
    Q^{ ( \frac{1}{2} - \alpha ) }
    B( v ) \left( -A \right)^{ -\vartheta }
  \right\|_{ HS( H ) }
\\&=
  \left( \kappa \pi^2 \right)^{ ( \frac{1}{2} - \alpha ) }
  \left\|
    \left(
      \frac{Q}{ \kappa \pi^2 }
    \right)^{ ( \frac{1}{2} - \alpha ) }
    B( v ) \left( -A \right)^{ -\vartheta }
  \right\|_{ HS( H ) }
=
  \left( \kappa \pi^2 \right)^{ ( \frac{1}{2} - \alpha ) }
  \left\|
    B( v ) \left( -A \right)^{ -\vartheta }
  \right\|_{ 
    HS( H, 
      H_{ 
        ( \alpha - 1 / 2 ) 
      } 
    ) 
  }
\end{align*}
and using inequality~(20) 
in Section~4 in \cite{jr11} 
and estimate~\eqref{eq:helpsec41}
in this article then yields
\begin{align*}
  \left\|
    \left( -A \right)^{ -\vartheta }
    B( v ) Q^{ -\alpha }
  \right\|_{ HS( U_0, H ) }
&=
  \left( \kappa \pi^2 
  \right)^{ ( \frac{1}{2} - \alpha ) 
  }
  \left(
    \sup_{ K \in \mathbb{N} }
    \mathbb{E}\!\left[
    \left\|
      B( v ) \chi^{ K, \vartheta }
    \right\|_{ 
      H_{ 
        ( \alpha - 1 / 2 ) 
      } 
    }^2
    \right]
  \right)^{ \frac{1}{2} }
\\&\leq
  C_{ (\alpha - \frac{1}{2} ) }
  \left( 1 + \kappa^{ -1 } \right)
  \left(
    \sup_{ K \in \mathbb{N} }
    \mathbb{E}\left\|
      B( v ) \chi^{ K, \vartheta }
    \right\|_{ 
      W^{2\alpha-1, 2}( (0,1), \mathbb{R} ) 
    }^2
  \right)^{ \frac{1}{2} }
\\&\leq
  \frac{ 
    12 C_{ (\alpha - \frac{1}{2} ) }^2
    q \left( 1 + \kappa^{ -1 } \right)^2
  }{\left( s + 1 - 2 \alpha \right)^2}
  \left(
    \sum_{ i=1 }^{ \infty }
    i^{ ( 2 s - 4 \vartheta ) }
  \right)
  \left( 1 + 
    \left\| v 
    \right\|_{ H_{ 
        ( \alpha - \frac{1}{2} ) 
      } 
    } 
  \right)
<
  \infty
\end{align*}
for all $ \vartheta \in (\frac{s}{2} + \frac{1}{4}, \frac{1}{2} ) $,
$ s \in (2\alpha-1, \frac{1}{2} ) $, 
$ 
  v \in H_{ (\alpha-\frac{1}{2}) } 
$ 
and all 
$ 
  \alpha 
  \in ( \frac{1}{2}, \frac{3}{4} ) 
$. 
Therefore, 
estimate~\eqref{condition_AAA} 
is satisfied for 
all $ \alpha \in ( 0, \frac{3}{4} ) $
and all 
$ \gamma \in 
( \frac{1}{2}, \frac{3}{4} ) $. 
This finally shows that Assumptions~\ref{semigroup}-\ref{initial}
are fullfilled for 
the SPDE~\eqref{eq:SPDE_41}
for all $ \alpha \in ( 0, \frac{3}{4} ) $,
$ \beta = \frac{1}{5} $ 
and all 
$ \gamma \in 
( \frac{1}{2}, \frac{3}{4} ) $. 

Theorem~\ref{thm:mainresult}
therefore
yields the existence 
of real numbers 
$ C_r \in ( 0, \infty) $, 
$ r \in (0, \frac{3}{4} ) $, 
such that
\begin{equation}
\label{eq:est_constC}
  \left(
    \mathbb{E}\!\left[
    \int_0^1
    \left|
      X_T( x )
      -
      Y_{ M }^{ N, M, K }( x )
    \right|^2 dx
    \right]
  \right)^{ \! \frac{1}{2} }
  \leq
  C_r 
  \left( 
    N^{ ( r - \frac{ 3 }{ 2 } ) } + 
    K^{ ( r - \frac{ 3 }{ 2 } ) } + 
    M^{ ( r - \frac{ 3 }{ 4 } ) } 
  \right)
\end{equation}
for all 
$ N, M, K \in \mathbb{N} $
and all arbitrarily 
small $ r \in (0,\frac{3}{4}) $. 
In order to balance the 
error terms on the right 
hand side of \eqref{eq:est_constC} 
we choose $ N^2 = K^2 = M $ in \eqref{eq:est_constC} and 
obtain the existence of 
real numbers $ C_r \in ( 0, \infty ) $, 
$ r \in ( 0, \frac{3}{2} ) $, 
such that
\begin{equation}
\label{eq:est_constC2}
  \left(
    \mathbb{E}\!\left[
    \int_0^1
    \left|
      X_T( x )
      -
      Y_{ N^2 }^{ N, N^2, N }( x )
    \right|^2 dx
    \right]
  \right)^{ \frac{1}{2} }
  \leq
  C_r \cdot N^{ ( r - \frac{ 3 }{ 2 } ) }
\end{equation}
for all 
$ N \in \mathbb{N} $
and all arbitrarily 
small $ r \in (0,\frac{3}{2}) $. Estimate~\eqref{eq:est_constC2} 
is nothing else than inequality~\eqref{eq:realNumberCepsilon} in the introduction. We also refer to Figure~\ref{fig1} in the introduction 
for a numerical result 
illustrating \eqref{eq:est_constC2}.
\subsection{A one-dimensional
stochastic reaction
diffusion equation
with $ A Q \neq Q A $}
\label{sec:reacdiff2}
In Subsection~\ref{sec:reacdiff} 
we assumed that the eigenfunctions 
of the dominating linear 
operator $ A $ and 
of the covariance 
operator $ Q $ of the 
driving Wiener process 
$ 
  W \colon [0,T] \times \Omega 
  \rightarrow H 
$ 
of the 
SPDE~\eqref{eq:SPDE_reduced2} 
coincide and in particular, 
we assumed in 
Subsection~\ref{sec:reacdiff} that
\begin{equation}
\label{eq:cond_coi}
  A Q v = Q A v 
\end{equation}
holds for all $ v \in D( A ) $. 
However, our general setting 
in Section~\ref{sec:settings} 
does not need 
condition~\eqref{eq:cond_coi} 
to be fulfilled. To illustrate 
this fact we consider in this 
subsection an example in 
which \eqref{eq:cond_coi} fails 
to hold. More precisely, 
in 
this subsection let $ d = 1 $, 
$ T = 1 $, $ \kappa = \frac{1}{20} $, 
let 
$ 
  x_0 \colon [0,1] \rightarrow \mathbb{R} 
$ 
be given by $ x_0( x ) = 0 $ 
for all $ x \in [0,1] $, 
let 
$ 
  f, b \colon (0,1) \times 
  \mathbb{R} \rightarrow \mathbb{R} 
$ 
be given by 
$ 
  f( x, y ) = 1 - y 
$ 
and 
$ 
  b( x, y ) = \frac{ y }{ 1 + y^2 } 
$ 
for all $ x \in (0,1) $, 
$ y \in \mathbb{R} $,
let 
$ \mathcal{J} = \{ 0, 1, 2, \ldots \} $, 
let $ \mathcal{J}_K
= \left\{ 0, 1, \dots, K \right\} $
for all $ K \in \mathbb{N} $,
let $ \eta_0 = 0 $, 
$ \eta_j = \frac{1}{j^3} $ and 
let 
$ 
  g_j \colon (0,1) \rightarrow \mathbb{R} 
$ 
be given by $ g_0( x ) = 1 $, 
$ g_j( x ) = \sqrt{2} \cos( j \pi x ) $ 
for all $ x \in (0,1) $ and 
all $ j \in \mathbb{N} $. 
The SPDE~\eqref{eq:SPDE_reduced2} 
thus reduces to
\begin{equation}
\label{eq:SPDE_42}
  d X_t( x )
  =
  \left[
    \frac{1}{20}
    \frac{\partial^2}{\partial x^2} X_t( x )
    +
    1 - X_t( x )
  \right] dt
  +
  \frac{ X_t( x ) }{ 1 + X_t( x )^2 }
  \, dW_t( x )
\end{equation}
with $ X_t( 0 ) = X_t( 1 ) = 0 $ 
and $ X_0( x ) = 0 $ 
for $ x \in (0,1) $ 
and $ t \in [0, 1] $. 
Of course, the SPDE~\eqref{eq:SPDE_42} can 
also be written as
\begin{equation*}
  d X_t( x )
  =
  \left[
    \frac{1}{20}
    \frac{\partial^2}{\partial x^2} X_t( x )
    +
    1 - X_t( x )
  \right] dt
  +
  \sum_{ j = 1 }^{ \infty }
  \frac{ 
    X_t( x ) 
  }{ 
    1 + X_t( x )^2 
  }
  \frac{ \sqrt{ 2 } }{ j^{1.5} } 
  \cos( j \pi x ) \, d\beta_t^j
\end{equation*}
with $ X_t( 0 ) = X_t( 1 ) = 0 $ and $ X_0( x ) = 0 $ for $ x \in (0,1) $ and $ t \in [0, 1] $. Estimate~\eqref{condition_AA} is here fulfilled for all $ \delta \in (0, \frac{1}{2}) $ 
due to Subsection~4.2 in \cite{jr11}. 
Moreover, as in Subsection~\ref{sec:reacdiff} 
it can be shown that inequality~\eqref{condition_AAA} holds for all 
$ \alpha \in (0,\frac{2}{3}) $
and all 
$ \gamma \in (\frac{1}{2},1) $.
This finally shows 
that Assumptions~\ref{semigroup}-\ref{initial}
are fullfilled for 
the SPDE~\eqref{eq:SPDE_42}
for all 
$ \alpha \in ( 0, \frac{2}{3} ) $,
$ \beta = \frac{1}{5} $ 
and all 
$ \gamma \in (\frac{1}{2},1) $. 

Theorem~\ref{thm:mainresult}
therefore
yields the existence of real 
numbers $ C_r \in (0, \infty ) $, 
$ r \in (0,1) $, such that
\begin{equation}
\label{eq:C_r_estimate}
  \left(
    \mathbb{E}\!\left[
    \int_0^1
    \left| 
      X_T( x ) - 
      Y_M^{ N, M, K }( x ) 
    \right|^2 dx
    \right]
  \right)^{ \! \frac{1}{2} }
  \leq
  C_r 
  \left( 
    N^{ (r-2) } 
    + K^{ (r-2) } 
    + M^{ (r-1) } 
  \right)
\end{equation}
for all $ N, M, K \in \mathbb{N} $
and all arbitrarily 
small $ r \in (0,1) $. 
Choosing $ N^2 = K^2 = M $ in \eqref{eq:C_r_estimate} hence gives
the existence of real numbers 
$ 
  C_r \in (0,\infty) 
$, 
$ 
  r \in (0, 2 ) 
$, 
such that
\begin{equation}
\label{eq:C_r_estimate2}
  \left(
    \mathbb{E}\!\left[
    \int_0^1
    \left| 
      X_T( x ) - 
      Y_{N^2}^{ N, N^2, N }( x ) 
    \right|^2 dx
    \right]
  \right)^{ \! \frac{1}{2} }
  \leq
  C_r \cdot N^{ (r-2) }
\end{equation}
for all $ N \in \mathbb{N} $
and all arbitrarily 
small $ r \in (0,2) $. 
The Milstein
type approximation 
$ 
  Y_{ N^2 }^{ N, N^2, N } 
$ 
thus converges in the root mean 
square sense to $ X_T $ 
with order $ 2- $ as $ N $ goes to 
infinity. 
Since $ P_N( H ) \subset H $ 
is $ N $-dimensional 
and since $ N^2 $ time steps are used to simulate $ Y_{ N^2 }^{ N, N^2, N } $,
$ O( N^3 \log(N) ) $
computational operations and random
variables are needed to
simulate $ Y_{ N^2 }^{ N, N^2, N } $
here.
Combining the computational effort
$ O( N^3 \log(N) ) $
and the convergence 
order $ 2- $ 
in \eqref{eq:C_r_estimate2}
shows that 
the Milstein type
algorithm~\eqref{eq:red_alg}
with $ N^2 = K^2 = M $
needs about 
$ 
  O( \varepsilon^{ - \frac{3}{2} } ) 
$ 
computational
operations and random variables
to achieve a root mean square
precision $ \varepsilon > 0 $.

The linear implicit
Euler scheme combined with spectral
Galerkin methods
which we denote by
$ \mathcal{F} $/$ 
\mathcal{B}(H) $-measurable
mappings 
$ 
  Z^N_n \colon \Omega
  \rightarrow H 
$, 
$ 
  n \in \left\{ 0,1,\dots, N^4
  \right\} 
$, 
$ 
  N \in \mathbb{N} 
$,
is given by
$ 
  Z^N_0 = 0 
$ 
and
\begin{equation}
\label{eq:EulerAQ}
  Z^N_{ n + 1 }
=
  P_N 
  \left( 
    I - \frac{T}{N^4} A
  \right)^{ \! - 1 }
  \left(
    Z^N_n
    + \frac{T}{ N^4 } \cdot
    f( \cdot, Z^N_n )
    +
    b( \cdot, Z^N_n ) \cdot
    \Delta W^{N^4,N}_n
  \right)
\end{equation}
for all
$ 
  n \in \left\{ 0, 1, \dots, N^4 - 1
  \right\} 
$ 
and all 
$ 
  N \in \mathbb{N} 
$ 
here.

\begin{figure}[htp]
\begin{center}
\includegraphics[width=10cm]{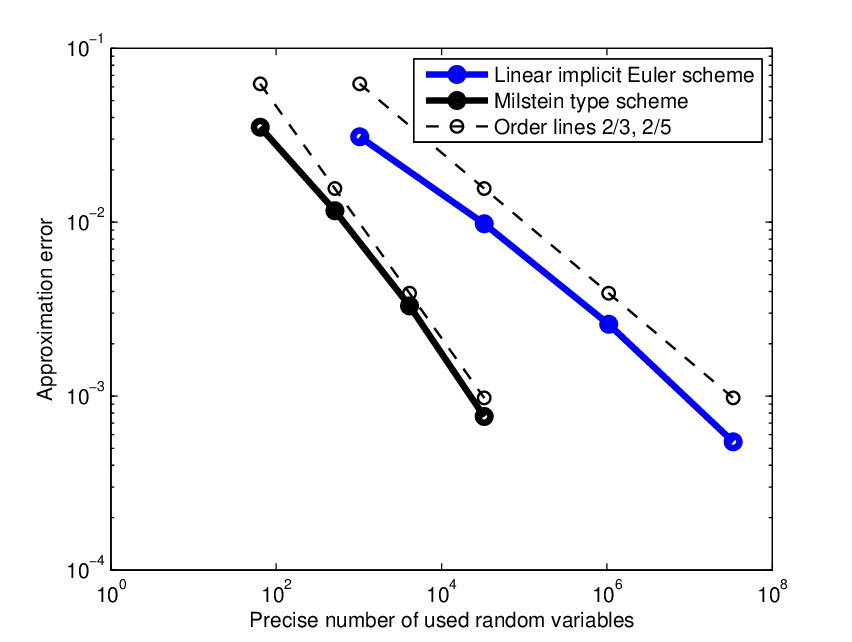}
\caption{\label{fig_heat}SPDE~\eqref{eq:SPDE_42}:
Root mean square
approximation error
$ 
  \big( 
    \mathbb{E}\big[ \|
      X_T - Z^N_{ N^4 }
    \|_H^2
    \big]
  \big)^{ 1/2 } 
$
of the linear implicit Euler
approximation $ Z^N_{N^4} $
(see \eqref{eq:EulerAQ}) 
and
root mean square
approximation error
$ 
  \big( 
    \mathbb{E}\big[
      \|
        X_T - Y^{ N, N^2, N }_{ N^2 }
      \|_H^2
    \big]
  \big)^{ 1 / 2 } 
$
of the Milstein type
approximation
$ Y^{ N, N^2, N }_{ N^2 } $
(see \eqref{eq:scheme}
and \eqref{eq:red_alg})
against
the precise 
number of independent 
standard normal random 
variables needed to 
compute the corresponding
approximation for 
$ 
  N \in \left\{ 4, 8, 16, 32 
  \right\} 
$.
}
\end{center}
\end{figure}

In Figure~\ref{fig_heat} 
the root mean square
approximation error
$ 
  \big( 
    \mathbb{E}\big[ \|
      X_T - Z^N_{ N^4 }
    \|_H^2
    \big]
  \big)^{ 1/2 } 
$
of the linear implicit Euler
approximation 
$ Z^N_{N^4} $
(see \eqref{eq:EulerAQ}) 
and
the root mean square
approximation error
$ 
  \big( 
    \mathbb{E}\big[
      \|
        X_T - Y^{ N, N^2, N }_{ N^2 }
      \|_H^2
    \big]
  \big)^{ 1 / 2 } 
$
of the Milstein type
approximation
$ Y^{ N, N^2, N }_{ N^2 } $
(see \eqref{eq:scheme}
and \eqref{eq:red_alg})
is plotted against
the precise 
number of independent 
standard normal random variables
needed to compute the corresponding
approximation for 
$ N \in \left\{ 4, 8, 16, 32 
\right\} $.

\subsection{A two-dimensional stochastic heat equation and splitting-up
approximations}
\label{sec:twoheat}

In this subsection 
the Milstein type 
algorithm~\eqref{eq:scheme}
is compared and related 
to certain 
splitting-up type
approximations in
the case of a two-dimensional
linear stochastic
heat equation with multiplicative
noise.
More formally, in this subsection let $ d = 2 $, $ T = 1 $, $ \kappa = \frac{1}{50} $,
let 
$ 
  x_0 \colon 
  [0,1]^2 \rightarrow \mathbb{R} 
$ 
be given by 
$ 
  x_0( x_1, x_2 ) 
  = 2 \sin( \pi x_1 ) \sin( \pi x_2 ) 
$ 
for all 
$ 
  x_1, x_2 \in [0,1] 
$, 
let 
$ 
  f, b \colon (0,1)^2 \times \mathbb{R} 
  \rightarrow \mathbb{R} 
$ 
be given by 
$ 
  f( x_1, x_2, y ) = 0 
$ 
and 
$ 
  b( x_1, x_2, y ) = y 
$ 
for all 
$ x_1, x_2 \in (0,1) $, 
$ y \in \mathbb{R} $, 
let $ \mathcal{J} = \mathbb{N}^2 $, 
let 
$ \mathcal{J}_K =
\left\{ 1, 2, \dots, K \right\}^2 $
for all $ K \in \mathbb{N} $,
let $ \eta_{ (j_1, j_2) } = ( j_1 + j_2 )^{-4} $ and let $ g_{ ( j_1, j_2 ) } = e_{ ( j_1, j_2 ) } $ for all $ j_1, j_2 \in \mathbb{N} $. 
The SPDE~\eqref{eq:SPDE_reduced2} 
thus reduces to
\begin{equation}
\label{eq:SPDE_2D}
  d X_t( x_1, x_2 )
  =
  \left[
    \frac{1}{50}
    \left(
      \frac{ \partial^2 }{ \partial x_1^2 }
      +
      \frac{ \partial^2 }{ \partial x_2^2 }
    \right) X_t( x_1, x_2 )
  \right] dt
  +
  X_t( x_1, x_2 ) \, dW_t( x_1, x_2 )
\end{equation}
with 
$ X_t|_{ \partial( 0,1 )^2 } \equiv 0 $ 
and
$ X_0( x_1, x_2 ) 
= 2 \sin( \pi x_1 ) \sin( \pi x_2 ) $ 
for $ x_1, x_2 \in (0,1) $ 
and $ t \in [0,1] $. 
In view of \eqref{eq:SPDE19} 
the SPDE~\eqref{eq:SPDE_2D} 
can also be written as
\begin{equation}
  d X_t( x_1, x_2 )
  =
  \left[
    \frac{1}{50}
    \left(
      \frac{ \partial^2 }{ \partial x_1^2 }
      +
      \frac{ \partial^2 }{ \partial x_2^2 }
    \right) X_t( x_1, x_2 )
  \right] dt
  +
  \sum_{ j_1, j_2 = 1 }^{ \infty }
  \frac{ X_t( x_1, x_2 ) 
  }{ ( j_1 + j_2 )^2 }
  2 \sin( j_1 \pi x_1 ) 
  \sin( j_2 \pi x_2 ) \, 
  d\beta_t^{ ( j_1, j_2 ) }
\end{equation}
with 
$ X_t|_{ \partial( 0,1 )^2 } \equiv 0 $ 
and
$ X_0( x_1, x_2 ) 
= 2 \sin( \pi x_1 ) \sin( \pi x_2 ) $ 
for $ x_1, x_2 \in (0,1) $ 
and $ t \in [0,1] $.

Due to Subsection~4.3 
in \cite{jr11} 
inequality~\eqref{condition_AA} 
holds for all 
$ \delta \in (0,\frac{1}{2}) $ here. 
In order to verify \eqref{condition_AAA} 
the notation
\begin{equation*}
  \left\| v \right\|_{ 
    L^{ \infty }( (0,1)^2, \mathbb{R} ) 
  }
  :=
  \inf \left\{ 
    R \in [0,\infty] 
    \colon
    \lambda \! \left( 
      \left\{ 
        x \in (0,1)^2 
        \colon 
        v( x ) > R 
      \right\}
    \right) = 0 
  \right\}
  \in [0,\infty]
\end{equation*}
is used for all 
$ 
  \mathcal{B}( (0,1)^2 ) 
$/$ 
  \mathcal{B}( \mathbb{R} ) 
$-measurable mappings 
$ 
  v \colon (0,1)^2 \rightarrow \mathbb{R} 
$ 
in this subsection. 
Then
\begin{equation}
\label{eq:two_v}
\begin{split}
  \left\|
    v
  \right\|_{ 
    L^{ \infty }( (0,1)^2, \mathbb{R} ) 
  }
& \leq
  \sum_{ i \in \mathbb{N}^2 }
  \left| \left< e_i, v \right>_H \right|
  \left\| e_i 
  \right\|_{ 
    C( (0,1)^2, \mathbb{R} ) 
  }
\\ & \leq
  2
  \left(
    \sum_{ i \in \mathbb{N}^2 }
    \left( \lambda_i \right)^{ -2r }
  \right)^{ \! \frac{1}{2} }
  \left(
    \sum_{ i \in \mathbb{N}^2 }
    \left( \lambda_i \right)^{ 2r }
    \left| \left< e_i, v \right>_H \right|^2
  \right)^{ \! \frac{1}{2} }
=
  2
  \left(
    \sum_{ i \in \mathbb{N}^2 }
    \left( \lambda_i \right)^{ -2r }
  \right)^{ \! \frac{1}{2} }
  \left\| v \right\|_{ H_r }
\end{split}
\end{equation}
for all 
$ 
  v \in H_r 
$ 
and all 
$ 
  r \in ( \frac{1}{2}, \infty ) 
$.
Moreover, we have
\begin{align*}
  \left\|
    \left( -A \right)^{ -\vartheta }
    B( v ) Q^{ -\alpha }
  \right\|_{ HS( U_0, H ) }
&=
  \left\|
    \left( -A \right)^{ -\vartheta }
    B( v ) Q^{ ( \frac{1}{2} - \alpha ) }
  \right\|_{ HS( H ) }
  =
  \left\|
    \left( -A \right)^{ -\vartheta }
    B( v ) Q^{ ( \frac{1}{2} - \alpha ) }
  \right\|_{ S_2( H ) }
\\&\leq
  \left\|
    \left( -A \right)^{ -\vartheta }
  \right\|_{ S_{ \frac{1}{\alpha} }( H ) }
  \left\|
    B( v )
  \right\|_{ L( H ) }
  \left\|
    Q^{ ( \frac{1}{2} - \alpha ) }
  \right\|_{ S_{ \frac{2}{( 1 - 2\alpha) } }( H ) }
\\&\leq 
  \left(
    \sum_{ i \in \mathbb{N}^2 }
    \left( \lambda_i \right)^{ -\frac{ \vartheta }{ \alpha } }
  \right)^{ \alpha }
  \left\|
    b( \cdot, v )
  \right\|_{ L^{ \infty }( (0,1)^2, \mathbb{R} ) }
  \left(
    \sum_{ j \in \mathcal{J} }
    \eta_j^{ ( \frac{1}{2} - \alpha) \frac{2}{(1-2\alpha)} }
  \right)^{ \frac{ (1-2\alpha)}{2} }
\end{align*}
and using \eqref{eq:two_v} shows
\begin{align*}
&
  \left\|
    \left( -A \right)^{ -\vartheta }
    B( v ) Q^{ -\alpha }
  \right\|_{ HS( U_0, H ) }
\leq 
  \left(
    \sum_{ i \in \mathbb{N}^2 }
    \left( \lambda_i \right)^{ -\frac{ \vartheta }{ \alpha } }
  \right)^{ \alpha }
  \left(
    q\left\|
      v
    \right\|_{ L^{ \infty }( (0,1)^2, \mathbb{R} ) }
    +
    q
  \right)
  \left( \text{Tr}( Q ) \right)^{ ( \frac{1}{2} - \alpha ) }
\\&\leq
  q \left( 1 + \text{Tr}( Q ) \right)
  \left(
    \sum_{ i \in \mathbb{N}^2 }
    \left( \lambda_i \right)^{ -\frac{ \vartheta }{ \alpha } }
  \right)^{ \alpha }
  \left(
    1
    +
    \left\|
      v
    \right\|_{ L^{ \infty }( (0,1)^2, \mathbb{R} ) }
  \right)
\\&\leq
  2 q \left( 1 + \text{Tr}( Q ) \right)
  \left(
    \sum_{ i \in \mathbb{N}^2 }
    \left( \lambda_i \right)^{ -\frac{ \vartheta }{ \alpha } }
  \right)^{ \alpha }
  \left(
    1
    +
    \sum_{ i \in \mathbb{N}^2 }
    \left( \lambda_i \right)^{ -2\gamma }
  \right)
  \left(
    1
    +
    \left\|
      v
    \right\|_{ H_{ \gamma } }
  \right)
  < \infty
\end{align*}
for all $ \vartheta \in (\alpha, \frac{1}{2}) $, $ \alpha \in (0,\frac{1}{2}) $, 
$ 
  v \in H_{ \gamma } 
$ 
and all 
$ \gamma \in (\frac{1}{2},1) $. Inequality~\eqref{condition_AAA} thus
holds for 
all $ \alpha \in (0,\frac{1}{2}) $
and all 
$ \gamma \in (\frac{1}{2},1) $
here.
This finally shows 
that
Assumptions~\ref{semigroup}-\ref{initial}
are fulfilled for 
the SPDE~\eqref{eq:SPDE_2D}
for all 
$ \alpha \in (0,\frac{1}{2}) $,
$ \beta = \frac{2}{5} $
and all 
$ \gamma \in (\frac{1}{2},1) $.

Theorem~\ref{thm:mainresult}
therefore yields the
existence of real numbers
$ C_r \in (0,\infty) $,
$ r \in (0,1) $,
such that
\begin{equation}
\label{eq:SPDE_2Dest}
  \left(
    \mathbb{E}\!\left[
    \int_0^1
    \int_0^1
    \big|
      X_T(x_1,x_2)
      -
      Y^{ N, M, K 
      }_{ M }(x_1,x_2)
    \big|^2 \,
    dx_1 \, dx_2
    \right]
  \right)^{ 
    \! 1 / 2 
  }
\leq
  C_r \left(
    N^{ \left( r - 2 \right) }
    +
    K^{ \left( r - 2 \right) }
    +
    M^{ \left( r - 1 \right) }
  \right)
\end{equation}
holds for all 
$ N, M, K \in \mathbb{N} $
and all arbitrarily 
small $ r \in (0,1) $.
In order to balance
the error terms
on the right hand side 
of \eqref{eq:SPDE_2Dest}
we choose $ M = N^2 = K^2 $
in \eqref{eq:SPDE_2Dest}
and obtain the existence
of real numbers 
$ C_r \in (0,\infty) $,
$ r \in (0,2) $, such that
\begin{equation}\label{eq:SPDE_2Dest2}
  \left(
    \mathbb{E}\!\left[
    \int_0^1
    \int_0^1
    \big|
      X_T(x_1,x_2)
      -
      Y^{ N, N^2, N 
      }_{ N^2 }(x_1,x_2)
    \big|^2
    \,
    dx_1 \, dx_2
    \right]
  \right)^{ 
    \! 1 / 2 
  }
\leq
  C_r \cdot
  N^{ (r - 2 ) }
\end{equation}
holds for 
all $ N \in \mathbb{N} $
and all arbitrarily 
small $ r \in (0,2) $.
The approximation
$ Y^{ N, N^2, N }_{ N^2 } $
thus converges in the
root mean square sense 
to $ X_T $ 
with order 
$ 2- $ as $N$ goes to
infinity.
The numerical approximations
$ 
  Y^{ N, N^2, N }_n
  \colon \Omega
  \rightarrow H 
$,
$ n \in 
\left\{ 0, 1, \dots, N^2
\right\} $, 
$ N \in \mathbb{N} $,
(see \eqref{eq:scheme}
and \eqref{eq:red_alg})
are here given by
$ Y^{ N, N^2, N }_0 = x_0 $
and
\begin{align}\label{eq:scheme2D}
&
  Y^{ N, N^2, N }_{ n + 1 }
= 
  P_N \, e^{ A \frac{T}{N^2} }
  \left(
    \Big[
      1 + 
      \Delta W^{ N^2, N }_n
      +
      \frac{1}{2}
      \left(
        \Delta W^{ N^2, N }_n
      \right)^2
      -
      \frac{T}{2 N^2}
      \sum_{ 
        j \in \mathcal{J}_K 
      }
      \eta_j ( g_j )^2
    \Big] \cdot
    Y^{ N, N^2, N }_n
  \right)
\end{align}
for all 
$ n \in \left\{ 0,1,\dots, N^2-1 \right\}$
and all $ N \in \mathbb{N} $.
Since $ P_N(H) \subset H $
is $ N^2 $-dimensional here
and since $ N^2 $ time steps
are used to simulate
$ Y^{ N, N^2, N }_{ N^2 } $,
$ O( N^4 \log(N) ) $
computational operations
and random variables are needed
to simulate 
$ Y^{ N, N^2, N }_{ N^2 } $.
Combining the computational effort
$ O( N^4 \log(N) ) $
and the convergence
order~$ 2 - $ in \eqref{eq:SPDE_2Dest2}
shows that the 
algorithm~\eqref{eq:scheme2D} 
in this article needs
about $ O( \varepsilon^{ - 2 } ) $ computational operations 
and random variables
to achieve a root mean square
precision $ \varepsilon > 0 $.

The linear implicit
Euler scheme combined with spectral
Galerkin methods
which we denote by
$ 
  \mathcal{F} 
$/$ 
  \mathcal{B}(H) 
$-measurable
mappings 
$ 
  Z^N_n \colon \Omega
  \rightarrow H 
$, 
$ n \in \left\{ 0,1,\dots, N^4
\right\} $, $ N \in \mathbb{N} $,
is given by
$ Z^N_0 = x_0 $ and
\begin{equation}\label{eq:Euler2D}
  Z^N_{ n + 1 }
=
  P_N \left( I - \frac{T}{N^4} A
  \right)^{ - 1 }
  \left(
    \Big[
      1 + \Delta W^{ N^4, N }_n
    \Big]
    \cdot
    Z^N_n
  \right)
\end{equation}
for all
$ n \in \left\{ 0, 1, \dots, N^4 - 1
\right\} $ and 
all $ N \in \mathbb{N} $ here.

Moreover, since the 
SPDE~\eqref{eq:SPDE_2D} is 
linear here, the splitting-up method 
can be used in 
order to solve \eqref{eq:SPDE_2D} approximatively. 
The idea of the splitting-up 
approach is to split the SPDE~\eqref{eq:SPDE_2D} into 
the explicit solvable 
subequations
\begin{equation}
\label{eq:split1}
  d \tilde{X}_t( x_1, x_2 )
  =
  \left[
    \frac{1}{50}
    \left(
      \frac{ \partial^2 }{ \partial x_1^2 }
      +
      \frac{ \partial^2 }{ \partial x_2^2 }
    \right)
    \tilde{X}_t( x_1, x_2 )
  \right] dt,
  \qquad
  \tilde{X}_t|_{ \partial( 0,1 )^2 } 
  \equiv 0
\end{equation}
and
\begin{equation}\label{eq:split2}
  d \tilde{\tilde{X}}_t( x_1, x_2 )
  =
  \tilde{\tilde{X}}_t( x_1, x_2 ) \,
  d W_t( x_1, x_2 )
\end{equation}
for $ t \in [0,1] $ 
and $ x_1, x_2  \in (0,1) $. 
For the solution processes
$ 
  \tilde{X}, 
  \tilde{\tilde{X}} 
  \colon 
  [0,T] \times \Omega
  \rightarrow H 
$
of \eqref{eq:split1}
and \eqref{eq:split2}
we obtain
$
  \tilde{X}_t = e^{ A t } 
  \tilde{X}_0
$
and
$
  \tilde{\tilde{X}}_t = 
  e^{ 
    \left(
      W_t -  
      \frac{t}{2}
      \sum_{ j \in \mathcal{J} }
      \eta_j ( g_j )^2
    \right)
  } 
  \cdot
  \tilde{\tilde{X}}_0
$
$
  \mathbb{P}
$-a.s.\ for 
all $ t \in [0,1] $.
This suggests the 
splitting-up approximation
$$
  X_t
\approx
  e^{ A t }
  \left(
  e^{ 
    \left(
      W_t -  
      \frac{t}{2}
      \sum_{ j \in \mathcal{J} }
      \eta_j ( g_j )^2
    \right)
  } 
  \cdot
  X_0
  \right)
$$
for $ t \in [0,1] $ where
$ 
  X \colon [0,T] \times
  \Omega \rightarrow H 
$
is the solution process 
of the SPDE~\eqref{eq:SPDE_2D}.
The resulting
splitting-up method which we denote 
by $ \mathcal{F} $/$ \mathcal{B}( H ) $-measurable mappings 
$ 
  \tilde{Z}_n^N \colon 
  \Omega \rightarrow H 
$, 
$ n \in \{ 0,1,\ldots,N^2 \} $, 
$ N \in \mathbb{N} $, is then 
given by 
$ \tilde{Z}_0^N = x_0 $ and
\begin{equation}
\label{eq:split}
  \tilde{Z}_{ n+1 }^N
=
  P_N \,
  e^{ A \frac{ T }{ N^2 } }
  \left(
    e^{ 
      \left(
        \Delta W^{N^2,N}_n - 
        \frac{T}{2 N^2}
        \sum_{ 
          j \in \mathcal{J}_K 
        }
        \eta_j ( g_j )^2
      \right)
    }
    \cdot
    \tilde{Z}_n^N
  \right)
\end{equation}
for all 
$ n \in 
\{ 0, 1, \ldots, N^2 - 1 \} $ 
and all $ N \in \mathbb{N} $.
Using the Talyor approximation
$ e^x \approx 1 + x 
+ \frac{x^2}{2} $
for all $ x \in \mathbb{R}$ 
then yields
\begin{align}
\label{eq:tayappr}
&
  e^{ 
    \left(
      \Delta W^{N^2,N}_n - 
      \frac{T}{2 N^2}
      \sum_{ 
        j \in \mathcal{J}_K 
      }
      \eta_j ( g_j )^2
    \right)
  }
\nonumber
\\&\approx
\nonumber
  1
  +
  \Delta W^{N^2,N}_n - 
  \frac{T}{2 N^2}
  \sum_{ 
    j \in \mathcal{J}_K 
  }
  \eta_j ( g_j )^2
  +
  \frac{1}{2}
  \left(
    \Delta W^{N^2,N}_n - 
    \frac{T}{2 N^2}
    \sum_{ 
      j \in \mathcal{J}_K 
    }
    \eta_j ( g_j )^2
  \right)^{\!2}
\\&\approx
  1
  +
  \Delta W^{N^2,N}_n 
  +
  \frac{1}{2}
  \left(
    \Delta W^{N^2,N}_n 
  \right)^2
  - 
  \frac{T}{2 N^2}
  \sum_{ 
    j \in \mathcal{J}_K 
  }
  \eta_j ( g_j )^2
\end{align}
for all $ n \in 
\left\{ 0, 1, \dots, N^2 - 1 \right\} $
and all $ N \in \mathbb{N} $.
Using approximation 
\eqref{eq:tayappr}
in \eqref{eq:split}
finally shows
$$
  \tilde{Z}^{ N }_{ n + 1 }\!
\approx
  P_N \, e^{ A \frac{T}{N^2} } \!
  \left(
    \Big[
      1 + 
      \Delta W^{ N^2, N }_n
      +
      \frac{1}{2}
      \left(
        \Delta W^{ N^2, N }_n
      \right)^2
      -
      \frac{T}{2 N^2}
      \sum_{ 
        j \in \mathcal{J}_K 
      }
      \eta_j ( g_j )^2
    \Big] \cdot
    \tilde{Z}^{ N }_n \!
  \right)
$$
for all
$ n \in \left\{ 0,1,\dots, N^2-1
\right\} $
and all $ N \in \mathbb{N} $
which is nothing else
than the recursion 
for $ Y^{N,N^2,N}_n $,
$ n \in \left\{0,1,\dots,N^2\right\}$,
$ N \in \mathbb{N} $,
in \eqref{eq:scheme2D}.
So, in the case of the 
linear SPDE~\eqref{eq:SPDE_2D}
an alternative way
for deriving
the Milstein type
algorithm~\eqref{eq:scheme2D}
is to apply an 
appropriate 
Taylor approximation
for the exponential function
(see \eqref{eq:tayappr}
for details)
to the splitting-up
approximation~\eqref{eq:split}.
More results on splitting-up
methods can be found in
A.\ Bensoussan, 
R.\ Glowinski and
A.\ Rascanu~\cite{bgr90,bgr92},
P.\ Florchinger and
F.\ Le Gland~\cite{fl91},
I.\ Gy{\"o}ngy and 
N.\ Krylov~\cite{gk03b,gk03a,gk05}
and K.\ Ito and 
B.\ L.\ Rozovskii~\cite{ir00}
and the references therein.

\begin{figure}
\begin{center}
\includegraphics[width=10cm]{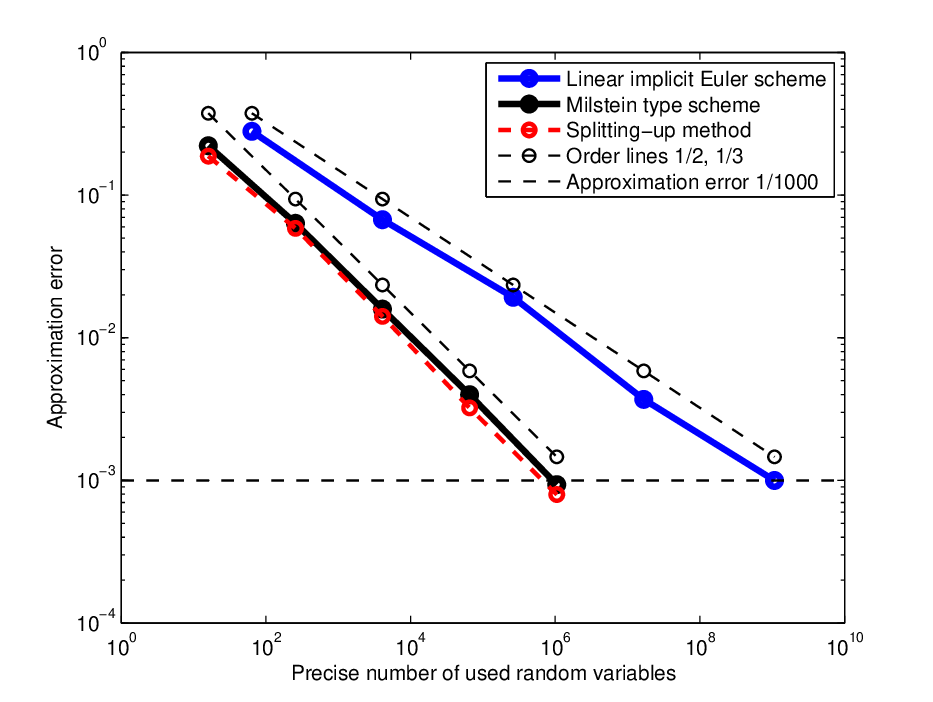}
\caption{SPDE~\eqref{eq:SPDE_2D}:
Root mean square
approximation error
$ 
  \big( 
    \mathbb{E}\big[ \|
      X_T - Z^N_{ N^4 }
    \|_H^2
    \big]
  \big)^{ 1 / 2 } 
$
of the linear implicit Euler
approximation $ Z^N_{N^4} $
(see \eqref{eq:Euler2D}),
root mean square
approximation error
of the linear implicit Euler
approximation $ Z^N_{N^4} $
(see \eqref{eq:Euler2D}),
the root mean square
approximation error
$ 
  \big( 
    \mathbb{E}\big[ \|
      X_T - Y^{ N, N^2, N }_{ N^2 }
    \|_H^2
    \big]
  \big)^{ 1 / 2 } 
$
of the Milstein type approximation
$ Y^{ N, N^2, N }_{ N^2 } $
(see \eqref{eq:scheme2D})
and 
root mean square
approximation error
$ 
  \big( 
    \mathbb{E}\big[ \|
      X_T - \tilde{Z}^N_{ N^2 }
    \|_H^2
    \big]
  \big)^{ 1 / 2 } 
$
of the
splitting-up 
approximation
$ \tilde{Z}^N_{ N^2 } $
(see \eqref{eq:split})
against
the precise 
number of independent 
standard normal random variables
needed to compute the corresponding
approximation for 
$ N \in \left\{ 2, 4, 8, 16, 32 \right\} $.
}
\label{fig_heat2d}
\end{center}
\end{figure}

In Figure~\ref{fig_heat2d}
the root mean square
approximation error
$ 
  \big( 
    \mathbb{E}\big[ \|
      X_T - Z^N_{ N^4 }
    \|_H^2
    \big]
  \big)^{ 1 / 2 } 
$
of the linear implicit Euler
approximation $ Z^N_{N^4} $
(see \eqref{eq:Euler2D}),
the root mean square
approximation error
$ 
  \big( 
    \mathbb{E}\big[ \|
      X_T - Y^{ N, N^2, N }_{ N^2 }
    \|_H^2
    \big]
  \big)^{ 1 / 2 } 
$
of the Milstein type
approximation
$ Y^{ N, N^2, N }_{ N^2 } $
(see \eqref{eq:scheme2D})
and 
the root mean square
approximation error
$ 
  \big( 
    \mathbb{E}\big[ \|
      X_T - \tilde{Z}^N_{ N^2 }
    \|_H^2
    \big]
  \big)^{ 1 / 2 } 
$
of the
splitting-up 
approximation
$ \tilde{Z}^N_{ N^2 } $
(see \eqref{eq:split})
is plotted against
the precise 
number of independent 
standard normal random variables
needed to compute the corresponding
approximation for 
$ N \in \left\{ 2, 4, 8, 16, 32 \right\} $:
It turns out 
that $ Z_{ 32^4 }^{ 32 } $ 
($ 32^6 = 1\,073\,741\,824 $
random variables)
in the case of the linear 
implicit Euler 
scheme~\eqref{eq:Euler2D},
that $ Y^{ 32, 32^2, 32 }_{ 32 } $
($ 32^4 = 1\,048\,576 $
random variables) 
in the case of the
algorithm~\eqref{eq:scheme2D} 
and that
$ \tilde{Z}^{ 32 }_{ 32^2 } $
($ 32^4 = 1\,048\,576 $
random variables) 
in the case of the
splitting-up method~\eqref{eq:split}
achieve a root mean square precision
$ \varepsilon = \frac{ 1 }{ 1000 } $ 
for the SPDE~\eqref{eq:SPDE_2D}.

\section{Proof of Theorem~\ref{thm:mainresult}}\label{secproofs}
Throughout this section 
the notation
\begin{equation}
  \left\|
    Z
  \right\|_{L^p (\Omega; E)}
  :=
  \Big(
    \mathbb{E}\big[
      \|
        Z
      \|_E^p
    \big]
  \Big)^{ 1 / p }
  \
  \in [0,\infty]
\end{equation}
is used 
for an  $ \mathbb{R} $-Banach
space $ \left( E, \left\| \cdot \right\|_E \right) $, 
an 
$ 
  \mathcal{F} 
$/$ 
  \mathcal{B}(E) 
$-measurable
mapping 
$ 
  Z \colon \Omega \rightarrow E 
$ 
and a real
number $ p \in [1,\infty) $.
We also use the following 
simple lemma 
(see, e.g., Theorem~37.5 in \cite{sy02}).
\begin{lemma}\label{lemmaRef4}
Let Assumptions~\ref{semigroup}-\ref{initial} in Section~\ref{sec:settings} be fulfilled. 
Then 
\begin{equation}
  \big\| 
    ( - t A )^r
    e^{ A t }
  \big\|_{ L(H) } 
  \leq 1
  \qquad
  \text{and}
  \qquad
  \big\|
    ( - t A )^{ - r }
    ( e^{ A t } - I )
  \big\|_{ L ( H ) }
  \leq
  1
\end{equation}
for all $ t \in (0,\infty) $
and all $ r \in [0,1] $.  
\end{lemma}
We now prove
Theorem~\ref{thm:mainresult}.
First of all, note that the exact solution of 
the SPDE~\eqref{eq:SPDE} satisfies
\begin{align}\label{eq:exactrep}
  X_{ m h } 
&=
  e^{ A m h } \xi
  +
  \int_0^{ m h }
  e^{ A(m h - s) }
  F\!\left( X_s \right)
  ds
  +
  \int_0^{ m h }
  e^{ A( m h - s) }
  B\!\left( X_s \right)
  dW_s
\\&=
\nonumber
  e^{ A m h } \xi
  +
  \sum_{ l=0 }^{ m-1 }
  \int_{ l h }^{ (l+1)h }
  e^{ A(m h - s) }
  F\!\left( X_s \right)
  ds
  +
  \sum_{ l=0 }^{ m-1 }
  \int_{ l h }^{ (l+1)h }
  e^{ A(m h - s) }
  B\!\left( X_s \right)
  dW_s
\end{align}
$ \mathbb{P} $-a.s.\ for all
$ m \in \{0,1,\ldots,M\} $ and 
all $ M \in \mathbb{N} $. Here and below 
$ h $ is the time stepsize 
$ h = h_M = \frac{T} {M} $ 
with $ M \in \mathbb{N} $.  
In particular, \eqref{eq:exactrep} shows
\begin{align}\label{eq:PNX}
  P_N\!\left( X_{ m h } \right)
&=
  e^{ A m h } P_N( \xi )
  +
  P_N\!\left(
  \sum_{ l=0 }^{ m-1 }
  \int_{ l h }^{ (l+1)h }
  e^{ A( m h - s) }
  F\!\left( X_s \right)
  ds
  \right)
\nonumber
\\&\quad+
  P_N\!\left(
  \sum_{ l=0 }^{ m-1 }
  \int_{ l h }^{ (l+1)h }
  e^{ A(m h-s) }
  B\!\left( X_s \right)
  dW_s \right)
\end{align}
$ \mathbb{P} $-a.s.\ for all
$ m \in \{0,1,\ldots,M\} $ and 
all $ N, M \in \mathbb{N} $.
In order to estimate the difference of the 
exact solution~\eqref{eq:exactrep} and 
the numerical solution~\eqref{eq:scheme} we rewrite
the numerical method~\eqref{eq:scheme} 
in some sense.
More precisely, the identity
\begin{equation}
\label{eq:centralob}
\begin{split}
&
   \frac{1}{2} 
    B'\!\left( Y_m^{N,M,K} \right) \!
    \bigg( \!
      B\!\left( Y_m^{N,M,K} \right) \!
      \Delta W_m^{M,K}
    \bigg) 
    \Delta W_m^{M,K}
    - \frac{T}{2 M} 
      \sum_{ 
        \substack{j \in \mathcal{J}_K 
        \\
       \eta_j \neq 0 } 
    }
        \eta_j
        B'\!\left( Y_m^{N,M,K} \right) \!
        \bigg( \!
          B\!\left( Y_m^{N,M,K} \right) 
          g_j
        \bigg) 
        g_j
\\ & =
   \int_{\frac{m T}{M}}^{\frac{(m+1)T}{M}}
    B'\!\left( Y_m^{N,M,K} \right) \!
    \bigg( \!
      \int_{\frac{m T}{M}}^s
      B\!\left( Y_m^{N,M,K} \right) dW_u^K \!
    \bigg) dW_s^K
\end{split}
\end{equation}
$ \mathbb{P} $-a.s.\ holds 
for all $ m \in \{0,1,\ldots,M-1\} $ 
and all $ N,M,K \in \mathbb{N} $. 
The proof of \eqref{eq:centralob} can be found
in Subsection~\ref{sec:centralob}.
Using \eqref{eq:centralob} shows that the numerical 
solution~\eqref{eq:scheme} fulfills
\begin{equation}
\label{eq:scheme2}
\begin{split}
  Y_{m+1}^{N,M,K} 
& = 
  P_N \, e^{ A \frac{ T }{ M } }
  \Bigg(
    Y_m^{N,M,K} 
    +
    \frac{T}{M} \cdot F\!\left( Y_m^{N,M,K} \right)
    +
    \int_{\frac{m T}{M}}^{\frac{(m+1)T}{M}}
    B\!\left( Y_m^{N,M,K} \right) dW_s^K
\\ & \quad
    +
    \int_{\frac{m T}{M}}^{\frac{(m+1)T}{M}}
    B'\!\left( Y_m^{N,M,K} \right) \!
    \bigg( \!
      \int_{\frac{m T}{M}}^s
      B\!\left( Y_m^{N,M,K} \right) dW_u^K \!
    \bigg) dW_s^K
  \Bigg)
\end{split}
\end{equation}
$ \mathbb{P} $-a.s.\ for 
all $ m \in \left\{ 0, 1, \dots, M-1 \right\} $
and all $ N, M, K \in \mathbb{N} $.
Therefore, the numerical solution~\eqref{eq:scheme} satisfies
\begin{align}\label{eq:defYNM}
  Y_{ m }^{ N,M,K }
&=
  e^{ A m h } P_N\!\left( \xi \right)
  +
  P_N\!\left(
    \sum_{ l=0 }^{ m-1 }
    \int_{ l h }^{ (l+1)h }
    e^{ A(m-l)h }
    F\!\left( Y_{ l }^{ N,M,K } \right)
    ds
  \right)
\nonumber
\\&\quad+
  P_N\!\left(
    \sum_{ l=0 }^{ m-1 }
    \int_{ l h }^{ (l+1)h }
    e^{ A(m-l)h }
    B\!\left( Y_{ l }^{ N,M,K } \right)
    dW_s^K
  \right)
\\&\quad+
\nonumber
  P_N\!\left(
    \sum_{ l=0 }^{ m-1 }
    \int_{ l h }^{ (l+1)h }
    e^{ A(m-l)h }
    B'\!\left( Y_{ l }^{ N,M,K } \right)
    \!\left(
      \int_{ l h }^{ s }
      B\!\left( Y_{ l }^{ N,M,K } \right)
      dW_u^K
    \right)
    dW_s^K
  \right)
\end{align}
$ \mathbb{P} $-a.s.\ for all 
$ m \in \{0,1,\ldots,M\} $ and 
all $ N, M, K \in \mathbb{N} $.
In order to estimate 
$ 
  \mathbb{E}\!\left[
    \| X_{ m h } - Y_m^{ N,M,K } 
    \|_H^2 
  \right]
$ 
for $ m \in \{0,1,\ldots,M \} $ 
and 
$ N, M, K \in \mathbb{N} $, 
we define the 
$ 
  \mathcal{F} 
$/$ 
  \mathcal{B}( H ) 
$-measurable 
mappings 
$ 
  Z_m^{N, M, K} 
  \colon \Omega \rightarrow H 
$, 
$ m \in \{ 0,1,\ldots,M \} $, 
$ N, M, K \in \mathbb{N} $, 
by
\begin{align}\label{eq:defYM}
  Z_{ m }^{ N, M, K }
&:=
  e^{ A m h } P_N\!\left( \xi \right)
  +
  P_N\!\left( 
  \sum_{ l=0 }^{ m-1 }
  \int_{ l h }^{ (l+1)h }
  e^{ A(m-l)h }
  F\!\left( X_{ l h } \right)
  ds 
  \right)
\nonumber
\\&\quad+
  P_N\!\left( 
  \sum_{ l=0 }^{ m-1 }
  \int_{ l h }^{ (l+1)h }
  e^{ A(m-l)h }
  B\!\left( X_{ l h } \right)
  dW_s^K
  \right)
\\&\quad+
\nonumber
  P_N\!\left( 
  \sum_{ l=0 }^{ m-1 }
  \int_{ l h }^{ (l+1)h }
  e^{ A(m-l)h }
  B'\!\left( X_{ l h } \right)
  \left(
    \int_{ l h }^{ s }
    B\!\left( X_{ l h } \right)
    dW_u^K
  \right)
  dW_s^K
  \right)
\end{align}
$ \mathbb{P} $-a.s.\ for 
all $ m \in \{0,1,\ldots,M\} $ 
and all $ N,M,K \in \mathbb{N} $.
The inequality
\begin{equation}\label{eq:squareest}
  \left( a_1 + \ldots + a_n
  \right)^2
\leq
  n \left( 
    \left( a_1 \right)^2 
    + \ldots +
    \left( a_n \right)^2 
  \right)  
\end{equation}
for all $ a_1, \dots, a_n \in \mathbb{R} $
and all $ n \in \mathbb{N} $ then shows
\begin{align}\label{eq:beginproof}
&
  \mathbb{E}\!\left[
    \|
      X_{ m h } - Y_m^{ N,M,K }
    \|_H^2
  \right]
\\ & \leq
\nonumber
  3 \cdot 
  \mathbb{E}\!\left[
  \|
    X_{ m h } - P_N\!\left( X_{ m h } \right)
  \|_H^2
  \right]
  +
  3 \cdot 
  \mathbb{E}\!\left[
  \|
    P_N\!\left( X_{ m h } \right) 
    - Z_m^{ N, M,K }
  \|_H^2
  \right]
  +
  3 \cdot 
  \mathbb{E}\!\left[
    \|
      Z_m^{ N, M,K } - Y_m^{ N,M,K }
    \|_H^2
  \right]
\end{align}
for all $ m \in 
\left\{ 0, 1, \dots, M \right\} $
and all $ N, M, K \in \mathbb{N} $.
In order to estimate
the expressions
$
  \mathbb{E}\!\left[
  \|
    X_{ m h } - P_N( X_{ m h } )
  \|_H^2
  \right]
$,
$
  \mathbb{E}\!\left[
  \|
    P_N\!\left( X_{ m h } \right) 
    - Z_m^{ N, M,K }
  \|_H^2
  \right]
$
and
$
  \mathbb{E}\!\left[
    \|
      Z_m^{ N, M,K } - Y_m^{ N,M,K }
    \|_H^2
  \right]
$
for $ m \in \left\{ 0,1,\dots, M \right\} $
and $ N, M, K \in \mathbb{N} $, 
a real number $ R \in (0,\infty) $ 
satisfying
$$
  \mathbb{E}\big[ 
    \| B(X_t) \|_{ 
      HS( U_0, H_{ \delta } ) 
    }^2
  \big]
  \leq R,
  \qquad
  \left\| F'\!\left( v \right) \right\|_{ L( H ) }
  \leq R,
  \qquad
  \left\| F''\!\left( v \right) 
  \right\|_{ L^{(2)}( H_{\beta}, H ) }
  \leq R ,
$$
$$
  \mathbb{E}\big[
    \| F\!\left( X_t \right) 
    \|_{ H }^2
  \big]
  \leq R ,
  \;\; 
  \left\| B'\!\left( v \right) 
  \right\|_{ L( H, HS( U_0, H ) ) } 
  \leq R ,
  \;\;
  \left\| 
    B''\!\left( v \right) 
  \right\|_{ 
    L^{(2)}( H_{\beta}, HS( U_0, H ) ) 
  } 
  \leq R ,
$$
$$
  \mathbb{E}\big[
    \| 
      ( - A )^{ \gamma } 
      X_t 
    \|_{ H }^2
  \big]
  =
  \mathbb{E}\big[
    \| 
      X_t 
    \|_{ H_{\gamma} }^2
  \big]
  \leq R, 
  \;\;
  \mathbb{E}\big[
    \| X_{t_2} - X_{t_1} 
    \|_{ H_{\beta} }^4
  \big]
  \leq 
  R \left| t_2 - t_1 
  \right|^{ 
    \min\left(  
      4 \left( \gamma - \beta \right),
      2 
    \right)
  } ,
$$
$$
  c +
  \frac{ 1 }{ \left( 1 - \gamma \right) }
  +
  \frac{ 1 }{ \left( 1 - 2 \vartheta \right) }
  +
  \frac{ 1 }{ \left( 1 - 2 \delta \right) }
  + T 
  + \left\|A^{-1} \right\|_{ L(H) } \leq R
$$
for all $ v \in H_{\beta} $
and all $ t, t_1, t_2 \in [0,T] $
is used throughout 
this proof.
Due to Assumptions~\ref{semigroup}-\ref{initial} in Section~\ref{sec:settings}
and Proposition~\ref{lem:existence} 
such a real number indeed exists.
For the spatial discretization error 
$ 
  \mathbb{E}\big[
  \| 
    X_{ m h }
    - P_N( 
      X_{ m h }
    )
  \|_H^2 
  \big] 
$
we then obtain
\begin{equation}
\label{eq:spatialerror}
\begin{split}
  \mathbb{E}\!\left[
  \left\| 
    X_{ m h }
    - P_N\!\left( 
      X_{ m h }
    \right)
  \right\|_H^2 
  \right]
&=
  \mathbb{E}\!\left[
  \left\| 
    \left( I - P_N \right) X_{ m h }
  \right\|_H^2 
  \right]
=
  \mathbb{E}\!\left[
  \left\| 
    \left( - A \right)^{ - \gamma }
    \left( I - P_N \right) 
    \left( - A \right)^{ \gamma }
    X_{ m h }
  \right\|_H^2 
  \right]
\\&\leq
  \left\|
    \left( - A \right)^{ - \gamma }
    \left( I - P_N \right) 
  \right\|_{ L(H) }^2
  \mathbb{E}\!\left[
    \left\| 
      X_{ m h }
    \right\|_{H_{\gamma}}^2
  \right] 
\leq
  R
  \left( r_N \right)^2
\end{split}
\end{equation}
for all $ m \in \left\{ 0, 1, \dots, M \right\} $
and all $ N, M \in \mathbb{N} $
where here and below the real numbers 
$   
  \left( r_N \right)_{ N \in \mathbb{N} } 
  \subset \mathbb{R} 
$ 
are given by
\begin{equation}
  r_N 
:= 
  \left\|
    \left( -A \right)^{ -\gamma }
    \left( I - P_N \right)
  \right\|_{ L( H ) }
=
  \left(
    \inf_{ 
      i \in \mathcal{I} \backslash
      \mathcal{I}_N
    }
    \lambda_i
  \right)^{ \! -\gamma }
\end{equation}
for all $ N \in \mathbb{N} $.
The rest of this proof is then divided into six parts.
In the first part (see Subsection~\ref{sec:driftproof})
we establish
\begin{equation}\label{eq:driftproof}
  \mathbb{E}\!\left[
  \left\|
    \sum_{ l=0 }^{ m-1 }
    \int_{ l h }^{ (l+1)h }
    \left(
      e^{ A(m h-s) }
      F\!\left( X_{ s } \right)
      -
      e^{ A(m-l)h }
      F\!\left( X_{ l h } \right)
    \right) ds
  \right\|_H^2
  \right]
  \leq
  \frac{ 36 R^{8} }{
  M^{ \min\left( 
    4 \left( \gamma - \beta \right) , 2\gamma 
  \right) } }
\end{equation}
for all $ m \in \{ 0, 1, \ldots, M \} $
and all $ M \in \mathbb{N} $.
We show
\begin{equation}\label{eq:pre_noiseI}
  \mathbb{E}\!\left[
  \left\|
    \sum_{ l = 0 }^{ m - 1 }
    \int_{ l h }^{ (l+1) h }
    e^{ A( m h - s ) }
    B\!\left( X_{ s } \right) 
    d\!\left( W_s - W_s^K \right)
  \right\|_H^2
  \right]
  \leq
  4 R^5
     \bigg(
      \sup_{ 
        j \in \mathcal{J} 
        \backslash \mathcal{J}_K 
      } 
      \eta_j 
      \bigg)^{ \! 2 \alpha }
\end{equation}
for all $ m \in \{ 0, 1, \ldots, M \} $
and all $ M, K \in \mathbb{N} $
in the second part 
(see Subsection~\ref{sec:pre_noiseI}) and
\begin{equation}\label{eq:noiseI}
  \mathbb{E}\!\left[
  \left\|
    \sum_{ l=0 }^{ m-1 }
    \int_{ l h }^{ (l+1)h }
    \left(
      e^{ A(m h-s) }
      -
      e^{ A(m-l)h }
    \right) 
    B\!\left( X_{ s } \right) dW_s^K
  \right\|_{ H }^2
  \right]
  \leq
  \frac{ 
    3 R^4  
  }{ 
    M^{ \left( 1 + 2\delta \right) } 
  }
\end{equation}
for all $ m \in \{ 0, 1, \ldots, M \} $
and all $ M, K \in \mathbb{N} $
in the third part 
(see Subsection~\ref{sec:noiseI}).
The fourth part (see Subsection~\ref{sec:noiseII})
gives
\begin{multline}\label{eq:noiseII}
  \mathbb{E}\!\left[
  \Bigg\|
    \sum_{ l=0 }^{ m-1 }
    \int_{ l h }^{ (l+1)h }
    e^{ A(m-l)h }
    \int_0^1
    B''\!\left( 
      X_{ l h } + r \left( X_{ s } - X_{ l h } 
    \right) \right)
    \left( 
      X_{ s } - X_{ l h }, X_{ s } - X_{ l h } 
    \right)
    \left( 1 - r \right) dr \, dW_s^K
  \Bigg\|_H^2
  \right]
\\
\leq
  \frac{ R^6 }{ 
    M^{ 
      \min\left( 
        4 \left( \gamma - \beta \right),
        2
      \right)    
    }
  }
\end{multline}
and in the fifth part (see 
Subsection~\ref{sec:noiseIII})
we obtain
\begin{multline}\label{eq:noiseIII}
  \mathbb{E}\!\left[
  \left\|
    \sum_{ l=0 }^{ m-1 }
    \int_{ l h }^{ (l+1)h }
    e^{ A(m-l)h }
    B'\!\left( X_{ l h } \right) 
    \left( 
      X_{ s } - X_{ l h }
      -
      \int_{ l h }^{ s }
      B\!\left( X_{ l h } \right) dW_u^K
    \right) dW_s^K
  \right\|_H^2  
  \right]
\\ \leq
  \frac{ 20 R^{ 13 } }{
    M^{ 
      \min\left(  
        4 \left( \gamma -
        \beta \right) ,
        2 \gamma
      \right) 
    }
  }
  +
  20 R^{ 11 }
    \bigg(
      \sup_{ j \in \mathcal{J} 
        \backslash \mathcal{J}_K 
      }
      \eta_j
    \bigg)^{ \! 2 \alpha }
\end{multline}
for all $ m \in \{ 0, 1, \ldots, M \} $
and all $ M, K \in \mathbb{N} $.
The
inequalities~\eqref{eq:driftproof}-\eqref{eq:noiseIII}
are used below to estimate
$
  \mathbb{E}\big[
  \|
    P_N\!\left( X_{ m h } \right) -
    Z_m^{N,M,K}
  \|_H^2
  \big]
$
for $ m \in \left\{ 0,1, \dots, M \right\} $
and $ N, M, K \in \mathbb{N} $ 
in \eqref{eq:beginproof}.
In the sixth part
(see Subsection~\ref{sec:lipproof})
we estimate
\begin{equation}
\label{eq:lipproof}
  \mathbb{E}\!\left[
  \left\|
    Z_m^{ N, M, K }
    -
    Y_m^{ N, M, K }
  \right\|_{ H }^2
  \right]
\leq
  \frac{ 9 R^4 }{ M } 
  \left(
    \sum_{ l=0 }^{ m - 1 }
    \mathbb{E}\!\left[
      \left\|
        X_{ l h }
        -
        Y_l^{ N,M,K }
      \right\|^2_H 
    \right]
  \right)
\end{equation}
for all $ m \in \{ 0, 1, \ldots, M \} $ and all
$ N, M, K \in \mathbb{N} $
by using the global Lipschitz continuity
of the coefficients 
$
  F \colon H_{\beta} \rightarrow H 
$
(see Assumption~\ref{drift}) 
and 
$
  B \colon H_{\beta} \rightarrow HS(U_0, H) 
$
(see Assumption~\ref{diffusion}).
Combining 
\eqref{eq:beginproof},
\eqref{eq:spatialerror}
and \eqref{eq:lipproof}
then yields 
\begin{equation}
\begin{split}
&
  \mathbb{E}\!\left[
    \left\|
      X_{ m h } - Y_m^{ N,M,K }
    \right\|_H^2
  \right]
\\ & \leq
  3 R \left( r_N \right)^2
  +
  3 \cdot 
  \mathbb{E}\!\left[
    \left\|
      P_N\!\left( X_{ m h } \right) -
      Z_m^{N,M,K}
    \right\|_H^2
  \right]
  +
  \frac{ 27 R^4 }{ M }
  \left(
    \sum_{ l=0 }^{ m-1 }
    \mathbb{E}\!\left[
    \left\|
      X_{ l h } - Y_l^{ N, M, K }
    \right\|_H^2
    \right]
  \right)
\end{split}
\end{equation}
for all $ m \in \{0,1,\ldots,M\} $ and 
all $ N, M, K \in \mathbb{N} $. Hence, \eqref{eq:PNX}, 
\eqref{eq:defYM},
\eqref{eq:squareest} 
and the fact
$ \left\| P_N(v) \right\|_H \leq 
\left\| v \right\|_H $ for all
$ v \in H $
show 
\begin{align*}
  \mathbb{E}\!\left[
  \left\|
    X_{ m h } - Y_m^{ N,M,K }
  \right\|_H^2
  \right]
&\leq
  9 \cdot 
  \mathbb{E}\!\left[
  \left\|
    \sum_{ l=0 }^{ m-1 }
    \int_{ l h }^{ (l+1)h }
    \left(
      e^{ A(m h-s) }
      F\!\left( X_s \right)
      -
      e^{ A(m-l)h }
      F\!\left( X_{l h} \right)
    \right) ds
  \right\|_H^2
  \right]
\\&\quad+
  9 \cdot \mathbb{E}\!\left[
  \left\|
    \sum_{ l=0 }^{ m-1 }
    \int_{ l h }^{ (l+1)h }
    e^{ A(m h-s) }
    B\!\left( X_s \right)
    d\!\left( W_s - W_s^K \right)
  \right\|_H^2
  \right]
\\&\quad+
  9 \cdot 
  \mathbb{E}\Bigg[
  \Bigg\|
    \sum_{ l=0 }^{ m-1 }
    \int_{ l h }^{ (l+1)h }
    \left( e^{ A(m h-s) }
    B\!\left( X_s \right) 
    -
    e^{ A(m-l)h }
    B\!\left( X_{l h} \right)
    \right) dW_s^K
\\&\quad\quad-
    \sum_{ l=0 }^{ m-1 }
    \int_{ l h }^{ (l+1)h }
    e^{ A(m-l)h }
    B'\!\left( X_{l h} \right) 
    \left(
      \int_{ l h }^s
      B\!\left( X_{l h} \right) dW_u^K
    \right) dW_s^K
  \Bigg\|_H^2
  \Bigg]
\\&\quad+
  3 R \left( r_N \right)^2
  +
  \frac{ 27 R^4 }{ M }
  \left(
    \sum_{ l=0 }^{ m-1 }
    \mathbb{E}\!\left[
    \left\|
      X_{ l h } - Y_l^{ N, M, K }
    \right\|_H^2
    \right]
  \right)
\end{align*}
for all $ m \in \{0,1,\ldots,M\} $ and 
all $ N, M, K \in \mathbb{N} $. 
Therefore,
\eqref{eq:driftproof} and \eqref{eq:pre_noiseI} yield
\begin{align*}
 &\mathbb{E}\!\left[
  \left\|
    X_{ m h } - Y_m^{ N,M,K }
  \right\|_H^2
  \right]
  \leq
    \frac{ 324 R^{ 8 } }{ 
      M^{ 
        \min\left( 
          4 \left( \gamma - \beta \right),
          2\gamma \right) 
      } 
    }
  +
  36 R^5
     \bigg(
      \sup_{ 
        j \in \mathcal{J} 
        \backslash \mathcal{J}_K 
      } 
      \eta_j 
      \bigg)^{ \! 2 \alpha }
\\ & \quad +
  18 \cdot 
  \mathbb{E}\!\left[
  \left\|
    \sum_{ l=0 }^{ m-1 }
    \int_{ l h }^{ (l+1)h }
    \left(
      e^{ A(m h-s) }
      -
      e^{ A(m-l)h }
    \right)
    B\!\left( X_{s} \right)
    dW_s^K
  \right\|_H^2
  \right]
\\&\quad+
  18 \cdot 
  \mathbb{E}\Bigg[ 
  \Bigg\|
    \sum_{ l=0 }^{ m-1 }
    \int_{ l h }^{ (l+1)h }
    e^{ A(m-l)h }
    \left(
      B\!\left( X_s \right)
      -
      B\!\left( X_{l h} \right) 
    \right) dW_s^K
\\&\qquad\qquad\quad-
    \sum_{ l=0 }^{ m-1 }
    \int_{ l h }^{ (l+1)h }
    e^{ A(m-l)h }
    B'\!\left( X_{l h} \right) 
    \left(
      \int_{ l h }^s
      B\!\left( X_{l h} \right) dW_u^K
    \right) dW_s^K
  \Bigg\|_H^2
  \Bigg]
\\&\quad+
  3 R \left( r_N \right)^2
  +
  \frac{ 27 R^4 }{ M }
  \left(
    \sum_{ l=0 }^{ m-1 }
    \mathbb{E}\!\left[
      \left\|
        X_{ l h } - Y_l^{ N, M, K }
      \right\|_H^2
    \right]
  \right)
\end{align*}
and \eqref{eq:noiseI} shows
\begin{align*}
  &
  \mathbb{E}\!\left[
  \left\|
    X_{ m h } - Y_m^{ N,M,K }
  \right\|_H^2
  \right]
  \leq
    \frac{
      324 R^{ 8 }
    }{ 
      M^{ 
        \min\left( 
          4 \left( \gamma - \beta \right) , 
          2 \gamma
        \right) 
      } 
    }
  +
  36 R^5
     \bigg(
      \sup_{ 
        j \in \mathcal{J} 
        \backslash \mathcal{J}_K 
      } 
      \eta_j 
    \bigg)^{ \! 2 \alpha }
  +
  \frac{54 R^4 }{ 
    M^{ \left( 1 + 2\delta \right) } }
\\&\quad+
  18 \cdot 
  \mathbb{E}\Bigg[
  \Bigg\|
    \sum_{ l=0 }^{ m-1 }
    \int_{ l h }^{ (l+1)h }
    e^{ A(m-l)h }
    \left(
      B\!\left( X_s \right)
      -
      B\!\left( X_{l h} \right) 
    \right) dW_s^K
\\&\qquad\qquad\quad-
    \sum_{ l=0 }^{ m-1 }
    \int_{ l h }^{ (l+1)h }
    e^{ A(m-l)h }
    B'\!\left( X_{l h} \right) 
    \left(
      \int_{ l h }^s
      B\!\left( X_{l h} \right) dW_u^K
    \right) dW_s^K
  \Bigg\|_H^2
  \Bigg]
\\&\quad+
  3 R \left( r_N \right)^2
  +
  \frac{ 27 R^4 }{ M }
  \left(
    \sum_{ l=0 }^{ m-1 }
    \mathbb{E}\!\left[
    \left\|
      X_{ l h } - Y_l^{ N, M, K }
    \right\|_H^2
    \right]
  \right)
\end{align*}
for all $ m \in \{0,1,\ldots,M\} $ and 
all $ N, M, K \in \mathbb{N} $. The fact
\begin{equation}
  B( X_s ) - B( X_{ l h } )
=
  B'( X_{ l h }) ( X_s - X_{ l h } )
  +
  \int^1_0 B''( X_{ l h } + r ( X_s - X_{ l h } ) )
  ( X_s - X_{ l h } , X_s - X_{ l h } ) (1-r) \, dr
\end{equation}
for all $ s \in [ l h, (l + 1) h ] $,
$ l \in \left\{ 0, 1, 
\dots, M-1 \right\} $
and all $ M \in \mathbb{N} $
then yields
\begin{align*}
  &
  \mathbb{E}\!\left[
  \left\|
    X_{ m h } - Y_m^{ N,M,K }
  \right\|_H^2
  \right]
\leq
  \frac{ 
    324 R^{ 8 } 
  }{ 
    M^{ 
      \min\left( 
        4 \left( \gamma - \beta \right) ,   
        2\gamma
      \right) 
    } 
  }
  +
  36 R^5
     \bigg(
      \sup_{ 
        j \in \mathcal{J} 
        \backslash \mathcal{J}_K 
      } 
      \eta_j 
    \bigg)^{ \! 2 \alpha }
  +
  \frac{ 
    54 R^4 
  }{ 
    M^{ \left( 1 + 2\delta \right) } 
  }
\\ & +
  36 \cdot 
  \mathbb{E}\!\left[
  \left\|
    \sum_{ l=0 }^{ m-1 }
    \int_{ l h }^{ (l+1)h } \!\!\!
    e^{ A(m-l)h }
    B'\!\left( X_{ l h} \right) \!
    \left(
      \!
      X_s - X_{ l h }
      - 
      \int^s_{ l h } 
      B\!\left( X_{l h} \right) dW_u^K
      \!
    \right) dW_s^K
  \right\|_H^2
  \right]
\\&+
  36 \cdot 
  \mathbb{E}\!\left[
  \Bigg\|
    \sum_{ l=0 }^{ m-1 }
    \int_{ l h }^{ (l+1)h } \!\!\!
    e^{ A(m-l)h } 
    \int_0^1
    B''\!\left( X_{ l h} + r \left(
        X_s - X_{ l h }
      \right) \right) \!
    \big(
      \!
      X_s - X_{ l h }, 
      X_s - X_{ l h }
      \!
    \big) \, (1-r) \, dr \, dW_s^K
  \Bigg\|_H^2
  \right]
\\&+
  3R \left( r_N \right)^2
  +
  \frac{ 27 R^4 }{ M }
  \left(
    \sum_{ l=0 }^{ m-1 }
    \mathbb{E}\!\left[ 
    \left\|
      X_{ l h } - Y_l^{ N, M, K }
    \right\|_H^2
    \right]
  \right)
\end{align*}
for all $ m \in \left\{ 0, 1, 
\dots, M-1 \right\} $
and all $ N, M, K \in \mathbb{N} $.
Therefore, \eqref{eq:noiseII}
and \eqref{eq:noiseIII} give
\begin{align*}
  &
  \mathbb{E}\!\left[
  \left\|
    X_{ m h } - Y_m^{ N,M,K }
  \right\|_H^2
  \right]
\leq
  \frac{ 324 R^{ 8 } }{
    M^{ \min\left( 
      4 \left( \gamma - \beta \right) ,
      2 \gamma 
      \right) 
    } 
  }
  +
  756 R^{ 11 }
     \bigg(
      \sup_{ 
        j \in \mathcal{J} \backslash 
        \mathcal{J}_K 
      } 
      \eta_j 
      \bigg)^{ \! 2 \alpha }
  +
  \frac{54 R^4 }{ 
    M^{ \left( 1 + 2\delta \right) } }
\\&\quad+
  \frac{ 
    720 R^{ 13 }
  }{
    M^{ \min\left(  
        4 \left( \gamma - \beta \right) ,
        2 \gamma
      \right) 
    }
  }
  +
  \frac{ 36 R^{ 6 } }{ 
    M^{ 
      \min\left(  
        4 \left( \gamma - \beta \right),
        2
      \right)
    } 
  }
  +
  3 R \left( r_N \right)^2
  +
  \frac{ 27 R^4 }{ M }
  \left(
    \sum_{ l=0 }^{ m-1 }
    \mathbb{E}\!\left[
    \left\|
      X_{ l h } - Y_l^{ N, M, K }
    \right\|_H^2
    \right]
  \right)
\end{align*}
and hence
\begin{align*}
 &\mathbb{E}\!\left[
  \left\|
    X_{ m h } - Y_m^{ N,M,K }
  \right\|_H^2
  \right]
  \leq
  \left(
    324 R^{ 8 } + 54 R^4
    + 720 R^{ 13 } + 36 R^6 
  \right)
  \frac{ 1 }{
    M^{ \min\left( 4 \left( \gamma - \beta \right) , 
      2\gamma \right) }
  }
\\&\quad
  +
  756 R^{11}
     \bigg(
      \sup_{ j \in \mathcal{J} \backslash \mathcal{J}_K } 
      \eta_j 
      \bigg)^{ \! 2 \alpha }
  +
  3 R \left( r_N \right)^2
  +
  \frac{ 27 R^4 }{ M }
  \left(
    \sum_{ l=0 }^{ m-1 }
    \mathbb{E}\!\left[
    \left\|
      X_{ l h } - Y_l^{ N, M, K }
    \right\|_H^2
    \right]
  \right)
\end{align*}
for all $ m \in \{0,1,\ldots,M\} $ and 
all $ N, M, K \in \mathbb{N} $.
Gronwall's lemma thus shows
\begin{align*}
&
  \mathbb{E}\!\left[
  \left\|
    X_{ m h } - Y_m^{ N,M,K }
  \right\|_H^2
  \right]
\\&\leq
  e^{ 27 R^4 }
  \Bigg(
    \frac{
      \left(
        324 R^{ 8 } + 54 R^4
        + 720 R^{ 13 } + 36 R^6 
      \right) 
    }{
      M^{ \min\left( 
      4 \left( \gamma - \beta \right), 
      2 \gamma
      \right) 
    }
  }
  +
  756 R^{ 11 }
     \bigg(
      \sup_{ j \in \mathcal{J} \backslash \mathcal{J}_K } 
      \eta_j 
      \bigg)^{ \! 2 \alpha }
  +
  3 R \left( r_N \right)^2
  \Bigg)
\\&\leq
  1134 R^{ 13 }
  e^{ 27 R^4 }
  \left(
    \left( r_N \right)^2
  +
     \bigg(
      \sup_{ j \in \mathcal{J} \backslash \mathcal{J}_K } 
      \eta_j 
      \bigg)^{ \! 2 \alpha }
  +
      M^{ - \min\left( 4 \left( \gamma - \beta \right),
        2 \gamma \right) }
  \right)
\end{align*}
for all $ m \in \{0,1,\ldots,M\} $ and 
all $ N, M, K \in \mathbb{N} $.
Finally, we obtain
\begin{align*}
   \left(
    \mathbb{E}\!\left[
    \left\|
      X_{ m h } - Y_m^{ N,M,K }
    \right\|_H^2
    \right]
  \right)^{ \frac{1}{2} }
&\leq
  34 R^7 e^{ 14 R^4 }
  \left(
    \bigg(
      \inf_{ 
        i \in \mathcal{I} \backslash
        \mathcal{I}_N
      }
      \lambda_i
    \bigg)^{ \! -\gamma }
    +
     \bigg(
      \sup_{ j \in \mathcal{J} \backslash \mathcal{J}_K } 
      \eta_j 
      \bigg)^{ \! \alpha }
    +
      M^{ - \min\left( 2 \left( \gamma - \beta \right), 
        \gamma \right) } 
  \right)
\\&\leq
  e^{ 20 R^4 }
  \left(
    \bigg(
      \inf_{ 
        i \in \mathcal{I} \backslash
        \mathcal{I}_N
      }
      \lambda_i
    \bigg)^{ \! -\gamma }
    +
     \bigg(
      \sup_{ j \in \mathcal{J} \backslash \mathcal{J}_K } 
      \eta_j 
      \bigg)^{ \! \alpha }
    +
      M^{ - \min\left( 2 \left( \gamma - \beta \right), 
        \gamma \right) } 
  \right)
\end{align*}
for all $ m \in \{0,1,\ldots,M\} $ and 
all $ N, M, K \in \mathbb{N} $.
\subsection{Temporal discretization 
error: Proof of \eqref{eq:driftproof}}
\label{sec:driftproof}
First of all, we have
\begin{align*}
& \left\|
    \sum_{ l=0 }^{ m-1 }
    \int_{ l h }^{ (l+1)h }
    \left(
      e^{ A(m h-s) }
      F\!\left( X_{ s } \right)
      -
      e^{ A(m-l)h }
      F\!\left( X_{ l h } \right)
   \right) ds
  \right\|_{ L^2( \Omega; H ) }
\\&\leq
  \sum_{ l=0 }^{ m-1 }
  \int_{ l h }^{ (l+1)h }
  \left\|
    \left(
      e^{ A(m h-s) }
      -
      e^{ A(m-l)h }
    \right)
    F\!\left( X_{ s } \right)
  \right\|_{ L^2( \Omega; H ) } ds
\\&\quad+
  \left\|
    \sum_{ l=0 }^{ m-1 }
    \int_{ l h }^{ (l+1)h }
    e^{ A(m-l)h }
    \left(
      F\!\left( X_{ s } \right)
      -
      F\!\left( X_{ l h } \right)
   \right) ds
  \right\|_{ L^2( \Omega; H ) }
\end{align*}
for all $ m \in \{0,1,\ldots,M\} $ and all $ M \in \mathbb{N} $. Using 
\begin{equation}
  F\!\left( X_s \right)
  -
  F\!\left( X_{l h} \right)
=
  F'\!\left( X_{l h} \right)
  \left( X_s - X_{l h} \right)
  +
  \int_0^1
  F''\!\left( X_{l h}  + r\left( X_s - X_{l h} \right) \right)
  \left( X_s - X_{l h}, X_s - X_{l h} \right)
  \left( 1 - r \right) dr
\end{equation}
for all 
$ s \in [ l h, (l+1)h ] $, 
$ l \in \{ 0,1,\ldots, M-1 \} $ 
and all 
$ M \in \mathbb{N} $ 
then shows
\begin{align*}
& \left\|
    \sum_{ l=0 }^{ m-1 }
    \int_{ l h }^{ (l+1)h }
    \left(
      e^{ A(m h-s) }
      F\!\left( X_{ s } \right)
      -
      e^{ A(m-l)h }
      F\!\left( X_{ l h } \right)
   \right) ds
  \right\|_{ L^2( \Omega; H ) }
\\&\leq
  \sum_{ l=0 }^{ m-1 }
  \int_{ l h }^{ (l+1)h }
  \left\|
    e^{ A(m h-s) }
    -
    e^{ A(m-l)h }
  \right\|_{ L( H ) }
  \left\|
    F\!\left( X_{ s } \right)
  \right\|_{ L^2( \Omega; H ) } ds
\\&+
  \left\|
    \sum_{ l=0 }^{ m-1 }
    \int_{ l h }^{ (l+1)h }
    e^{ A(m-l)h }
    F'\!\left( X_{ l h } \right)
    \left(
      X_{ s }
      -
      X_{ l h }
   \right) ds
  \right\|_{ L^2( \Omega; H ) }
\\&+
  \sum_{ l=0 }^{ m-1 }
  \int_{ l h }^{ (l+1)h }
  \int_0^1
  \left\| 
    F''\!\left( X_{ l h } + r \left( X_{ s } - X_{ l h } \right) \right)
    \left(
      X_{ s }
      -
      X_{ l h },
      X_{ s }
      -
      X_{ l h }
    \right)
  \right\|_{ L^2( \Omega; H ) }
  \left( 1 - r \right)dr \, ds
\end{align*}
and hence
\begin{align*}
& \left\|
    \sum_{ l=0 }^{ m-1 }
    \int_{ l h }^{ (l+1) h }
    \left(
      e^{ A(m h-s) }
      F\!\left( X_{ s } \right)
      -
      e^{ A(m-l)h }
      F\!\left( X_{ l h } \right)
   \right) ds
  \right\|_{ L^2( \Omega; H ) }
\\&\leq
  R \left( 
    2h +
    \sum_{ l=0 }^{ m-2 }
    \int_{ l h }^{ (l+1)h }
    \left\|
      e^{ A(m h-s) }
      -
      e^{ A(m-l)h }
    \right\|_{ L( H ) } ds
  \right)
\\&\quad+
  \left\|
    \sum_{ l=0 }^{ m-1 }
    \int_{ l h }^{ (l+1)h }
    e^{ A(m-l)h }
    F'\!\left( X_{ l h } \right)
    \left(
      X_{ s }
      -
      X_{ l h }
   \right) ds
  \right\|_{ L^2( \Omega; H ) }
\\&\quad+
  \sum_{ l=0 }^{ m-1 }
  \int_{ l h }^{ (l+1)h }
  \int_0^1
  \left\| 
    R
    \left\|
      X_{ s }
      -
      X_{ l h }
    \right\|_{ H_{\beta} }^2
  \right\|_{ L^2( \Omega; \mathbb{R} ) } 
  \left( 1 - r \right) dr \, ds
\end{align*}
for all $ m \in \{0,1,\ldots,M\} $ and all $ M \in \mathbb{N} $.
Therefore, we obtain
\begin{align*}
& \left\|
    \sum_{ l=0 }^{ m-1 }
    \int_{ l h }^{ (l+1)h }
    \left(
      e^{ A(m h-s) }
      F\!\left( X_{ s } \right)
      -
      e^{ A(m-l)h }
      F\!\left( X_{ l h } \right)
   \right) ds
  \right\|_{ L^2( \Omega; H ) }
\\&\leq
  R \left( 
    2h +
    \sum_{ l=0 }^{ m-2 }
    \int_{ l h }^{ (l+1)h }
    \left\| A e^{ A(m h-s) } \right\|_{L(H)}
    \left\| A^{-1} \left( e^{ A(s-l h) } - I \right) \right\|_{L(H)} ds
  \right)
\\&\quad+
  \left\|
    \sum_{ l=0 }^{ m-1 }
    \int_{ l h }^{ (l+1)h }
    e^{ A(m-l)h }
    F'\!\left( X_{ l h } \right)
    \left(
      X_{ s }
      -
      X_{ l h }
   \right) ds
  \right\|_{ L^2( \Omega; H ) }
\\&\quad+
  \frac{ R }{ 2 }
  \left(
    \sum_{ l=0 }^{ m-1 }
    \int_{ l h }^{ (l+1)h }
    \left( \mathbb{E}\left\|
      X_{ s }
      -
      X_{ l h }
    \right\|_{ H_{\beta} }^4 
    \right)^{ \frac{1}{2} } ds
  \right)
\end{align*}
and Lemma~\ref{lemmaRef4} gives
\begin{align*}
& \left\|
    \sum_{ l=0 }^{ m-1 }
    \int_{ l h }^{ (l+1)h }
    \left(
      e^{ A(m h-s) }
      F\!\left( X_{ s } \right)
      -
      e^{ A(m-l)h }
      F\!\left( X_{ l h } \right)
   \right) ds
  \right\|_{ L^2( \Omega; H ) }
\\&\leq
  R \left( 
    2h +
    \sum_{ l=0 }^{ m-2 }
    \int_{ l h }^{ (l+1)h }
    \frac{ \left( s-l h \right) }{ \left( m h - s \right) } \, ds
  \right)
\\&\quad+
  \sum_{ l=0 }^{ m-1 }
  \int_{ l h }^{ (l+1)h }
  \left\|
    e^{ A(m-l)h }
    F'\!\left( X_{ l h } \right)
    \left(
      \left( e^{ A(s - l h) } - I \right) 
      X_{ l h }
   \right) 
  \right\|_{ L^2( \Omega; H ) } ds
\\&\quad+
  \sum_{ l=0 }^{ m-1 }
  \int_{ l h }^{ (l+1)h }
  \left\|
    e^{ A(m-l)h }
    F'\!\left( X_{ l h } \right)
    \left(
      \int_{ l h }^s
      e^{ A(s - u) }
      F\!\left( X_{ u } \right) du
   \right)
  \right\|_{ L^2( \Omega; H ) } ds
\\&\quad+
  \left\|
    \sum_{ l=0 }^{ m-1 }
    \int_{ l h }^{ (l+1)h }
    e^{ A(m-l)h }
    F'\!\left( X_{ l h } \right)
    \left(
      \int_{ l h }^s
      e^{ A(s - u) }
      B\!\left( X_{ u } \right) dW_u
   \right) ds
  \right\|_{ L^2( \Omega; H ) }
\\&\quad+
  \frac{ R }{ 2 }
  \left(
    \sum_{ l=0 }^{ m-1 }
    \int_{ l h }^{ (l+1)h }
    \left(
      R \left( s-l h \right)^{ 
        \min\left( 
          4 \left( \gamma - \beta \right),
          2 
        \right)
      }
    \right)^{\frac{1}{2} } ds
  \right)
\end{align*}
for all $ m \in \{0,1,\ldots,M\} $ and all $ M \in \mathbb{N} $.
This shows
\begin{align*}
& \left\|
    \sum_{ l=0 }^{ m-1 }
    \int_{ l h }^{ (l+1)h }
    \left(
      e^{ A(m h-s) }
      F\!\left( X_{ s } \right)
      -
      e^{ A(m-l)h }
      F\!\left( X_{ l h } \right)
   \right) ds
  \right\|_{ L^2( \Omega; H ) }
\\&\leq
  R \left( 
    2h +
    \sum_{ l=0 }^{ m-2 }
    \int_{ l h }^{ (l+1)h }
    \frac{ \left( s-l h \right) }{ 
    ( m-l-1 ) h } \, ds
  \right)
\\&\quad+
  \sum_{ l=0 }^{ m-1 }
  \int_{ l h }^{ (l+1)h }
  \left\|
    F'\!\left( X_{ l h } \right)
    \left(
      \left( e^{ A(s - l h) } - I \right) 
      X_{ l h }
   \right) 
  \right\|_{ L^2( \Omega; H ) } ds
\\&\quad+
  \sum_{ l=0 }^{ m-1 }
  \int_{ l h }^{ (l+1)h }
  \left\|
    F'\!\left( X_{ l h } \right)
    \left(
      \int_{ l h }^s
      e^{ A(s - u) }
      F\!\left( X_{ u } \right) du
   \right)
  \right\|_{ L^2( \Omega; H ) } ds
\\&\quad+
  \left\{
    \sum_{ l=0 }^{ m-1 }
    \mathbb{E}\left\|
      \int_{ l h }^{ (l+1)h } \!\!\!
      e^{ A(m-l)h }
      F'\!\left( X_{ l h } \right)
      \left(
        \int_{ l h }^s \!\!
        e^{ A(s - u) }
        B\!\left( X_{ u } \right) dW_u
     \right) ds
    \right\|_{ H }^2
  \right\}^{ \frac{1}{2} }
\\&\quad+
  \frac{ R^2 }{ 2 }
  \left(
    \sum_{ l=0 }^{ m-1 }
    h^{ \left( 1 + 
      \min\left( 
        2 \left( \gamma - \beta \right),
        1
      \right)
    \right) }
  \right)
\end{align*}
and
\begin{align*}
& \left\|
    \sum_{ l=0 }^{ m-1 }
    \int_{ l h }^{ (l+1)h }
    \left(
      e^{ A(m h-s) }
      F\!\left( X_{ s } \right)
      -
      e^{ A(m-l)h }
      F\!\left( X_{ l h } \right)
   \right) ds
  \right\|_{ L^2( \Omega; H ) }
\\&\leq
  R \left( 
    2h +
    \sum_{ l=0 }^{ m-2 }
    \frac{ h }{ 2 ( m-l-1 ) }
  \right)
  +
  \frac{1}{2}
    R^2 T 
    h^{ 
      \min\left( 
        2 \left( \gamma - \beta \right), 1
      \right)
    }
\\&\quad+
  R \left(
    \sum_{ l=0 }^{ m-1 }
    \int_{ l h }^{ (l+1)h }
    \left\|
      \left( e^{ A(s - l h) } - I \right) 
      X_{ l h } 
    \right\|_{ L^2( \Omega; H ) } ds
  \right)
\\&\quad+
  R \left(
    \sum_{ l=0 }^{ m-1 }
    \int_{ l h }^{ (l+1)h }
    \left\|
      \int_{ l h }^s
      e^{ A(s - u) }
      F\!\left( X_{ u } \right) du
    \right\|_{ L^2( \Omega; H ) } ds
  \right)
\\&\quad+
  \sqrt{h}
  \left\{
    \sum_{ l=0 }^{ m-1 }
    \int_{ l h }^{ (l+1)h }
    \mathbb{E}\left\|
      F'\!\left( X_{ l h } \right)
      \left(
        \int_{ l h }^s
        e^{ A(s - u) }
        B\!\left( X_{ u } \right) dW_u
     \right)
    \right\|_{ H }^2 ds
  \right\}^{ \frac{1}{2} }
\end{align*}
for all $ m \in \{0,1,\ldots,M\} $ and all $ M \in \mathbb{N} $.
Hence, we have
\begin{align*}
& \left\|
    \sum_{ l=0 }^{ m-1 }
    \int_{ l h }^{ (l+1)h }
    \left(
      e^{ A(m h-s) }
      F\!\left( X_{ s } \right)
      -
      e^{ A(m-l)h }
      F\!\left( X_{ l h } \right)
   \right) ds
  \right\|_{ L^2( \Omega; H ) }
\\&\leq
  R \left( 
    2 h + \frac{ h }{ 2 }
    \left(
      \sum_{ l=1 }^{ m-1 }
      \frac{ 1 }{ l }
    \right)
  \right)
  +
  \frac{1}{2}
  R^3 h^{ 
    \min\left( 
      2 \left( \gamma - \beta \right), 1
    \right)
  }
\\&\quad+
  R \left(
    \sum_{ l=0 }^{ m-1 }
    \int_{ l h }^{ (l+1)h }
    \left\|
      \left( -A \right)^{ -\gamma }
      \left( e^{ A(s-l h) } - I \right)
    \right\|_{ L( H ) }
    \left\|
      \left( -A \right)^{ \gamma }
      X_{l h}
    \right\|_{ L^2( \Omega; H ) } ds
  \right)
\\&\quad+
  R \left(
    \sum_{ l=0 }^{ m-1 }
    \int_{ l h }^{ (l+1)h }
    \int_{ l h }^s
    \left\|
      F\!\left( X_{ u } \right) 
    \right\|_{ L^2( \Omega; H ) } du \, ds
  \right)
\\&\quad+
  R \sqrt{h}
  \left\{
    \sum_{ l=0 }^{ m-1 }
    \int_{ l h }^{ (l+1)h }
    \mathbb{E}\left\|
      \int_{ l h }^s
      e^{ A(s - u) }
      B\!\left( X_{ u } \right) dW_u
    \right\|_{ H  }^2 ds
  \right\}^{ \frac{1}{2} }
\end{align*}
and Lemma~\ref{lemmaRef4} 
shows
\begin{align*}
& \left\|
    \sum_{ l=0 }^{ m-1 }
    \int_{ l h }^{ (l+1)h }
    \left(
      e^{ A(m h-s) }
      F\!\left( X_{ s } \right)
      -
      e^{ A(m-l)h }
      F\!\left( X_{ l h } \right)
   \right) ds
  \right\|_{ L^2( \Omega; H ) }
\\&\leq
  R \left( 
    2 h + \frac{ h }{ 2 }
    \left(
      1 + \ln( M )
    \right)
  \right)
  +
  \frac{1}{2}
    R^4 
    M^{ 
      - \min\left(  
        2 \left( \gamma - \beta \right),
        1
      \right)
    }
\\&\quad+
  R^2 \left(
    \sum_{ l=0 }^{ m-1 }
    \int_{ l h }^{ (l+1)h }
    \left( s - l h \right)^{ \gamma } ds
  \right)
  +
  \frac{1}{2} R^2 M h^2
\\&\quad+
  R \sqrt{h}
  \left\{
    \sum_{ l=0 }^{ m-1 }
    \int_{ l h }^{ (l+1)h }
    \int_{ l h }^s
    \mathbb{E}\left\|
      e^{ A(s - u) }
      B\!\left( X_{ u } \right)
    \right\|_{ HS( U_0, H ) }^2 du \, ds
  \right\}^{ \frac{1}{2} }
\end{align*}
for all $ m \in \{0,1,\ldots,M\} $ and all $ M \in \mathbb{N} $.
Therefore, we obtain
\begin{align*}
& \left\|
    \sum_{ l=0 }^{ m-1 }
    \int_{ l h }^{ (l+1)h }
    \left(
      e^{ A(m h-s) }
      F\!\left( X_{ s } \right)
      -
      e^{ A(m-l)h }
      F\!\left( X_{ l h } \right)
   \right) ds
  \right\|_{ L^2( \Omega; H ) }
\\&\leq
  Rh \left( 
    \frac{5}{2} + \frac{ 1 }{ 2 } \ln( M )
  \right)
  +
  \frac{1}{2} R^4
    M^{ 
      - \min\left(  
        2 \left( \gamma - \beta \right),
        1
      \right)
    }
  +
  R^2 M h^{ (1 + \gamma) }
  +
  \frac{1}{2} R^3 h
\\&\quad+
  R \sqrt{h}
  \left\{
    \sum_{ l=0 }^{ m-1 }
    \int_{ l h }^{ (l+1)h }
    \int_{ l h }^s
    \mathbb{E}\left\|
      B\!\left( X_{ u } \right)
    \right\|_{ HS( U_0, H ) }^2 du \, ds
  \right\}^{ \frac{1}{2} }
\\&\leq
  R^2 M^{-1} \left( 
    \frac{5}{2} + \frac{ 1 }{ 2 } \ln( M )
  \right)
  +
  \frac{1}{2}
  R^4 M^{ 
    - \min\left(  
      2 \left( \gamma - \beta \right), 1
    \right)
  }
  +
  R^4 M^{ -\gamma }
  +
  \frac{1}{2} R^4 M^{ -1 }
\\&\quad+
  R \sqrt{h}
  \left\{
    \sum_{ l=0 }^{ m-1 }
    \int_{ l h }^{ (l+1)h }
    \int_{ l h }^s
    \left\| \left( - A \right)^{ - \delta }
    \right\|_{ L(H) }^2 
    \mathbb{E}\left\|
      \left( - A \right)^{ \delta } B\!\left( X_{ u } \right)
    \right\|_{ HS( U_0, H ) }^2 du \, ds
  \right\}^{ \frac{1}{2} }
\\&\leq
  \frac{5}{2} R^2 M^{-1} \left( 
     1 + \ln( M )
  \right)
  +
  \frac{1}{2}
  R^4 M^{ 
    - \min\left(  
      2 \left( \gamma - \beta \right), 1
    \right)
  }
  +
  R^4 M^{ -\gamma }
  +
  \frac{1}{2} R^4 M^{ -1 }
  +
  R^2 \sqrt{h}
  \left( \frac{1}{2} M h^2 \right)^{ \frac{1}{2} }
\end{align*}
for all $ m \in \{0,1,\ldots,M\} $ and all $ M \in \mathbb{N} $.
The estimate
$$
  1 + \ln(x)
  =
  1 + \int_1^x \frac{1}{s} \, ds
  \leq 
  1 + \int_1^x \frac{ 1 }{ s^{ ( 1 - r ) } } \, ds
  =
  1 + \frac{ \left( x^r - 1 \right) }{ r }
  =
  \frac{ x^r }{ r }
  - \frac{ \left( 1 - r \right) }{ r }
  \leq 
  \frac{ x^r }{ r }
$$
for all $ r \in (0,1] $ and 
all $ x \in [1,\infty) $ 
then shows
\begin{align*}
& \left\|
    \sum_{ l=0 }^{ m-1 }
    \int_{ l h }^{ (l+1)h }
    \left(
      e^{ A(m h-s) }
      F\!\left( X_{ s } \right)
      -
      e^{ A(m-l)h }
      F\!\left( X_{ l h } \right)
   \right) ds
  \right\|_{ L^2( \Omega; H ) }
\\&\leq
  \frac{5}{2} R^2 
  \frac{ M^{\left( 1 - \gamma \right)} 
  }{ M \left( 1 - \gamma \right) }
  +
  \frac{ R^4 }{ 
    2 
    M^{ 
      \min\left(  
        2 \left( \gamma - \beta \right), 1
      \right)
    } 
  }
  +
  \frac{ R^4 } { M^{ \gamma } }
  +
  \frac{ R^4 } { 2 M }
  +
  R^{ 2 } \sqrt{h}
  \left( T h \right)^{ \frac{1}{2} }
\\&\leq
  \frac{ 5 R^4 }{ 2 M^{ \gamma } }
  +
  \frac{ R^4 } { 2 
    M^{ 
      \min\left( 
        2 \left( \gamma - \beta \right), 1
      \right)
    } 
  }
  +
  \frac{ R^4 } { M^{ \gamma } }
  +
  \frac{ R^4 }{ 2 M }
  +
  \frac{ R^4 }{ M }
\end{align*}
and finally
\begin{align*}
& \left\|
    \sum_{ l=0 }^{ m-1 }
    \int_{ l h }^{ (l+1)h }
    \left(
      e^{ A(m h-s) }
      F\!\left( X_{ s } \right)
      -
      e^{ A(m-l)h }
      F\!\left( X_{ l h } \right)
   \right) ds
  \right\|_{ L^2( \Omega; H ) }
\\&\leq
  \left( \frac{5}{2} + \frac{1}{2} + 1 + \frac{1}{2}
   + 1 \right)
  \frac{ R^4  }{ M^{ \min\left( 
    2 \left( \gamma - \beta \right), \gamma 
  \right) } }
  \leq
  \frac{ 6 R^4  }{ M^{ \min\left( 
    2 \left( \gamma - \beta \right), \gamma 
  \right) } }
\end{align*}
for all $ m \in \{0,1,\ldots,M\} $ and all $ M \in \mathbb{N} $.
\subsection{Noise discretization error: 
Proof of \eqref{eq:pre_noiseI}}
\label{sec:pre_noiseI}
We have
\begin{align*}
&
  \mathbb{E}\left\|
    \int_s^t e^{ A ( t - u ) } B( X_u ) \,
    d\!\left( W_u - W^K_u \right)
  \right\|_H^2
=
  \mathbb{E}\left\|
    \sum_{ \substack{j \in \mathcal{J} \backslash \mathcal{J}_K  \\ 
    \eta_j \neq 0 } }
    \int_s^t e^{ A ( t - u ) } B( X_u ) g_j \,
    d\!\left< g_j, W_u \right>_U
  \right\|_H^2
\\&=
    \sum_{ \substack{ j \in \mathcal{J} \backslash \mathcal{J}_K  \\ \eta_j \neq 0 } }
    \eta_j 
    \int_s^t 
      \mathbb{E}
      \left\| 
        e^{ A ( t - u ) } B( X_u ) g_j
      \right\|_H^2
    du 
=
    \sum_{ \substack{j \in \mathcal{J} \backslash \mathcal{J}_K  \\ \eta_j \neq 0 } }
    \eta_j 
    \int_s^t 
      \mathbb{E}
      \left\| 
        e^{ A ( t - u ) } B( X_u ) Q^{ - \alpha }
        \left( Q^{ \alpha } g_j \right)
      \right\|_H^2
    du 
\\&=
    \sum_{ \substack{ j \in \mathcal{J} \backslash \mathcal{J}_K  \\ \eta_j \neq 0 } }
    \left( \eta_j \right)^{ \left( 1 + 2 \alpha 
    \right) } \left(
    \int_s^t 
      \mathbb{E}
      \left\| 
        e^{ A ( t - u ) } B( X_u ) Q^{ - \alpha }
        g_j
      \right\|_H^2
    du \right)
\end{align*}
for all $ s, t \in [0,T] $ with $ s \leq t $
and all $ K \in \mathbb{N} $.
This shows
\begin{align*}
&
  \mathbb{E}\left\|
    \int_s^t e^{ A ( t - u ) } B( X_u ) \,
    d\!\left( W_u - W^K_u \right)
  \right\|_H^2
\\&\leq
  \bigg( \sup_{ j \in \mathcal{J} 
    \backslash \mathcal{J}_K }
    \eta_j
  \bigg)^{ \! 2 \alpha }
  \left(
    \sum_{ \substack{ j \in \mathcal{J} \backslash \mathcal{J}_K  \\ \eta_j \neq 0 } }
    \eta_j 
    \left(
    \int_s^t 
      \mathbb{E}
      \left\| 
        e^{ A ( t - u ) } B( X_u ) Q^{ - \alpha }
        g_j
      \right\|_H^2
    du \right)
  \right)
\\&\leq
  \bigg( 
    \sup_{ 
      j \in \mathcal{J} \backslash \mathcal{J}_K 
    }
    \eta_j
  \bigg)^{ \! 2 \alpha }
  \left(
    \sum_{ j \in \mathcal{J} }
    \eta_j 
    \int_s^t 
      \mathbb{E}
      \left\| 
        e^{ A ( t - u ) } B( X_u ) Q^{ - \alpha }
        g_j
      \right\|_H^2
    du 
  \right)
\\&\leq
  \bigg( \sup_{ j \in \mathcal{J} 
    \backslash \mathcal{J}_K }
    \eta_j
  \bigg)^{ \! 2 \alpha }
  \left(
    \int_s^t 
      \mathbb{E}
      \left\| 
        e^{ A ( t - u ) } B( X_u ) Q^{ - \alpha }
      \right\|_{ HS( U_0, H ) }^2
    du 
  \right)
\\&\leq
  \bigg( \sup_{ j \in \mathcal{J} 
    \backslash \mathcal{J}_K }
    \eta_j
  \bigg)^{ \! 2 \alpha }
  \left(
    \int_s^t 
      \left( t - u \right)^{ - 2 \vartheta }
      \mathbb{E}
      \left\| 
        \left( - A \right)^{ - \vartheta } 
        B( X_u ) Q^{ - \alpha }
      \right\|_{ HS( U_0, H ) }^2
    du 
  \right)
\end{align*}
for all $ s, t \in [0,T] $ with $ s \leq t $
and all $ K \in \mathbb{N} $.
Therefore, we obtain
\begin{align*}
&
  \mathbb{E}\left\|
    \int_s^t e^{ A ( t - u ) } B( X_u ) \,
    d\!\left( W_u - W^K_u \right)
  \right\|_H^2
\\&\leq
  c^2
  \bigg( \sup_{ j \in \mathcal{J} 
    \backslash \mathcal{J}_K }
    \eta_j
  \bigg)^{ \! 2 \alpha }
  \left(
    \int_s^t 
      \left( t - u \right)^{ - 2 \vartheta }
      \mathbb{E}\left[
        \left( 1 + \left\| X_u 
          \right\|_{ H_{ \gamma } }
        \right)^2
      \right]
    du 
  \right)
\\&\leq
  2 c^2
  \bigg( \sup_{ j \in \mathcal{J} 
    \backslash \mathcal{J}_K }
    \eta_j
  \bigg)^{ \! 2 \alpha }
  \left(
    \int_s^t 
      \left( t - u \right)^{ - 2 \vartheta }
      \left( 1 + \mathbb{E} \left\| X_u 
        \right\|_{ H_{ \gamma } }^2
      \right)
    du 
  \right)
\\&\leq
  4 R^3
  \bigg( \sup_{ j \in \mathcal{J} 
    \backslash \mathcal{J}_K }
    \eta_j
  \bigg)^{ \! 2 \alpha }
  \left(
    \int_s^t 
      \left( t - u \right)^{ - 2 \vartheta }
    du 
  \right)
=
  4 R^3
  \bigg( \sup_{ j \in \mathcal{J} 
    \backslash \mathcal{J}_K }
    \eta_j
  \bigg)^{ \! 2 \alpha }
  \left(
    \int_0^{ \left( t - s \right) } 
      u^{ - 2 \vartheta } \,
    du 
  \right)
\end{align*}
for all $ s, t \in [0,T] $ with $ s \leq t $
and all $ K \in \mathbb{N} $.
Hence, we have
\begin{multline}
  \mathbb{E}\left\|
    \int_s^t e^{ A ( t - u ) } B( X_u ) \,
    d\!\left( W_u - W^K_u \right)
  \right\|_H^2
\leq
  4 R^3
  \bigg( \sup_{ j \in \mathcal{J} 
    \backslash \mathcal{J}_K }
    \eta_j
  \bigg)^{ \! 2 \alpha }
  \left[
    \frac{ 
      u^{ \left( 1 - 2 \vartheta \right) } 
    }{
      \left( 1 - 2 \vartheta \right)
    }
  \right]^{ u = \left( t - s \right) }_{ u = 0 }
\\
\leq
  4 R^4
  \bigg( \sup_{ j \in \mathcal{J} 
    \backslash \mathcal{J}_K }
    \eta_j
  \bigg)^{ \! 2 \alpha }
  \left( t - s \right)^{ \left( 1 - 2 \vartheta \right) } 
\label{eq:noisedis}
\leq
  4 R^5
  \bigg( \sup_{ j \in \mathcal{J} 
    \backslash \mathcal{J}_K }
    \eta_j
  \bigg)^{ \! 2 \alpha }
\end{multline}
for all $ s, t \in [0,T] $ with $ s \leq t $
and all $ K \in \mathbb{N} $.
In particular, we obtain
\begin{multline*}
  \mathbb{E}\left\|
    \sum_{ l = 0 }^{ m - 1 }
    \int_{ l h }^{ (l+1) h }
    e^{ A( m h - s ) }
    B\!\left( X_{ s } \right) 
    d\!\left( W_s - W_s^K \right)
  \right\|_H^2
\\=
  \mathbb{E}\left\|
    \int_{ 0 }^{ m h }
    e^{ A( m h - s ) }
    B\!\left( X_{ s } \right) 
    d\!\left( W_s - W_s^K \right)
  \right\|_H^2
\leq
  4 R^5 
  \bigg( \sup_{ j \in \mathcal{J} 
    \backslash \mathcal{J}_K }
    \eta_j
  \bigg)^{ \! 2 \alpha }
\end{multline*}
for all $ m \in \{0,1,\ldots,M\} $ and all $ M, K \in \mathbb{N} $ which shows
\eqref{eq:pre_noiseI}.
\subsection{Temporal discretization error: 
Proof of \eqref{eq:noiseI}}
\label{sec:noiseI}
We have
\begin{align*}
& \mathbb{E}\left\|
    \sum_{ l=0 }^{ m-1 }
    \int_{ l h }^{ (l+1)h }
    \left(
      e^{ A(m h-s) }
      -
      e^{ A(m-l)h }
    \right) 
    B\!\left( X_{ s } \right) dW_s^K
  \right\|_{ H }^2
\\&\leq 
    \sum_{ l=0 }^{ m-1 }
    \int_{ l h }^{ (l+1)h }
    \mathbb{E}\left\|
      \left(
        e^{ A(m h-s) }
        -
        e^{ A(m-l)h }
      \right) 
      B\!\left( X_{ s } \right)
    \right\|_{ HS( U_0, H ) }^2 ds
\\&\leq 
    \sum_{ l=0 }^{ m -1 } \!
    \int_{ l h }^{ (l+1)h } \!
    \left\|
      \left( -A \right)^{ -\delta } \!
      \left( 
        e^{ A(m h-s) }
        -
        e^{ A(m-l)h }
      \right)
    \right\|_{ L( H ) }^2 \!\!
    \mathbb{E}\left\|
      \left( -A \right)^{ \delta } \!\!
      B\!\left( X_{ s } \right)
    \right\|_{ HS( U_0, H ) }^2 \!\!\! ds
\end{align*}
and hence
\begin{align*}
& \mathbb{E}\left\|
    \sum_{ l=0 }^{ m-1 }
    \int_{ l h }^{ (l+1)h }
    \left(
      e^{ A(m h-s) }
      -
      e^{ A(m-l)h }
    \right) 
    B\!\left( X_{ s } \right) dW_s^K
  \right\|_{ H }^2
\\&\leq
  R \left(
    \sum_{ l=0 }^{ m -1 }
    \int_{ l h }^{ (l+1)h }
    \left\|
      \left( -A \right)^{ -\delta }
      \left( 
        e^{ A(m h-s) }
        -
        e^{ A(m-l)h }
      \right)
    \right\|_{ L( H ) }^2 ds
  \right)
\\&\leq
  R
    \int_{ (m-1)h }^{ m h }
    \left\|
      \left( -A \right)^{ -\delta }
      \left( 
        e^{ A(m h-s) }
        -
        e^{ Ah }
      \right)
    \right\|_{ L( H ) }^2 ds
\\&+
  R \left(
    \sum_{ l=0 }^{ m -2 } \!
    \int_{ l h }^{ (l+1)h }\!
    \left\|
      \left( -A \right)^{ -1 } \!
      \left( 
        e^{ A(s-l h) }
        -
        I
      \right)
    \right\|_{ L( H ) }^2 
    \left\|
      \left( -A \right)^{ (1 - \delta) } \!
      e^{ A(m h-s) }
    \right\|_{ L( H ) }^2 \!\! ds
  \right)
\end{align*}
for all $ m \in \{0,1,\ldots,M\} $ and all $ M, K \in \mathbb{N} $.
Therefore, we obtain
\begin{align*}
& \mathbb{E}\left\|
    \sum_{ l=0 }^{ m-1 }
    \int_{ l h }^{ (l+1)h }
    \left(
      e^{ A(m h-s) }
      -
      e^{ A(m-l)h }
    \right) 
    B\!\left( X_{ s } \right) dW_s^K
  \right\|_{ H }^2
\\&\leq
  R 
    \int_{ (m-1)h }^{ m h }
    \left\|
      \left( -A \right)^{ -\delta }
      \left( 
        e^{ A(s-(m-1)h) }
        -
        I
      \right)
    \right\|_{ L( H ) }^2 ds
\\&\quad+
  R h^2 \left(
    \sum_{ l=0 }^{ m -2 } \!
    \int_{ l h }^{ (l+1)h }\!
    \left\|
      \left( -A \right)^{ (1 - \delta) } \!
      e^{ A(m h-s) }
    \right\|_{ L( H ) }^2 ds
  \right)
\\&\leq
  R
    \int_{ (m-1)h }^{ m h }
    \left( s-(m-1)h \right)^{ 2\delta } ds
  +
  R h^2 \left(
    \sum_{ l=0 }^{ m -2 }
    \int_{ l h }^{ (l+1)h }
    \left( m h - s \right)^{ 2(\delta - 1) } ds
  \right)
\\&\leq
  R
  h^{ ( 1 + 2\delta) }
  +
  R h^3 \left(
    \sum_{ l=0 }^{ m -2 }
      \left( m - l - 1 \right)^{ 2(\delta - 1) } 
      h^{ 2(\delta - 1) }
  \right)
\end{align*}
for all $ m \in \{0,1,\ldots,M\} $ and all $ M, K \in \mathbb{N} $. 
This implies
\begin{align*}
& \mathbb{E}\left\|
    \sum_{ l=0 }^{ m-1 }
    \int_{ l h }^{ (l+1)h }
    \left(
      e^{ A(m h-s) }
      -
      e^{ A(m-l)h }
    \right) 
    B\!\left( X_{ s } \right) dW_s^K
  \right\|_{ H }^2
\leq
  R h^{ ( 1 + 2\delta) }
  +
  R h^{ (1 + 2\delta) } \left(
    \sum_{ l=1 }^{ m -1 }
   l^{ 2(\delta - 1) }
  \right)
\\&\leq
  R h^{ (1 + 2\delta) } \left(
   2 + \sum_{ l=2 }^{ \infty }
   l^{ 2(\delta - 1) }
  \right)
\leq
  R h^{ (1 + 2\delta) } \left(
   2 + \int_{ 1 }^{ \infty }
   s^{ 2(\delta - 1) } ds
  \right)
\\&\leq
  R h^{ (1 + 2\delta) } \left(
   2 + 
   \left[
     \frac{ s^{ (2\delta - 1) } }{ (2\delta - 1) }
   \right]_{ s = 1}^{ s = \infty }
  \right)
  =
  R h^{ (1 + 2\delta) } \left(
   2 + 
   \frac{ 1 }{ (1 - 2 \delta) }
  \right)
\leq
  3 R^2 h^{ (1 + 2\delta) }
  \leq
  \frac{ 
    3 R^4
  }{ 
    M^{ \left( 1 + 2 \delta \right) } 
  }
\end{align*}
for all $ m \in \{0,1,\ldots,M\} $ and all $ M, K \in \mathbb{N} $. 
\subsection{Temporal discretization error: 
Proof of \eqref{eq:noiseII}}
\label{sec:noiseII}
We have
\begin{align*}
& 
  \mathbb{E}
  \Bigg\|
  \sum_{ l=0 }^{ m-1 }
    \int_{ l h }^{ (l+1)h }
      e^{ A(m-l)h }
      \int_0^1
        B''\!\left( X_{ l h } 
          + r \left( X_{ s } - X_{ l h } 
          \right)  
        \right)
        \left( 
          X_{ s } - X_{ l h }, 
          X_{ s } - X_{ l h } 
        \right)
        \left( 1 - r \right) 
      dr 
    \, dW_s^K
  \Bigg\|_{ H }^2
\\&\leq
    \sum_{ l=0 }^{ m-1 } 
    \int_{ l h }^{ (l+1)h } 
    \int_0^1 
    \mathbb{E}\big\|
      B''\!\left( X_{ l h } 
        + r \left( X_{ s } - X_{ l h } 
        \right) 
      \right) \!
      \left( X_{ s } - X_{ l h }, 
        X_{ s } - X_{ l h } 
      \right) \!
    \big\|_{ HS( U_0, H ) }^2 dr \, ds 
\\&\leq
    \sum_{ l=0 }^{ m-1 } 
    \int_{ l h }^{ (l+1)h } 
    \mathbb{E}\!\left[
      \left(
        R
        \big\| X_{ s } - X_{ l h } 
        \big\|_{ H_{ \beta } }^2
      \right)^2
    \right] \! ds
  =
  R^2 \! \left(
    \sum_{ l=0 }^{ m-1 } 
    \int_{ l h }^{ (l+1)h } 
    \mathbb{E}\big\|
      X_{ s } - X_{ l h }
    \big\|_{ H_{\beta} }^4 \, ds 
  \right)
\\&\leq
  R^2 \left(
    \sum_{ l=0 }^{ m-1 }
    \int_{ l h }^{ (l+1)h }
      R \left( s-l h \right)^{ 
        \min\left( 
          4 \left( \gamma - \beta \right) ,
          2
        \right)
      }
    ds
  \right)
  \leq
  R^3 \left(
    \sum_{ l=0 }^{ m-1 }
    h^{ 
      \left( 1 
        + 
        \min\left( 4 \left( \gamma - \beta 
          \right) , 2 
        \right)
      \right) 
    }
  \right)
\\&\leq
  R^3 
  M
  h^{ \left( 1 + 
        \min\left( 
          4 \left( \gamma - \beta \right) ,
          2
        \right)    
  \right) }
  =
  R^3 T h^{ 
    \min\left( 
      4 \left( \gamma - \beta \right) , 2
    \right)
  }
\leq
  \frac{ R^6 }{ 
    M^{ 
      \min\left( 
        4 \left( \gamma - \beta \right) ,
        2
      \right)    
    }
  }
\end{align*}
for all $ m \in \{0,1,\ldots,M\} $ and all $ M,K \in \mathbb{N} $.
\subsection{Temporal discretization 
error: Proof of \eqref{eq:noiseIII}}
\label{sec:noiseIII}
In order to 
show \eqref{eq:noiseIII}, we first estimate
$$
  \mathbb{E}\left\|
     X_t - X_s - \int_s^t B\!\left( X_s \right) dW_u^K
  \right\|_H^2
$$
for all $ s, t \in [0,T] $ with $ s \leq t $ 
and all $ K \in \mathbb{N} $.
More precisely, we have
\begin{align*}
&
  \mathbb{E}\left\|
    X_{ t } - X_{ s }
    -
    \int_s^t
    B\!\left( X_{ s }  \right) dW_u^K
  \right\|_H^2
\leq 
  5 \cdot \mathbb{E}\left\|
    \left( e^{ A(t-s) } - I \right)
    X_{ s }
  \right\|_H^2
\\&\quad
  +
  5 \cdot \mathbb{E}\left\|
    \int_s^t
    e^{ A(t-u) }
    F\!\left( X_{ u }  \right) du
  \right\|_H^2
  +
  5 \cdot \mathbb{E}\left\|
    \int_s^t
      e^{ A(t-u) }
      B\!\left( X_{ u } \right) 
    d\!\left( W_u - W^K_u \right)
  \right\|_H^2
\\&\quad
  +
  5 \cdot \mathbb{E}\left\|
    \int_s^t \!\!\!
    \left( e^{ A(t-u) } - I \right)
    B\!\left( X_{ u }  \right) dW^K_u
  \right\|_H^2
  +
  5 \cdot \mathbb{E}\left\|
    \int_s^t
    \left( 
      B\!\left( X_{ u } \right)
      -
      B\!\left( X_{ s } \right)
    \right) dW^K_u
  \right\|_H^2
\end{align*}
and using \eqref{eq:noisedis} shows
\begin{align*}
& \mathbb{E}\left\|
    X_{ t } \!-\! X_{ s }
    \!-\!\!
    \int_s^t \!\!
    B\!\left( X_{ s }  \right) dW_u^K
  \right\|_H^2 \!\!\!
\leq 
  5 \left\| 
    \left( -A \right)^{ -\gamma }
    \left( e^{ A(t-s) } - I \right)
  \right\|_{ L( H ) }^2
  \mathbb{E}\!\left\| 
    \left( -A \right)^{ \gamma }
    X_s
  \right\|_H^2
\\&\quad+
  5 \left( t-s \right)
  \left(
  \int_s^t 
  \mathbb{E}\left\|
    e^{ A(t-u) }
    F\!\left( X_{ u } \right)
  \right\|_{ H }^2 du
  \right)
  +
  20 R^5 \bigg(
    \sup_{ j \in \mathcal{J} \backslash \mathcal{J}_K }
    \eta_j
  \bigg)^{ \! 2 \alpha }
\\&\quad+
  5 \left(
    \int_s^t 
    \mathbb{E}\left\| 
      \left( e^{ A(t-u) } - I \right) 
      B\!\left( X_{ u }  \right) 
    \right\|_{ HS( U_0, H ) }^2 du
  \right)
\\&\quad+
  5 \left(
    \int_s^t 
    \mathbb{E}\left\|
      B\!\left( X_{ u } \right)
      -
      B\!\left( X_{ s } \right)
    \right\|_{ HS( U_0, H ) }^2 du
  \right)
\end{align*}
for all $ s,t \in [0,T] $ with $ s \leq t $
and all $ K \in \mathbb{N} $.
This implies
\begin{align*}
& \mathbb{E}\left\|
    X_{ t } - X_{ s }
    -
    \int_s^t 
    B\!\left( X_{ s }  \right) dW_u^K 
  \right\|_H^2 
\\&\leq 
  5 R \left( t - s \right)^{ 2 \gamma }
  +
  5 \left( t -s \right)
  \left(
  \int_s^t 
  \mathbb{E}\left\|
    F\!\left( X_{ u } \right) 
  \right\|_{ H }^2 du
  \right)
  +
  20 R^5 \bigg(
    \sup_{ j \in \mathcal{J} \backslash \mathcal{J}_K }
    \eta_j
  \bigg)^{ \! 2 \alpha }
\\&\quad+
  5
  \left(
  \int_s^t
  \left\|
    \left( -A \right)^{ -\delta }
    \left( e^{ A( t-u ) } - I \right)
  \right\|_{ L( H ) }^2
  \mathbb{E}\left\|
    \left( -A \right)^{ \delta }
    B\!\left( X_{ u } \right) 
  \right\|_{ HS( U_0, H ) }^2 du
  \right)
\\&\quad+
  5 R^2 \left( 
    \int_s^t
    \mathbb{E}\left\|
      X_u - X_s
    \right\|_{ H }^2 du
  \right)
\\&\leq
  5 R \left( t -s \right)^{ 2 \gamma }
  +
  5 R \left( t - s \right)^2
  +
  20 R^5 \bigg(
    \sup_{ 
      j \in \mathcal{J} 
      \backslash \mathcal{J}_K 
    }
    \eta_j
  \bigg)^{ \! 2 \alpha }
\\&\quad+
  5
  \left(
  \int_s^t
    \left( t - u \right)^{ 2 \delta }
  \mathbb{E}\left\|
    \left( -A \right)^{ \delta }
    B\!\left( X_{ u } \right) 
  \right\|_{ HS( U_0, H ) }^2 du
  \right)
\\&\quad+
  5 R^2 \left( 
    \int_s^t
    \left\| \left( -A \right)^{ -\beta } \right\|_{ L( H ) }^2
    \mathbb{E}\left\|
      X_u - X_s
    \right\|_{ H_{\beta} }^2 du
  \right)
\end{align*}
for all $ s,t \in [0,T] $ with $ s \leq t $
and all $ K \in \mathbb{N} $.
Therefore, we obtain
\begin{align*}
&
  \mathbb{E}\left\|
    X_{ t } - X_{ s }
    -
    \int_s^t 
    B\!\left( X_{ s }  \right) dW_u^K 
  \right\|_H^2 
\leq
  5 R \left( t -s \right)^{ 2 \gamma }
  +
  5 R \left( t - s \right)^2
\\&\quad
  +
  20 R^5 \bigg(
    \sup_{ j \in \mathcal{J} \backslash \mathcal{J}_K }
    \eta_j
  \bigg)^{ \! 2 \alpha }
  +
  5 R
  \left(
  \int_s^t
    \left( t - u \right)^{ 2 \delta } du
  \right)
  +
  5 R^4 \left( 
    \int_s^t
    \mathbb{E}\left\|
      X_u - X_s
    \right\|_{ H_{\beta} }^2 du
  \right)
\\&\leq
  10 R^3 \left( t -s \right)^{ 2 \gamma }
  +
  20 R^5 \bigg(
    \sup_{ j \in \mathcal{J} \backslash \mathcal{J}_K }
    \eta_j
  \bigg)^{ \! 2 \alpha }
  +
  5 R
  \left( t - s \right)^{ \left( 1 + 2 \delta \right) }
  +
  5 R^4 \left( 
    \int_s^t
    \mathbb{E}\left\|
      X_u - X_s
    \right\|_{ H_{\beta} }^2 du
  \right)
\\&\leq
  15 R^3 \left( t -s \right)^{ 2 \gamma }
  +
  20 R^5 \bigg(
    \sup_{ j \in \mathcal{J} 
      \backslash \mathcal{J}_K }
    \eta_j
  \bigg)^{ \! 2 \alpha }
  +
  5 R^4 \left( 
    \int_s^t
    \left(
    \mathbb{E}\left\|
      X_u - X_s
    \right\|_{ 
      H_{\beta} 
    }^4 
  \right)^{  
    \frac{1}{2} 
  } du
  \right)
\\&\leq
  15 R^3 \left( t - s \right)^{ 2 \gamma }
  +
  20 R^5 \bigg(
    \sup_{ j \in \mathcal{J} \backslash \mathcal{J}_K }
    \eta_j
  \bigg)^{ \! 2 \alpha }
  +
  5 R^4 \left( 
    \int_s^t
    \left(
      R \left( u - s 
      \right)^{ 
        \min\left( 
          4 \left( \gamma - \beta \right) ,
          2
        \right) 
      }
    \right)^{ \frac{1}{2} } du
  \right)
\end{align*}
and hence
\begin{align}
\label{eq:estimate_1}
&
  \mathbb{E}\left\|
    X_{ t } - X_{ s }
    -
    \int_s^t 
    B\!\left( X_{ s }  \right) dW_u^K 
  \right\|_H^2 
\nonumber
\\&\leq
\nonumber
  15 R^3 \left( t - s \right)^{ 2 \gamma }
  +
  20 R^5 \bigg(
    \sup_{ j \in \mathcal{J} 
      \backslash \mathcal{J}_K }
    \eta_j
  \bigg)^{ \! 2 \alpha }
  +
  5 R^5 \left( 
    \int_s^t
      \left( u - s 
      \right)^{ 
        \min\left( 
          2 \left( \gamma - \beta \right),
          1
        \right) 
      }
    du
  \right)
\\&\leq
  15 R^3 \left( t - s \right)^{ 2 \gamma }
  +
  20 R^5 \bigg(
    \sup_{ j \in \mathcal{J} 
      \backslash \mathcal{J}_K }
    \eta_j
  \bigg)^{ \! 2 \alpha }
  +
  5 R^5 \left( t - s \right)^{
    \left( 1 + 
      \min\left( 
        2 \left( \gamma - \beta \right),
        1 
      \right)
    \right)
  }
\\&\leq
  15 R^3 \left( t - s \right)^{ 2 \gamma }
  +
  20 R^5 \bigg(
    \sup_{ j \in \mathcal{J} 
      \backslash \mathcal{J}_K 
    }
    \eta_j
  \bigg)^{ \! 2 \alpha }
  +
  5 R^6 \left( t - s \right)^{
    \min\left(  
      4 \left( \gamma - \beta \right), 
      2
    \right)
  }
\nonumber
\\&\leq
\nonumber
  20 R^8 \left( t -s \right)^{ 
    \min\left( 
      4 ( \gamma - \beta ) ,
      2 \gamma
    \right) 
  }
  +
  20 R^5 \bigg(
    \sup_{ j \in \mathcal{J} \backslash \mathcal{J}_K }
    \eta_j
  \bigg)^{ \! 2 \alpha }
\end{align}
for all $ s,t \in [0,T] $ with $ s \leq t $
and all $ K \in \mathbb{N} $.
Now we prove \eqref{eq:noiseIII}. 
To this end we note that
\begin{align*}
& \mathbb{E}\left\|
    \sum_{ l=0 }^{ m-1 } 
    \int_{ l h }^{ (l+1)h }
    e^{ A(m-l)h }
    B'\!\left( X_{ l h } \right) 
    \left( 
      X_{ s } - X_{ l h }
      -
      \int_{ l h }^{ s }
      B\!\left( X_{ l h } \right) dW^K_u 
      \right) dW^K_s
  \right\|_{ H }^2
\\&\leq
    \sum_{ l=0 }^{ m-1 }
    \int_{ l h }^{ (l+1)h }
    \mathbb{E}\left\|
      B'\!\left( X_{ l h } \right) 
      \left(
        X_{ s } - X_{ l h }
        -
        \int_{ l h }^{ s }
        B\!\left( X_{ l h } \right) dW^K_u
      \right)
    \right\|_{ HS( U_0, H ) }^2 ds
\\&\leq
  R^2 \left(
    \sum_{ l=0 }^{ m-1 }
    \int_{ l h }^{ (l+1)h }
    \mathbb{E}
    \left\|
      X_{ s } - X_{ l h }
      -
      \int_{ l h }^{ s }
      B\!\left( X_{ l h } \right) dW^K_u
    \right\|_{ H }^2 ds
  \right)
\end{align*}
holds 
for all $ m \in \{0,1,\ldots,M\} $ and 
all $ M, K \in \mathbb{N} $.
Hence, \eqref{eq:estimate_1} yields
\begin{align*}
& \mathbb{E}\left\|
    \sum_{ l=0 }^{ m-1 } 
    \int_{ l h }^{ (l+1)h }
    e^{ A(m-l)h }
    B'\!\left( X_{ l h } \right) 
    \left( 
      X_{ s } - X_{ l h }
      -
      \int_{ l h }^{ s }
      B\!\left( X_{ l h } \right) dW^K_u 
      \right) dW^K_s
  \right\|_{ H }^2
\\&\leq
  20 R^{ 10 } \left(
    \sum_{ l=0 }^{ m-1 }
    \int_{ l h }^{ (l+1)h }
      \left(
  \left( s - l h \right)^{ 
    \min\left( 
      4 ( \gamma - \beta ) , 
      2 \gamma
    \right) }
  +
  \bigg(
    \sup_{ j \in \mathcal{J} \backslash \mathcal{J}_K }
    \eta_j
  \bigg)^{ \! 2 \alpha }
      \right)
    ds
  \right)
\\&\leq
  20 R^{ 10 }
  \left(
    \sum_{ l=0 }^{ m-1 }
    \int_{ l h }^{ (l+1)h }
    \left( s - l h 
    \right)^{ 
      \min\left( 
        4 ( \gamma - \beta ) , 
        2 \gamma
      \right)
    } ds
  \right)
  +
  20 R^{ 10 }
    T \bigg(
    \sup_{ 
      j \in \mathcal{J} 
      \backslash \mathcal{J}_K 
    }
    \eta_j
    \bigg)^{ \! 2 \alpha }
\\&\leq
  20 R^{ 10 }
    M h^{ \left( 1 + 
      \min\left( 4 \left( \gamma -
        \beta \right), 
        2 \gamma 
      \right) 
    \right) 
  }
  +
  20 R^{ 11 }
    \bigg(
    \sup_{ 
      j \in \mathcal{J} 
      \backslash \mathcal{J}_K 
    }
    \eta_j
    \bigg)^{ \! 2 \alpha }
\\&\leq
  20 R^{ 11 }
    h^{ \min\left( 
      4 \left( \gamma -
      \beta \right), 2 \gamma
    \right) }
  +
  20 R^{ 11 }
    \bigg(
    \sup_{ j \in \mathcal{J} 
      \backslash \mathcal{J}_K 
    }
    \eta_j
    \bigg)^{ \! 2 \alpha }
\leq
  \frac{ 20 R^{ 13 } }{
    M^{ 
      \min\left( 
        4 \left( \gamma -
        \beta \right), 
        2 \gamma 
      \right) 
    }
  }
  +
  20 R^{ 11 }
    \bigg(
    \sup_{ 
      j \in \mathcal{J} 
      \backslash \mathcal{J}_K 
    }
    \eta_j
    \bigg)^{ \! 2 \alpha }
\end{align*}
for all $ m \in \{0,1,\ldots,M\} $ and 
all $ M, K \in \mathbb{N} $.
\subsection{Lipschitz estimates: 
Proof of \eqref{eq:lipproof}}
\label{sec:lipproof}
Before we estimate 
$ \mathbb{E}\left\| Z_m^{ N, M, K } - 
Y_m^{ N, M, K } \right\|_{ H }^2 $ 
for $ m \in \{0,1,\ldots,M\} $ 
and for $ N, M, K \in \mathbb{N} $, 
we need
some preparations.
More precisely, we have
\begin{align*}
&
  \mathbb{E}\left\|
    B'\!\left( X_{ l h } \right)
      \int_{ l h }^s 
      B\!\left( X_{ l h } \right) dW_u^K
     - 
    B'\!\left( Y_l^{ N,M,K } \right)
      \int_{ l h }^s 
      B\!\left( Y_l^{ N,M,K } \right) dW_u^K
  \right\|_{ HS( U_0, H ) }^2
\\&=
  \mathbb{E}\Bigg\|
    B'\!\left( X_{ l h } \right)
    \left( 
      \sum_{ \substack{ j \in \mathcal{J}_K \\ \eta_j \neq 0 } }
      \int_{ l h }^s
      B\!\left( X_{ l h } \right)
      g_j \, d\!\left< g_j, W_u \right>_U
    \right)
\\&\quad
    -
    B'\!\left( Y_l^{ N,M,K } \right)
    \left( 
      \sum_{ \substack{ j \in \mathcal{J}_K  \\ \eta_j \neq 0 } }
      \int_{ l h }^s
      B\!\left( Y_l^{ N,M,K } \right)
      g_j \, d\!\left< g_j, W_u \right>_U
    \right)
  \Bigg\|_{ HS( U_0, H ) }^2
\\&=
  \mathbb{E}\Bigg\|
    B'\!\left( X_{ l h } \right)
    \left( 
      \sum_{ \substack{ j \in \mathcal{J}_K  \\ \eta_j \neq 0 } }
      B\!\left( X_{ l h } \right)
      g_j 
      \left< g_j, W_s - W_{ l h } \right>_U
    \right)
\\&\quad
    -
    B'\!\left( Y_l^{ N,M,K } \right)
    \left( 
      \sum_{ \substack{ j \in \mathcal{J}_K \\ \eta_j \neq 0 } }
      B\!\left( Y_l^{ N,M,K } \right)
      g_j 
      \left< g_j, W_s - W_{ l h } \right>_U
    \right)
  \Bigg\|_{ HS( U_0, H ) }^2
\\&=
  \mathbb{E}\Bigg\|
    \sum_{ \substack{ j \in \mathcal{J}_K \\ \eta_j \neq 0 } }
    \bigg\{
      B'\!\left( X_{ l h } \right)
      \left( 
        B\!\left( X_{ l h } \right)
        g_j
      \right)
      -
      B'\!\left( Y_l^{ N,M,K } \right)
      \left( 
        B\!\left( Y_l^{ N,M,K } \right)
        g_j
      \right)
    \bigg\}
    \left< g_j, W_s - W_{ l h } \right>_U
  \Bigg\|_{ HS( U_0, H ) }^2
\end{align*}
for all $ s \in [l h, (l+1)h ] $, $ l \in \{0,1,\ldots,M-1\} $ and all $ M, K \in \mathbb{N} $.
This implies
\begin{align*}
&
  \mathbb{E}\left\|
    B'\!\left( X_{ l h } \right)
      \int_{ l h }^s 
      B\!\left( X_{ l h } \right) dW_u^K
     - 
    B'\!\left( Y_l^{ N,M,K } \right)
      \int_{ l h }^s 
      B\!\left( Y_l^{ N,M,K } \right) dW_u^K
  \right\|_{ HS( U_0, H ) }^2
\\&=
  \sum_{ \substack{ j \in \mathcal{J}_K \\ \eta_j \neq 0 } }
  \mathbb{E}\Bigg\|
    \bigg\{
      B'\!\left( X_{ l h } \right)
      \left( 
        B\!\left( X_{ l h } \right)
        g_j
      \right)
      -
      B'\!\left( Y_l^{ N,M,K } \right)
      \left( 
        B\!\left( Y_l^{ N,M,K } \right)
        g_j
      \right)
    \bigg\}
    \left< g_j, W_s - W_{ l h } \right>_U
  \Bigg\|_{ HS( U_0, H ) }^2
\\&=
  \sum_{ \substack{ j \in \mathcal{J}_K \\ \eta_j \neq 0 } }
  \mathbb{E}\left\|
    B'\!\left( X_{ l h } \right)
    \left( 
      B\!\left( X_{ l h } \right)
      g_j
    \right)
    -
    B'\!\left( Y_l^{ N,M,K } \right)
    \left( 
      B\!\left( Y_l^{ N,M,K } \right)
      g_j
    \right)
  \right\|_{ HS( U_0, H ) }^2
  \!\!\!\cdot
  \mathbb{E}\left|
    \left< g_j, W_s - W_{ l h } \right>_U
  \right|^2
\end{align*}
for all $ s \in [l h, (l+1)h ] $, $ l \in \{0,1,\ldots,M-1\} $ and all $ M, K \in \mathbb{N} $.
Hence, we obtain
\begin{align}
\label{eq:temporal_proof}
&
  \mathbb{E}\left\|
    B'\!\left( X_{ l h } \right)
      \int_{ l h }^s 
      B\!\left( X_{ l h } \right) dW_u^K
     - 
    B'\!\left( Y_l^{ N,M,K } \right)
      \int_{ l h }^s 
      B\!\left( Y_l^{ N,M,K } \right) dW_u^K
  \right\|_{ HS( U_0, H ) }^2
\nonumber
\\&=
\nonumber
    \sum_{ \substack{ j \in \mathcal{J}_K \\ \eta_j \neq 0 } } \eta_j \cdot
    \mathbb{E}\left\|
      B'\!\left( X_{ l h } \right)\!
      \left( 
        B\!\left( X_{ l h } \right)
        g_j
      \right)
      -
      B'\!\left( Y_l^{ N,M,K } \right)\!
      \left( 
        B\!\left( Y_l^{ N,M,K } \right)
        g_j
      \right)
    \right\|_{ HS( U_0, H ) }^2
  \left( s - l h \right)
\\&\leq
    \!\!
    \sum_{ \substack{ j, k \in \mathcal{J} 
      \\ \eta_j, \eta_k \neq 0 } } 
    \!\!
    \eta_j \, \eta_k \,
    \mathbb{E}\left\|
      B'\!\left( X_{ l h } \right) \!
      \left( 
        B\!\left( X_{ l h } \right)
        g_j
      \right) g_k
      -
      B'\!\left( Y_l^{ N,M,K } \right) \!\!
      \left( 
        B\!\left( Y_l^{ N,M,K } \right)
        g_j
      \right) g_k
    \right\|_{ H }^2
  \left( s - l h \right)
\\&=
  \left( s - l h \right)
  \cdot
  \mathbb{E}\left\|
    B'\!\left( X_{ l h } \right)
    B\!\left( X_{ l h } \right)
    -
    B'\!\left( Y_l^{ N,M,K } \right) 
    B\!\left( Y_l^{ N,M,K } \right)
  \right\|_{ HS^{(2)}( U_0, H ) }^2
\nonumber
\\&\leq
\nonumber
  c^2 \cdot \left( s - l h \right)
  \cdot
  \mathbb{E}\left\|
    X_{ l h }
    -
    Y_l^{ N,M,K }
  \right\|_{ H }^2
\leq
  R^3
  \cdot
  \mathbb{E}\left\|
    X_{ l h }
    -
    Y_l^{ N,M,K }
  \right\|_{ H }^2
\end{align}
for all $ s \in [l h, (l+1)h ] $, 
$ l \in \{0,1,\ldots,M-1\} $ and 
all $ M, K \in \mathbb{N} $. 
Additionally, 
\eqref{eq:squareest} 
and the fact
$ \| P_N(v) \|_H
\leq \| v \|_H $
for all $ v \in H $
show
\begin{equation}
\begin{split}
&
  \mathbb{E}\left\|
    Z_m^{ N, M, K }
    -
    Y_m^{ N, M, K }
  \right\|_{ H }^2
\\ &\leq
  3 \cdot \mathbb{E}\left\|
    \sum_{ l=0 }^{ m - 1 }
    \int_{ l h }^{ (l+1)h }
    e^{ A(m-l)h }
    \left(
      F\!\left( X_{ l h } \right)
      -
      F\!\left( Y_l^{ N,M,K } \right)
    \right) ds
  \right\|_{ H }^2
\\&\quad+
 3 \cdot \mathbb{E}\left\|
    \sum_{ l=0 }^{ m - 1 }
    \int_{ l h }^{ (l+1)h }
    e^{ A(m-l)h }
    \left(
      B\!\left( X_{ l h } \right)
      -
      B\!\left( Y_l^{ N,M,K } \right)
    \right) dW_s^K
  \right\|_{ H }^2
\\&\quad+
  3 \cdot \mathbb{E}\Bigg\|
    \sum_{ l=0 }^{ m - 1 }
    \int_{ l h }^{ (l+1)h }
    e^{ A(m-l)h }
    \bigg( 
      B'\!\left( X_{ l h } \right)
      \int_ { l h }^{ s }
      B\!\left( X_{ l h } \right) dW_u^K
      -
      B'\!\left( Y_l^{ N,M,K } \right)
      \int_ { l h }^{ s }
      B\!\left( Y_l^{ N,M,K } \right) dW_u^K
    \bigg) dW_s^K
  \Bigg\|_{ H }^2
\\&\leq
  3 M h^2 \left(
    \sum_{ l=0 }^{ m - 1 }
    \mathbb{E}\left\|
      e^{ A(m-l)h }
      \left(
        F\!\left( X_{ l h } \right)
        -
        F\!\left( Y_l^{ N,M,K } \right)
      \right)
    \right\|_{ H }^2
  \right)
\\&\quad+
 3 \left(
    \sum_{ l=0 }^{ m - 1 }
    \int_{ l h }^{ (l+1)h }
    \mathbb{E}\left\|
      e^{ A(m-l)h }
      \left(
        B\!\left( X_{ l h } \right)
        -
        B\!\left( Y_l^{ N,M,K } \right)
      \right)
    \right\|_{ HS( U_0, H ) }^2 ds
  \right)
\\&\quad+
  3 \Bigg(
    \sum_{ l=0 }^{ m - 1 }
    \int_{ l h }^{ (l+1)h }
    \mathbb{E}\bigg\|
      B'\!\left( X_{ l h }  \right)
      \int_ { l h }^{ s }
      B\!\left( X_{ l h }  \right) dW_u^K
      -
      B'\!\left( Y_l^{ N,M,K } \right)
      \int_ { l h }^{ s }
      B\!\left( Y_l^{ N,M,K } \right) dW_u^K
    \bigg\|_{ HS(U_0, H) }^2 ds
  \Bigg)
\end{split}
\end{equation}
and due to \eqref{eq:temporal_proof} 
we finally obtain
\begin{align*}
& \mathbb{E}\left\|
    Z_m^{ N, M, K }
    -
    Y_m^{ N, M, K }
  \right\|_{ H }^2
\leq
  3 T h \left(
    \sum_{ l=0 }^{ m - 1 }
    \mathbb{E}\left\|
      F\!\left( X_{ l h } \right)
      -
      F\!\left( Y_l^{ N,M,K } \right)
    \right\|_{ H }^2
  \right)
\\&\quad+
 3 h \left(
    \sum_{ l=0 }^{ m - 1 }
    \mathbb{E}\left\|
      B\!\left( X_{ l h } \right)
      -
      B\!\left( Y_l^{ N,M,K } \right)
    \right\|_{ HS( U_0, H ) }^2
  \right)
  +
  3 R^3 h \left(
    \sum_{ l=0 }^{ m - 1 }
    \mathbb{E}\left\|
      X_{ l h }
      -
      Y_l^{ N,M,K }
    \right\|_{ H }^2 
  \right)
\\&\leq
  9 R^3 h \left(
    \sum_{ l=0 }^{ m - 1 }
    \mathbb{E}\left\|
      X_{ l h }
      -
      Y_l^{ N,M,K }
    \right\|_H^2 
  \right)
  \leq
  \frac{ 9 R^4 }{ M } \left(
    \sum_{ l=0 }^{ m - 1 }
    \mathbb{E}\left\|
      X_{ l h }
      -
      Y_l^{ N,M,K }
    \right\|_H^2 
  \right)
\end{align*}
for all $ m \in \{0,1,\ldots,M\} $ and 
all $ N,M,K \in \mathbb{N} $.
\subsection{Iterated integral identity:
Proof of \eqref{eq:centralob}}
\label{sec:centralob}
First of all, we have
\begin{align*}
 &\int_{\frac{m T}{M}}^{\frac{(m+1)T}{M}}
  B'\!\left( Y_m^{N,M,K} \right) \!
  \bigg( 
    \int_{\frac{m T}{M}}^s
    B\!\left( Y_m^{N,M,K} \right) dW_u^K 
  \bigg) dW_s^K
\\&=
  \sum_{ \substack{j,k \in \mathcal{J}_K \\ \eta_j, \eta_k \neq 0 }  }
  \int_{\frac{m T}{M}}^{\frac{(m+1)T}{M}}
  B'\!\left( Y_m^{N,M,K} \right) \!
  \bigg( 
    \int_{\frac{m T}{M}}^s
    B\!\left( Y_m^{N,M,K} \right) 
    g_k \, d\!\left<g_k, W_u\right>_U \!\!
  \bigg) g_j \, d\!\left<g_j, W_s\right>_U
\\&=
  \sum_{ \substack{ j,k \in \mathcal{J}_K  
    \\ \eta_j, \eta_k \neq 0 }  } \!\!
  B'\!\left( Y_m^{N,M,K} \right) \!
  \bigg( 
    B\!\left( Y_m^{N,M,K} \right) 
    g_k 
  \bigg) \, g_j
  \cdot
    \int_{\frac{m T}{M}}^{\frac{(m+1)T}{M}} \!\!\!
    \int_{\frac{m T}{M}}^s 
    d\!\left<g_k, W_u\right>_U d\!\left<g_j, W_s\right>_U
\end{align*}
$\mathbb{P}$-a.s.\ and 
hence
\begin{align*}
 &\int_{\frac{m T}{M}}^{\frac{(m+1)T}{M}}
  B'\!\left( Y_m^{N,M,K} \right) \!
  \bigg( 
    \int_{\frac{m T}{M}}^s
    B\!\left( Y_m^{N,M,K} \right) dW_u^K 
  \bigg) dW_s^K
\\&=
  \sum_{ \substack{ j \in \mathcal{J}_K \\ \eta_j \neq 0 }  }
  B'\!\left( Y_m^{N,M,K} \right) \!
  \bigg( 
    B\!\left( Y_m^{N,M,K} \right) 
    g_j 
  \bigg) \, g_j
  \cdot
    \int_{\frac{m T}{M}}^{\frac{(m+1)T}{M}} \!\!\!
    \int_{\frac{m T}{M}}^s 
    d\!\left<g_j, W_u\right>_U d\!\left<g_j, W_s\right>_U
\\&+
  \sum_{ 
    \substack{ j,k \in \mathcal{J}_K 
      \\ \eta_j, \eta_k \neq 0 \\
      j \neq k
    }
  } \!
  B'\!\left( Y_m^{N,M,K} \right) \!
  \bigg( 
    B\!\left( Y_m^{N,M,K} \right) 
    g_k 
  \bigg) \, g_j
  \cdot
    \int_{\frac{m T}{M}}^{\frac{(m+1)T}{M}} \!\!\!
    \int_{\frac{m T}{M}}^s 
    d\!\left<g_k, W_u\right>_U d\!\left<g_j, W_s\right>_U
\end{align*}
$\mathbb{P}$-a.s.\ for 
all $ m \in \{0,1,\ldots, M-1 \} $ 
and all $ N,M,K \in \mathbb{N} $.
Moreover, since the bilinear 
operator 
$ B'\!\left( Y_m^{N,M,K} \right) B\!\left( Y_m^{N,M,K} \right) 
\in HS^{(2)}( U_0, H ) $ is symmetric 
(see Assumption~\ref{diffusion})
and since 
\[
  \int_{\frac{m T}{M}}^{\frac{(m+1)T}{M}}
  \int_{\frac{m T}{M}}^s
  d\!\left<g_j, W_u\right>_U \, d\!\left<g_j, W_s\right>_U
  =
  \frac{1}{2}
  \left( \!
    \left(
      \left<g_j, \Delta W_m^{M,K}\right>_U
    \right)^2
    -
    \frac{ T \eta_j }{ M }
  \right)
\]
$\mathbb{P}$-a.s.\ holds for 
all $ j \in \mathcal{J}_K $, 
$ m \in \{0,1,\ldots,M-1\} $ 
and all $ N,M,K \in \mathbb{N} $ 
(see (3.6) in Section~10.3 in \cite{kp92}),
we obtain
\begin{align*}
 &\int_{\frac{m T}{M}}^{\frac{(m+1)T}{M}}
  B'\!\left( Y_m^{N,M,K} \right) \!
  \bigg( 
    \int_{\frac{m T}{M}}^s
    B\!\left( Y_m^{N,M,K} \right) dW_u^K 
  \bigg) dW_s^K
\\&=
  \frac{1}{2}
  \sum_{ \substack{ j \in \mathcal{J}_K \\ \eta_j \neq 0 }  }
  B'\!\left( Y_m^{N,M,K} \right) \!
  \bigg( 
    B\!\left( Y_m^{N,M,K} \right)
    g_j
  \bigg) \, g_j
  \left( 
    \left(
      \left<g_j, \Delta W_m^{M,K}\right>_U
    \right)^2
    -
    \frac{ T \eta_j }{ M }
  \right)
\\&\quad+
  \frac{1}{2}
  \sum_{ \substack{ j,k \in \mathcal{J}_K 
    \\ \eta_j, \eta_k \neq 0 
    \\ j \neq k }
  }
  B'\!\left( Y_m^{N,M,K} \right) \!
  \bigg( 
    B\!\left( Y_m^{N,M,K} \right)
    g_k
  \bigg) \, g_j 
\\&\qquad\quad  \cdot
  \left(
      \int_{\frac{m T}{M}}^{\frac{(m+1)T}{M}} \!\!\!
      \int_{\frac{m T}{M}}^s 
      d\!\left<g_k, W_u\right>_U d\!\left<g_j, W_s\right>_U
    +
      \int_{\frac{m T}{M}}^{\frac{(m+1)T}{M}} \!\!\!
      \int_{\frac{m T}{M}}^s 
      d\!\left<g_j, W_u\right>_U d\!\left<g_k, W_s\right>_U
  \right)
\end{align*}
$\mathbb{P}$-a.s.\ for all $ m \in \{0,1,\ldots, M-1 \} $ 
and all $ N,M,K \in \mathbb{N} $.
The fact 
\begin{equation}
      \int_{\frac{m T}{M}}^{\frac{(m+1)T}{M}} \!\!\!
      \int_{\frac{m T}{M}}^s 
      d\!\left<g_k, W_u\right>_U d\!\left<g_j, W_s\right>_U
    +
      \int_{\frac{m T}{M}}^{\frac{(m+1)T}{M}} \!\!\!
      \int_{\frac{m T}{M}}^s 
      d\!\left<g_j, W_u\right>_U d\!\left<g_k, W_s\right>_U
  =
  \left<g_k, \Delta W_m^{M,K}\right>_U
  \left<g_j, \Delta W_m^{M,K}\right>_U
\end{equation}
$ \mathbb{P} $-a.s.\ for 
all $ j \in \mathcal{J}_K $, $ m \in \{0,1,\ldots,M-1\} $ 
and all $ N,M,K \in \mathbb{N} $ 
(see (3.15) in Section~10.3 in \cite{kp92}) 
then yields
\begin{align*}
 &\int_{\frac{m T}{M}}^{\frac{(m+1)T}{M}}
  B'\!\left( Y_m^{N,M,K} \right) \!
  \bigg( 
    \int_{\frac{m T}{M}}^s
    B\!\left( Y_m^{N,M,K} \right) dW_u^K 
  \bigg) dW_s^K
\\&=
  \frac{1}{2}
  \sum_{ \substack{ j, k \in \mathcal{J}_K \\ \eta_j, \eta_k \neq 0 }  }
  B'\!\left( Y_m^{N,M,K} \right) \!
  \bigg( \!
    B\!\left( Y_m^{N,M,K} \right)
    g_k
  \bigg) \, g_j
  \left<g_k, \Delta W_m^{M,K}\right>_U
  \left<g_j, \Delta W_m^{M,K}\right>_U
\\&\quad-
  \frac{T}{2M}
  \sum_{ \substack{ j \in \mathcal{J}_K \\ \eta_j \neq 0 }  }
  \eta_j B'\!\left( Y_m^{N,M,K} \right) \!
  \bigg( \!
    B\!\left( Y_m^{N,M,K} \right)
    g_j
  \bigg) \, g_j
\\&=
  \frac{1}{2}
  B'\!\left( Y_m^{N,M,K} \right) \!
  \bigg( \!
    B\!\left( Y_m^{N,M,K} \right)
    \Delta W_m^{ M,K }
  \bigg) \, \Delta W_m^{ M,K }
-
  \frac{T}{2M}
  \sum_{ \substack{ j \in \mathcal{J}_K \\ \eta_j \neq 0 }  }
  \eta_j
  B'\!\left( Y_m^{N,M,K} \right) \!
  \bigg( \!
    B\!\left( Y_m^{N,M,K} \right)
    g_j
  \bigg) \, g_j
\end{align*}
$ \mathbb{P} $-a.s.\ for 
all $ m \in \{0,1,\ldots, M-1 \} $ 
and all $ N,M,K \in \mathbb{N} $ 
which shows \eqref{eq:centralob}.
\subsubsection*{Acknowledgment}
This work has been 
supported by the 
BiBoS Research Center,
by the research project
``Numerical solutions of
stochastic differential equations
with non-globally Lipschitz
continuous coefficients''
funded by the German Research
Foundation,
by the Collaborative 
Research Centre $701$
``Spectral Structures 
and Topological Methods 
in Mathematics'' 
funded by the German
Research Foundation
and
by the
International Graduate 
School ``Stochastics and Real 
World Models''
funded by the German
Research Foundation.
The support of Issac Newton
Institute for Mathematical
Sciences in Cambridge is
also gratefully acknowledged
where part of this was done
during the special semester
on ``Stochastic Partial Differential
Equations''.
In addition, we are very 
grateful 
to Sebastian Becker for his
help with the 
numerical simulations
and to Carlo Marinelli for his
helpful advice concerning
Schatten norms.
\bibliographystyle{acm}
\bibliography{../Bib/bibfile}
\end{document}